\documentclass[11pt, reqno]{amsart}

\usepackage{amsmath, amsthm, amsfonts}
\usepackage{epsfig}
\newtheorem{prop}{Proposition}
\newtheorem{theorem}{Theorem}
\newtheorem{coroll}{Corollary}
\newtheorem{lemma}{Lemma}

\newtheorem{remark}{Remark}

\def\floor#1{\lfloor #1 \rfloor}
\def\ceil#1{\lceil #1 \rceil}
\newcommand{\A}{\mathbb{A}}
\newcommand{\B}{\mathbb{B}}
\newcommand{\C}{\mathbb{C}}
\newcommand{\Z}{\mathbb{Z}}

\newcommand{\un}{\mathbf{1}}
\renewcommand{\P}{\mathbb{P}}
\newcommand{\E}{\mathbb{E}}

\newcommand{\F}{\mathcal{F}}
\newcommand{\I}{\mathcal{I}}
\newcommand{\G}{\mathcal{G}}
\renewcommand{\H}{\mathcal{H}}
\renewcommand{\L}{\mathbb{L}}
\newcommand{\Q}{\mathbb{Q}}
\newcommand{\D}{\mathcal{D}}
\newcommand{\T}{\mathbb{T}}

\author{Jean B\'{e}rard$^{1}$ and Alejandro Ram\'{\i}rez$^{1,2}$}

\thanks{$^1$Partially supported by ECOS-Conicyt grant CO5EO2}

\thanks{$^2$Partially supported by Fondo Nacional de Desarrollo Cient\'\i fico
y Tecnol\'ogico grant 1060738 and by Iniciativa Cient\'\i fica Milenio P-04-069-F}

\address[Jean B\'{e}rard]{\noindent Universit\'e de Lyon ;
Universit\'e Lyon 1 ;
Institut Camille Jordan CNRS UMR 5208 ;
43, boulevard du 11 novembre 1918,
F-69622 Villeurbanne Cedex; France
\newline
e-mail:  \rm \texttt{jean.berard@univ-lyon1.fr}}

\address[Alejandro F. Ram\'\i rez]{Facultad de Matem\'aticas\\
Pontificia Universidad Cat\'olica de Chile\\
Vicu\~na Mackenna 4860, Macul\\
Santiago, Chile
\newline
e-mail:  \rm \texttt{aramirez@mat.puc.cl}}

\date{}

\title[Large deviations for a one dimensional model of $X+Y \to 2X$]{Large deviations of the front in a one dimensional model of $X+Y \to 2X$}

\begin{document}


\begin{abstract}  We investigate the probabilities of large deviations for the position of the front in a stochastic model of the reaction
$X+Y \to 2X$ on the integer lattice in which $Y$ particles do not move while $X$ particles move as independent simple continuous time
random walks of total jump rate $2$. For a wide class of initial conditions, we prove that a large deviations principle
holds and we show that the zero set of the rate function is the interval $[0,v]$, where
 $v$ is the velocity of the front given by the law of large numbers. 
 We also give more precise estimates for the rate of decay of the slowdown
probabilities. Our results indicate a gapless property of the generator of the process as seen from the front,
as it happens in the context of  nonlinear diffusion equations describing the propagation of a pulled
front into an unstable state.
 \end{abstract}



\maketitle

\section{Introduction}

We consider a microscopic model of a one-dimensional reaction-diffusion equation, with
a propagating front representing the passage from an unstable equilibrium to a stable one.
It is defined  as an interacting particle system on the integer lattice $\Z$ with two types of particles:
 $X$ particles, that move as independent, continuous time, symmetric, simple random walks with
total jump rate $D_X=2$; and $Y$ particles, which are inert and can be interpreted as
random walks with total jump rate $D_Y=0$. Initially, each site $x = 0,-1,-2,\ldots$ bears a certain number $ \eta(x) \geq 0$
of $X$ particles (with at least one site $x$ such that $\eta(x) \geq 1$), while each site $x=0,1,\ldots$ bears a fixed number $a$ of particles of type $Y$ (with $1 \leq a <+\infty$).
When a site $x = 1,2,\ldots$ is visited by an $X$ particle for the first time, all the 
$Y$ particles located at site $x$ are instantaneously turned into $X$ particles, and start moving. The {\it front} at time $t$ is defined as the rightmost site that has been visited by an $X$ particle
 up to time $t$, and is denoted by $r_{t}$, with the convention $r_{0}:=0$. This model can be interpreted as
an infection process, where  the $X$ and $Y$ particles represent ill and healthy individuals 
respectively.
It can also be interpreted as a combustion reaction, where the $X$ and $Y$ particles
correspond to heat units and reactive molecules respectively,  modeling the combustion of a propellant into
a stable stationary state. We will denote this model the {\it $X+Y\to 2X$ front propagation process} with jump rates $D_X$ and $D_Y$. Within the physics literature, a number of studies have been done both numerically and 
analytically of this process for different values of $D_X$ and $D_Y$ and of
corresponding variants where the infection of a $Y$ particle by an $X$
particle at the same site is not instantaneous,
 drawing analogies with continuous space
time  nonlinear reaction-diffusion equations having uniformly traveling wave solutions
\cite{Pan}, \cite{MaiSokBlu,MaiSokBlu2, MaiSokKuzBlu}, \cite{vSa}, 
\cite{EbevSa}. A particular well-known
 example is the  F-KPP equation studied by Fisher \cite{Fis} and Kolmogorov, Petrovsky and Piscounov \cite{KolPetPis}.

 Mathematically not too much is known. For the case $D_Y=0$,  when $\sum_{x \leq 0} \exp(\theta x) \eta(x)<+\infty$ for a small enough $\theta>0$,
a  law of large numbers with
a deterministic speed $0<v<+\infty$  not depending on the initial condition is
satisfied (see \cite{RamSid} and \cite{ComQuaRam2}):

\begin{equation}\label{e:lgn-presque-sure}\lim_{t\to +\infty} t^{-1} r_t=v  \hspace{2em} a.s.\end{equation}
  In~\cite{ComQuaRam2} it was proved  that the fluctuations around this speed satisfy 
a functional central limit theorem and that the marginal law of the particle
 configuration as seen from the front
 converges to a unique invariant measure as $t\to\infty$. Furthermore, a multi-dimensional version of
this process on the lattice $\mathbb Z^d$, with an initial configuration having one $X$ particle at the origin and
one $Y$ particle at every other site was studied in \cite{RamSid}, \cite{AlvMacPop}, proving an asymptotic
shape theorem as $t\to\infty$ for the set of visited sites. A similar result was proved by Kesten and
Sidoravicius \cite{KesSid}
for the case $D_X=D_Y>0$  with a product Poisson initial law. In particular, in dimension $d=1$ they prove
a law of large numbers for the front as in (\ref{e:lgn-presque-sure}). For the case $D_X> D_Y>0$, even the problem of proving 
a law of large numbers  in dimension $d=1$ remains open (see \cite{KesSid2}).


 Within a certain class
of one-dimensional nonlinear diffusion equations having uniformly traveling wave solutions describing the passage
from an unstable to a stable state,  it has been observed that for
certain initial conditions
the velocity of the front at a given time has a rate
of relaxation towards its asymptotic value which 
 is algebraic (see \cite{EbevSa}, \cite{Pan} and physics literature references therein). 
These are the so called {\it pulled} fronts, whose speed is determined
by a region of the profile linearized about the unstable solution.
For the F-KPP equation, Bramson \cite{Bra} proved   that the speed of the front at a given time is below
its asymptotic value and that the convergence is algebraic. 
In general,  the slow relaxation is due to a gapless property of a linear operator governing
the convergence of the centered front profile towards the stationary state. 
A natural question is wether such a behavior can be observed in the 
 $X+Y\to 2X$ front propagation type processes. Deviations from the
law of large numbers of a larger size than those given by central limit theorem should shed some light
on such a question: in particular it would be
reasonable to expect a large deviations principle with a degenerate rate function, reflecting a
slow convergence of the particle configuration as seen from equilibrium towards the unique
invariant measure \cite{ComQuaRam2}.
In this paper, we investigate for
the case $D_Y=0$ the
large time asymptotics of the distribution of $r_t/t$,
$$
\P\left[\frac{r_t}{t}\in \cdot\right].
$$
Our main result is that a full large deviations principle holds, with a degenerate rate function on
the interval $[0,v]$, when the initial condition
 satisfies the following  growth condition:

\noindent {\bf Assumption (G)}. For all $\theta>0$
 \begin{equation}\label{e:decroiss-petit-exp}
\sum_{x \leq 0} \exp(\theta x) \eta(x)<+\infty.\end{equation}

\bigskip

\begin{theorem}\label{t:ldp} {\bf Large Deviations Principle} There exists a rate function
$I \, : \, [0,+\infty) \to [0,+\infty)$ such that, for every initial condition satisfying {\bf (G)},

$$
\limsup_{t\to +\infty}\frac{1}{t}\log \P\left[\frac{r_t}{t}\in C\right]\le -\inf_{b\in C}I(b),\qquad
{\rm for}\quad C\subset [0,+\infty)\quad{\rm closed},
$$
and
$$
\liminf_{t\to +\infty}\frac{1}{t}\log \P\left[\frac{r_t}{t}\in G\right]\ge -\inf_{b\in G}I(b),\qquad
{\rm for}\quad G\subset [0,+\infty)  \quad{\rm open}.
$$
Furthermore, $I$ is identically zero on $[0,v]$, positive, convex and increasing
on $(v,+\infty)$.
\end{theorem}
It is interesting to notice that the rate function $I$ is independent of the initial conditions within
the class {\bf (G)}:  
the large deviations of the empirical distribution function of the process as 
seen from the front  appear to exhibit a
uniform behavior for such initial conditions.
Furthermore, this result seems to be  in agreement with the phenomenon of slow  relaxation  of the velocity in the so-called pulled reaction diffusion equations.
In \cite{EbevSa}, a nonlinear diffusion equation of the form

\begin{equation}
\label{oo}
\partial_t\phi=\partial^2_x\phi+f(\phi)
\end{equation}
is studied
where $f$ is a  function  chosen so that  $\phi=0$ is an unstable state and the equation develops  pulled fronts. It is argued that for steep enough initial conditions, the velocity  relaxes algebraically towards the asymptotic speed, providing an
explicit expansion up to order $O(1/t^2)$.  Such a
non-exponential decay is explained by the fact that the
linearization of  (\ref{oo}) 
around the uniformly translating front,  gives a linear equation for the perturbation governed by a gapless Schr\"odinger
operator.  The position of the front in the $X+Y\to 2X$ particle system
can be decomposed as  $r_t=\int_0^t L g(\eta_s)ds+M_t$, where $L$ is the generator of the 
centered dynamics, $g$ is an explicit function and $M_t$ is a martingale. 
The fact that under assumption {\bf (G)} the zero set of the large deviations
 principle of Theorem 1 is the interval $[0,v]$  is an indication that the 
symmetrization of  $L$
is a gapless operator.

The second result of this paper
gives more precise estimates for the probability of the slowdown deviations.
Let

$$U(\eta):= \limsup_{x \to -\infty}   \frac{1}{\log |x|} \log \left( \sum_{y=0}^{x} \eta(y) \right), \, u(\eta):= \liminf_{x \to -\infty}   \frac{1}{\log |x|} \log \left( \sum_{y=0}^{x} \eta(y) \right),$$
and
$$s(\eta):=\min(1, U(\eta)).$$
For the statement of the following theorem we will write $U,u,s$ instead
of $U(\eta), u(\eta), s(\eta)$.
\medskip
\begin{theorem}\label{t:slowdown} {\bf Slowdown deviations estimates.}
Let $\eta$ be an initial condition satisfying~{\bf (G)}. Then the
following statements are satisfied.
\begin{itemize}
\item[(a)] For all $0 \leq c < b < v$, as $t$ goes to infinity,
\begin{equation}\label{e:slowdown-lower-bound}
\P \left[c \leq  \frac{ r_{t}}{t}  \leq b  \right] \geq  \exp \left( -t^{ s/2  +o(1)}     \right). \end{equation}

\item[(b)] In the special case where $\eta(x)\geq a$ for all $x \leq 0$, one has that,
for every $0 \leq b < v$, as $t$ goes to infinity,
\begin{equation}
\label{e:slowdown-upper-bound} \P\left[\frac{r_t}{t}\leq b\right] \leq \exp \left( -t^{1/3+o(1)}     \right). 
\end{equation}

\item[(c)] When $u<+\infty$, as $t$ goes to infinity, 
\begin{equation}\label{e:slowdown-zero-bound}
  \exp \left( -t^{U/2+o(1)} \right)   \leq  \P \left[ r_{t} = 0  \right]\leq \exp \left( -t^{u/2+o(1)} \right) .
 \end{equation}

\end{itemize}

\end{theorem}
One may notice that the slowdown probabilities considered in ~(\ref{e:slowdown-lower-bound}) and in (~\ref{e:slowdown-zero-bound}) exhibit distinct behaviors 
when $u>1$. 
Furthermore, the results contained in Theorems~\ref{t:ldp} and~\ref{t:slowdown} should be compared with the case
 of the random walk in random environment with positive or zero drift \cite{PisPovZei,PisPov}.

A natural question is whether it is possible to relax assumption~{\bf (G)} in Theorem~\ref{t:ldp}.
It appears that even if assumption~{\bf (G)} is but mildly violated, the slowdown behavior is not in accordance with that described
 by Theorem~\ref{t:ldp}. Moreover, if assumption~{\bf (G)} is  strongly violated,
  the law of large numbers  with asymptotic velocity $v$ breaks down, so that the speedup part of Theorem~\ref{t:ldp} cannot hold either.

\begin{theorem}\label{t:optimal}
The following properties hold:
\begin{itemize}
\item[(i)] Assume there is a $\theta>0$ such that
 $$\liminf_{x \to -\infty}\eta(x) \exp(\theta x) = +\infty.$$
 Then there exists $b>0$ such that
$$\limsup_{t \to +\infty} \frac{1}{t} \log \P\left[\frac{r_t}{t}\leq b \right]  < 0.$$
\item[(ii)] There exists $\theta'>0$ and $v'>v$ such that, when
 $$\liminf_{x \to -\infty}\eta(x) \exp(\theta' x) = +\infty,$$
 then
  $$ \P\left[\liminf_{t \to +\infty} \frac{r_t}{t} \geq v'\right] =1.$$
\end{itemize}
\end{theorem}

It is important to stress that the proof of Theorem~\ref{t:ldp} would not be much simplified if we 
considered initial conditions with only a finite number of particles. Indeed, condition {\bf (G)}
is an assumption which delimits sensible initial data.
To prove Theorem~\ref{t:ldp} we first establish that for  initial conditions consisting only of a
single particle at the origin, 
 for all $b \geq 0$, the limit \begin{equation}\label{e:lim-subadd-posassoc} \lim_{t \to +\infty} t^{-1} \log \P( r_{t} \geq bt )\end{equation}
 exists.
The proof of this fact relies on a soft argument based
on the sub-additivity property of the hitting times. 
On the other hand, it is not difficult to show that for $b$ large enough the decay of $\P(r_t\geq bt)$ is
exponentially fast. Nevertheless, showing this for $b$ arbitrarily close to but larger than the speed $v$ is a subtler problem.
For example, it is not clear how the standard sub-additive arguments could help.
Our main tool to tackle this problem is the
regeneration structure of the process defined in~\cite{ComQuaRam2}. To overcome the fact that the
regeneration times and positions
have only polynomial tails, we couple the original process with one where the $X$ particles have a small bias to the right,
so that they jump to the right with probability $1/2+\epsilon$ for some small $\epsilon>0$, and the position of the front in the biased process dominates that of the 
front in the original process. We then 
use the regeneration structure to study the biased model and how it relates to the original one
 as   $\epsilon$
 tends to zero.  In particular, if $v_\epsilon$ is the speed of the biased front, we  establish via uniform
bounds on the moments of the regeneration times and positions that
$$
\lim_{\epsilon\to 0^+}v_\epsilon=v.
$$
Furthermore, we show that the regeneration times and positions of the biased model have exponentially decaying tails.
Combining these arguments proves that the limit in~(\ref{e:lim-subadd-posassoc}) is positive for any $b>v$.
We then establish that this limit  exists and has the same value for all
 initial conditions satisfying~{\bf (G)} by exploiting a comparison argument.

To show that the rate function vanishes on $[0,v]$ (and more precisely~(\ref{e:slowdown-upper-bound})), we first consider initial conditions  having a uniformly
bounded number of particles per site. In this case it is essentially enough to observe that
the probability that the front remains at zero up to time $t$ is  bounded from
below by $(1/\sqrt{t})^{t^{1/2+o(1)}}$, since there
are at most of the order of $t^{1/2+o(1)}$ random walks that yield a non-negligible contribution to this event. Similar 
estimates on hitting times of random walks are used to prove~(\ref{e:slowdown-zero-bound}) and Theorem~\ref{t:optimal}, while more refined arguments are needed 
to establish~(\ref{e:slowdown-lower-bound}) for arbitrary initial conditions within the
class {\bf (G)}. On the other hand,
 the proof of the upper bound for the slowdown probabilities~(\ref{e:slowdown-upper-bound}) in Theorem~\ref{t:slowdown} is more involved, and relies on arguments
using the sub-additivity property and
the positive association of the
hitting times, together with estimates on their tails and  their correlations, refining
an idea already used in~\cite{RamSid} in  a similar context.

The rest of the paper is organized as follows. In Section~\ref{s:construction}, we give a formal definition of the model and introduce its basic structural properties, including
 sub-additivity and monotonicity of hitting times. In Section~\ref{s:proof-large-dev}, we explain how Theorem~\ref{t:ldp} is
proved, building on results proved in other sections. Section~\ref{s:speedup} is devoted to the proof of the fact that speedup
large deviations events have exponentially small
probabilities. Section~\ref{s:slowdown} contains our estimates on slowdown probabilities, with the proofs of Theorems~\ref{t:slowdown} and~\ref{t:optimal}. Several
 appendices contain proofs that are not included in the core of the paper.

\section{Construction and basic properties}\label{s:construction}

Throughout the sequel we will use the convention  $\inf \emptyset = +\infty$.
\subsection{Construction of the process}

For our purposes, we have to define on the same probability space not only the original
model, but also models including random walks with an arbitrary bias defined through a parameter $0 \leq \epsilon < 1/2$.

In the sequel, we assume that we have a reference probability space $(\Omega, \F, \P)$ giving us access
to an i.i.d. family of random variables
$$\left[ (\tau_{n}(u,i), U_{n}(u,i))  ; \, n \geq 1, \, u \in \Z , \, 1 \leq i \leq a           \right],$$
 such that,
 for all $(n,u,i)$,  $\tau_{n}(u,i)$ has an exponential(2) distribution, and $ U_{n}(u,i)$ has the uniform distribution on $[0,1]$,
 and $\tau_{n}(u,i)$ and $U_{n}(u,i)$ are independent.

For every $n \geq 1$, $(x,i) \in \Z \times \{1,\ldots, a \}$
and $0 \leq \epsilon < 1/2$, we let $$ \varepsilon_{n}(x,i,\epsilon):= 2 ( \un(  U_{n}(x,i)
 \leq 1/2+\epsilon )  )  - 1. $$
 Let $(Y^{\epsilon}_{x,i,t})_{t \geq 0}$ be the continuous-time random walk started at $Y^{\epsilon}_{x,i,0}:=0$,
whose sequence of time steps is
$(\tau_{n}(x,i))_{n \geq 1}$, and whose sequence of space increments is $(\varepsilon_{n}(x,i,\epsilon))_{n \geq 0}$.

A configuration of particles is a triple $w=(F,r,A)$, where $r \in \Z$, $A$ is a non-empty subset of $\Z \times \{1,\ldots, a \}$ such that $\max \{ x ; \, (x,i) \in A \} \leq r$, and
$F \ : \ A \to \{-\infty,\ldots, r \}$ is a map.
To every index $(x,i) \in A$ corresponds the position $F(x,i)$ of an $X$ particle, and we say that $(x,i)$ is the birthplace of the corresponding particle.
We see that such a configuration carries more information than just the number of $X$ particles at each site, since every $X$ particle is labeled by its birthplace $(x,i)$.
Note that, to allow for various types of initial configurations, we do not require that the initial configuration of the model satisfies $F(x,i)=(x,i)$. In fact,
any distribution of $X$ particles on $\{ \ldots, -1, 0 \}$ with a finite number of particles at each site can be encoded by such a triple $w=(F,r,A).$

For $w=(F,r,A)$ and $(x,i) \in A$, we use the notation $w(x,i)$ to 
denote the configuration $(F',r',A')$ with $r'=r$, $A'=\{ (x,i) \}$, and $F'(x,i)=F(x,i)$.
For a configuration of particles $w=(F,r,A)$ and $q \in \{1,2,\ldots, \}$, we define a configuration
$w \oplus q=(F',r',A')$ by $A':=A \cup \{r+1,\ldots,r+q\} \times \{1, \ldots, a \}$, $r':=r+q$, $F':=F$ on $A$, and
$F'(x,i):=x$ for $(x,i) \in  \{r+1,\ldots,r+q\} \times \{1, \ldots, a \}$.

The following definitions list special kinds of configurations that are used in the sequel.
For $u \in \Z$ and $1 \leq i \leq a$, let $\delta_{u}$  be defined by $A:=\{(u,1) \}$,  $r:=u$ and $F(u,1):=u$;
let $a \delta_{u}$ be defined by $A:= \{ u \} \times \{1,\ldots, a \}$, $r:=u$,
$F(u,i):=u$ for every $1 \leq i \leq a$;
let $\I_{u}$ be defined by
$A:= \{-\infty, u \} \times \{1,\ldots, a \}$, $r:=u$,
$F(x,i):=x$ for every $(x,i) \in A$.

For $w=(F,r,A)$ and $\theta>0$, we let
$$f_{\theta}(w):=\sum_{(x,i) \in A} \exp( \theta (F(x,i) - r)  ).$$
and $\eta_{w}$ be the map defined on $ \{\ldots,r-1, r \}$ so that
$\eta_w(x)$ is the number of particles at site $x$ of the configuration $w$. Hence
$$\eta_{w}(x) := \# \{   (y,i) \in A; \, F(y,i)=x   \}.$$
As a consequence,
$f_{\theta}(w)=\sum_{x \leq r} \eta_{w}(x) \exp( \theta (x -r) )$.
When $r=0$, we define for $x\leq 0$

\begin{equation}
\label{def-h}
H_{w}(x):=\sum_{y = 0}^{x} \eta_{w}(y).
\end{equation}
Now, for every $\theta>0$,  let
$$\L_{\theta} :=   \{   w=(F,r,A); \,    f_{\theta}(w)   < +\infty  \}.$$
Observe that $\I_{u}$, $\delta_{u}$ and $a \delta_{u}$ belong to $\L_{\theta}$ for all $u \in \Z$ and $\theta>0$.

For $w=(F,r,A)$ and $(x,i) \in \Z \times \{1, \ldots, a \}$, let $\chi_{\theta}(w,x,i) :=  \un((x,i) \in A)   \exp( \theta(F(x,i)-r) ) $.
We equip $\L_{\theta}$ with the metric $d_{\theta}$ defined as follows:  for $w=(F,r,A)$ and $w'=(F',r',A')$,
$$ d_{\theta}(w, w')  := | r-r' | +
 \sum_{(x,i) \in \Z\times \{1,\ldots, a \} }    |  \chi_{\theta}(w,x,i)  - \chi_{\theta}(w',x,i) |.$$
The metric space $(L_{\theta}, d_{\theta})$ is a Polish space.
We let $\D( \L_{\theta})$ denote the space of c\`adl\`ag functions from
$[0,+\infty)$ to $\L_{\theta}$ equipped with the Skorohod topology and the corresponding Borel $\sigma-$field.

Now, for every $0 \leq \epsilon < 1/2$, and every $w=(F,r,A) \in \L_{\theta}$,
we define a collection of random variables
$(X_{t}^{\epsilon}(w))_{t \geq 0}=( F_{t}^{\epsilon}(w), r_{t}^{\epsilon}(w), A_{t}^{\epsilon}(w))_{t \geq 0}$
which describes the time-evolution of the configuration of particles.
 In order to alleviate notations, the dependence of $F_{t}^{\epsilon}$, $r_{t}^{\epsilon}$, $A_{t}^{\epsilon}$
with respect to $w$ will not explicitly mentioned in the sequel when there is no ambiguity. Moreover, we
shall often not mention the dependence
with respect to $\epsilon$ when $\epsilon=0$, and for example, use the notation $r_{t}$ instead of $r^{0}_{t}$.

The definition is done through the following inductive procedure.
Let $\sigma^{\epsilon}_{0}:=0$, $r^{\epsilon}_{0}:=r$, $A^{\epsilon}_{0}:=A$,
and for every $t \geq 0$ and $(x,i) \in  A^{\epsilon}_{0}$, let $F^{\epsilon}_{t}(x,i):=F(x,i)+Y^{\epsilon}_{x,i,t}$.
Assume that, for some $n \geq 1$, we have already defined $\sigma^{\epsilon}_{0} \leq \ldots \leq \sigma^{\epsilon}_{n-1}$,
$A^{\epsilon}_{t}$ and $r^{\epsilon}_{t}$ for every $0 \leq t \leq \sigma^{\epsilon}_{n-1}$,
and $F^{\epsilon}_{t}(x,i)$ for every $t \geq 0$ and $(x,i) \in A^{\epsilon}_{\sigma^{\epsilon}_{n-1}}$.
Let  $$\sigma^{\epsilon}_{n} := \inf \left\{  t > \sigma^{\epsilon}_{n-1}; \,    \mbox{  there is an }
(x,i) \in A^{\epsilon}_{\sigma^{\epsilon}_{n-1}} \mbox{ such that }F^{\epsilon}_{t}(x,i)=r^{\epsilon}_{\sigma^{\epsilon}_{n-1}}+1 \right\},$$
 Now, for $\sigma^{\epsilon}_{n-1} < t < \sigma^{\epsilon}_{n}$, let
 $r^{\epsilon}_{t}:=r^{\epsilon}_{\sigma^{\epsilon}_{n-1}}$,
 $A^{\epsilon}_{t}:=A^{\epsilon}_{\sigma_{n-1}}$,
 and let $r^{\epsilon}_{\sigma^{\epsilon}_{n}} :=
  r^{\epsilon}_{\sigma^{\epsilon}_{n-1} } + 1$
 and $A^{\epsilon}_{\sigma^{\epsilon}_{n}}:= A^{\epsilon}_{\sigma^{\epsilon}_{n}} \cup \{  ( r^{\epsilon}_{\sigma^{\epsilon}_{n}} ,i) ; \, 1 \leq i \leq a \}  $.
 Then, for $x=   r^{\epsilon}_{\sigma^{\epsilon}_{n}}$, $i \in \{ 1, \ldots, a \}$, and $t \geq \sigma^{\epsilon}_{n}$, let
  $F^{\epsilon}_{t}(x,i):=x+Y_{x,i,t-\sigma^{\epsilon}_{n}}$. We shall see that $\sup_n \sigma^{\epsilon}_{n}=+\infty$ a.s.

 From the results in~\cite{ComQuaRam2}, (where only the case $\epsilon=0$ is treated, but it is immediate to adapt them to the present
setting), the following results hold.
For any $0 \leq \epsilon<1/2$ and $w \in \L_{\theta}$, almost surely with respect to $\P$:
\begin{itemize}
\item  for every $n\geq 1$,  $\sigma^{\epsilon}_{n-1}<\sigma^{\epsilon}_{n}<+\infty$, and there is a unique
$(x,i) \in A_{\sigma^{\epsilon}_{n-1}}$ such that
$F^{\epsilon}_{ \sigma^{\epsilon}_{n}}(x,i)=r^{\epsilon}_{\sigma^{\epsilon}_{n}}$;
\item $\lim_{n \to +\infty} \sigma^{\epsilon}_n= +\infty$;
\item for all $t \geq 0$, $X_{t}^{\epsilon}(w)   \in \L_{\theta}$;
\item the map $t \mapsto X_{t}^{\epsilon}(w)$ belongs to $\D(\L_{\theta})$.
\end{itemize}

For any $0 \leq \epsilon<1/2$, $\theta>0$, and $w=(F,r,A) \in \L_{\theta}$,
let
$\Q^{\epsilon,\theta}_{w}$ denote the probability distribution of
the random process $(X_{t}^{\epsilon}(w))_{t \geq 0}$,
viewed as a random element of $\D(\L_{\theta})$.
Again, as in~\cite{ComQuaRam2},
\begin{prop}
For any $0 \leq \epsilon<1/2$ and $\theta>0$, the family of probability measures $(\Q^{\epsilon,\theta}_{w})_{w \in \L_{\theta}}$
defines a strong Markov process on $\L_{\theta}$.
\end{prop}

In the sequel, we use $\E$ to denote expectation with respect to $\P$ of random variables defined on $(\Omega, \F)$.
The notation $\E^{\epsilon,\theta}_{w}$ is used to denote the expectation with respect to $\Q^{\epsilon,\theta}_{w}$ of random variables
defined on $\D(\L_{\theta})$ equipped with its Borel $\sigma-$field.

\subsection{Properties of hitting times}

For  $w=(F,r,A) \in \L_{\theta}$, and $u \geq r$, we define the
first time that the front touches site $u$, given that the initial condition was $w$,
$$T^{\epsilon}_{w}(u):= \inf \{ t > 0; \, r^{\epsilon}_{t} = u \}.$$

For all $u,v \in \Z$ such that $u<v$, $1 \leq i \leq a$, and  $0 \leq \epsilon < 1/2$, let

\begin{equation}
\label{def-a}
\A^{\epsilon}(u,i,v) := \inf \left\{ \sum_{k=1}^{m} \tau_{k}(u,i);  \,  u +
  \sum_{k=1}^{m} \varepsilon_{k}(u,i,\epsilon)  = v ,  \, m \geq 1\right\}.
\end{equation}
 This represents the first time that the random walk born
at $(u,i)$ hits site $v$ (assuming that the walk starts at $u$ at time zero).

 \begin{prop}\label{p:charact-T}
Let $w = (F,r,A) \in \L_{\theta}$.
 \begin{itemize}

\item[(i)] For all  $u > r$ and  $0 \leq \epsilon < 1/2$, $\P-$a.s.
  \begin{equation*} T^{\epsilon}_{w}(u) = \inf \sum_{j=1}^{L-1} \A^{\epsilon}(x_{j},i_{j},x_{j+1}),\end{equation*}
 where the infimum is taken over all finite sequences with $L \geq 2$, $x_{1},\ldots,x_{L} \in \mathbb{Z}$ and $i_{1},\ldots, i_{L-1}$ such that
 $x_{1}=F(y_{1}, i_{1})$ for some $(y_{1},i_{1}) \in A$,  $ r<x_{2} < \cdots < x_{L-1}<u, \, x_{L}=u$, $i_{2},\ldots, i_{L-1} \in \{1, \ldots, a \}$.

\item[(ii)] For all  $u > r$ and  $0 \leq \epsilon < 1/2$, the following identity holds $\P-$a.s.
 \begin{equation*} T^{\epsilon}_{w}(u) = \inf_{(x,i) \in A}  T^{\epsilon}_{w(x,i)}(u)
\end{equation*}

 \item[(iii)] For all  $r<u<v$ and  $0 \leq \epsilon < 1/2$,  the following  sub-additivity property holds  $\P-$a.s.
 $$T^{\epsilon}_{w}(v) \leq T^{\epsilon}_{w}(u) + T^{\epsilon}_{w \oplus (u-r)}(v).$$


\item[(iv)] \label{p:temps-monotone}
For any $0 \leq \epsilon_{1} \leq \epsilon_{2} < 1/2$, and all $u>r$,
$\P-$almost surely, $T^{\epsilon_{1}}_{w}(u) \geq T^{\epsilon_{2}}_{w}(u)$.
\end{itemize}
 \end{prop}

 \begin{proof}
  The proof of (i) is quite similar to that in~\cite{RamSid}, and so is the proof that (iii) is a consequence
 of (i). Then (ii) is a simple consequence of (i). 
As for (iv), this is an easy consequence of the characterization in (i)
and of the fact that, for every $(x,i) \in \Z \times \{1, \ldots, a \}$ and $n \geq 1$,
$\varepsilon_{n}(x,i,\epsilon_{1}) \leq \varepsilon_{n}(x,i,\epsilon_{2})$.
\end{proof}

An immediate consequence of (iv) in the above proposition is the following result.
\begin{coroll}\label{c:position-monotone}
For any $w \in L_{\theta}$, $0 \leq \epsilon_{1} \leq \epsilon_{2} < 1/2$,
$\P-$almost surely, for all $t \geq 0$, $r^{\epsilon_{1}}_{t}(w) \leq r^{\epsilon_{2}}_{t}(w)$.
\end{coroll}

\section{Proof of the large deviations principle for $t^{-1}r_{t}$}\label{s:proof-large-dev}

\begin{prop}\label{p:ldp-partielle}
There exists a convex function $J \, : \, (0,+\infty) \to [0,+\infty)$ such that, for all $b \in (0,+\infty)$,
$$\lim_{n \to +\infty} n^{-1} \log \P( T^{0}_{\delta_{0}}(n) \leq b n          )  = -J(b).$$
\end{prop}

\begin{proof}

For any $b > 0$, and all $n \geq 1$, it is easily checked that $\P(T^{0}_{\delta_{0}}(n) \leq b n )>0$.
Then let $u_{n}(b):= \log \P( T^{0}_{\delta_{0}}(n) \leq b n )$.
Observe that, by subadditivity (part $(iii)$ of Proposition~\ref{p:charact-T}),
$T^{0}_{\delta_{0}}(n+m)   \leq  T^{0}_{\delta_{0}}(n)  +   T^{0}_{\delta_{0} \oplus n}(n+m) $.
Now, by part $(ii)$ of Proposition~\ref{p:charact-T}, $T^{0}_{\delta_{0} \oplus n}(n+m) \leq T^{0}_{\delta_{n}}(n+m)$, since the infimum characterizing
$T^{0}_{\delta_{0} \oplus n}(n+m) $ runs over a larger set than the infimum characterizing $T^{0}_{\delta_{n}}(n+m)$.
As a consequence, $T^{0}_{\delta_{0}}(n+m)   \leq  T^{0}_{\delta_{0}}(n)  +   T^{0}_{\delta_{n}}(n+m) $.
We deduce that, for all $m,n \geq 1$, and all $b,c > 0$,
\begin{equation}\label{e:subadd-evenements}\{ T^{0}_{\delta_{0}}(n)  \leq bn  \}  \cap     \{     T^{0}_{\delta_{n}}(n+m)  \leq cm  \}
\subseteq   \{ T^{0}_{\delta_{0}}(n+m)  \leq bn+cm  \}.\end{equation}

Now, observe that  $T^{0}_{\delta_{0}}(n)$ and  $ T^{0}_{\delta_{n}}(n+m)$ are independent random variables, since 
their definitions involve disjoint sets of independent random walks. As a consequence, 
\begin{equation}\label{e:ineg-temps-assoc} \P (  \{ T^{0}_{\delta_{0}}(n)  \leq bn  \}  \cap     \{     T^{0}_{\delta_{n}}(n+m)  \leq cm  \}         )   =   
   \P(  T^{0}_{\delta_{0}}(n)  \leq bn     )
\P (T^{0}_{\delta_{n}}(n+m)  \leq cm    ).\end{equation}
>From the above two relations~(\ref{e:subadd-evenements}),~(\ref{e:ineg-temps-assoc}), and the fact that, by translation invariance of the model,
$T^{0}_{\delta_{0}}(m) $ and  $T^{0}_{\delta_{n}}(n+m)$  possess the same distribution, we deduce that, for all $m,n \geq 1$, and all $b,c > 0$,
\begin{equation}\label{e:sous-add-u}u_{n+m} \left(   \frac{bn+cm}{n+m}  \right) \geq u_{n}(b) + u_{m}(c).\end{equation}
Applying Inequality~(\ref{e:sous-add-u}) above with $c=b$, we deduce that the sequence $(u_{n}(b))_{n \geq 1}$ is super-additive.
Since $u_{n}(b) \leq 0$ for all $n \geq 1$, we deduce from the standard subadditive lemma that there exists a non-negative real number $J(b)$ such that
$\lim_{n \to +\infty} n^{-1} u_{n}(b) = -J(b)$.
Moreover, by definition, $b \mapsto u_{n}(b)$ is non-decreasing, and so $b \mapsto J(b)$ is non-increasing.

To establish that $J$ is convex, consider $b,c$, such that $0 < b < c$, $t \in (0,1)$, $k\geq 1$, and
apply~(\ref{e:sous-add-u}) with $n_{k}:=\ceil{kt}$ and $m_{k}:=\floor{k(1-t)}$. For large enough $k$,
$  \frac{bn_{k}+cm_{k}}{n_{k}+m_{k}} \leq tb + (1-t)c$, so that
$u_{n_{k}+m_{k}}(    tb + (1-t)c ) \geq  u_{n_{k}}(b)    + u_{m_{k}}(c)$.
Taking the limit as $k$ goes to infinity,
we deduce that $J(tb+(1-t)c) \leq tJ(b)+(1-t)J(c)$.

\end{proof}

\begin{prop}\label{p:prop-fonction-de-taux}
The function $J$ defined  in Proposition~\ref{p:ldp-partielle} is  identically zero on $[v^{-1},+\infty)$, positive and decreasing on $(0, v^{-1})$.
\end{prop}

The proof of the above proposition makes use of the following result, which is the main result of Section \ref{s:speedup}.
\begin{prop}\label{p:speedup}
For any $c > v$,
\begin{equation*}\limsup_{t \to +\infty} t^{-1} \log \P( r^{0}_{t}(\I_{0}) \geq ct ) < 0.\end{equation*}
\end{prop}

\begin{proof}[ Proof Proposition~\ref{p:prop-fonction-de-taux}]
For $n \geq 1$, $(ii)$ of Proposition~\ref{p:charact-T} implies that
$T^{0}_{\I_{0}}(n) \leq T^{0}_{\delta_{0}}(n)$ $\P-$a.s.
In view of the immediate identity $\{T^{0}_{w}(n) \leq bn\}=\{r^{0}_{bn}(w)\geq n \}$,
we deduce that $$\P(T^{0}_{\delta_{0}}(n) \leq bn ) \leq \P(    r^{0}_{bn}(\I_{0})\geq n ).$$
>From Proposition~\ref{p:speedup}, we deduce that $J$
is positive on $(0, v^{-1})$.
 On the other hand, by the law of large numbers~(\ref{e:lgn-presque-sure}), we see that $J$ must be identically $0$
 on $(v^{-1},+\infty)$. The function $J$ being convex on $(0,+\infty)$,
it is also continuous, so that $J(v^{-1})=0$.  Moreover, as we have already noted, $J$ is non-increasing.
These facts imply that $J$ is decreasing on
$(0, v^{-1})$.
\end{proof}

Let $I$ be defined by $I(b):=bJ(b^{-1})$ for $b>0$ and $I(0):=0$.
>From the previous results on $J$, it is easy to deduce the following.
\begin{coroll}
The function $I$ is identically zero on $[0,v]$, positive, increasing and convex on $(v,+\infty)$.
\end{coroll}

\begin{proof}
Only the convexity of $I$ is not totally obvious. Note that, since $J$ is convex, $b \mapsto J(b^{-1})$ is convex on
$(0,+\infty)$ as the composition of two convex functions. Then, since $b \mapsto J(b^{-1})$ is also increasing and positive,
the convexity of $b \mapsto b J(b^{-1})$ on $(0,+\infty)$ follows easily.

\end{proof}

\begin{prop}\label{p:indep-cond-init}
Assume that the initial condition $w$ satisfies $r=0$ and~{\bf (G)}. Then,  for all $b > 0$,
$$\lim_{n \to +\infty} n^{-1} \log \P( T^{0}_{w}(n) \leq b n          )  = -J(b),$$
where $J$ is the function defined in Proposition~\ref{p:ldp-partielle}.
\end{prop}

The proof of the proposition makes use of the following lemma.
\begin{lemma}\label{l:neglig-temps-atteinte}
Let $w=(F,r,A) \in \L_{\theta}$.
For all $t \geq 0$, and all $\gamma>0$,
$$\P \left(   \sup_{ (x,i) \in A} \sup_{0 \leq s \leq t}   F^{0}_{s}(x,i) \geq  r + \gamma t   \right) \leq
  f_{\theta}(w) \exp \left[  -g_{\gamma}(\theta) t \right],$$
where $$g_{\gamma}(\theta) := \gamma \theta - 2 (\cosh \theta - 1).$$
\end{lemma}

\begin{proof}
Let $G:= \left\{        \sup_{ (x,i) \in A} \sup_{0 \leq s \leq t}   F^{0}_{s}(x,i) >  r + \gamma t      \right\}$.
For all integers $K \leq 0$, let $G_{K}:=  \bigcup_{(x,i) \in A; \, F(x,i) \geq K} \{     \sup_{0 \leq s \leq t}   F^{0}_{s}(x,i) >  r + \gamma t  \} $.
Clearly, $K_{1} \leq K_{2}$ implies that  $G_{K_{2}} \subset G_{K_{1}}$, and $\bigcup_{K \leq 0} G_{K} = G$, whence
$\P(G) = \lim_{K \to -\infty} \P(G_{K})$.
Now observe that, for all $K$, the process $(M_{K,s})_{s \geq 0}$ defined by
$M_{K,s}:=\sum_{(x,i) \in A; \, F(x,i) \geq K}   \exp \left(\theta  \left(  F^{0}_{s}(x,i) -r \right) -2 (\cosh \theta - 1) s   \right)$
is a c\`adl\`ag martingale. Then note that, for all $K$,
$G_{K} \subset \{  \sup_{0 \leq s \leq t}  M_{K,s}  \geq    \exp(   g_{\gamma}(\theta) t \}$, then apply the
martingale maximal inequality to deduce that, 
$$\P  \left(  \sup_{0 \leq s \leq t}  M_{K,s}  \geq    \exp(   g_{\gamma}(\theta) t   )    \right) \leq \sum_{(x,i) \in A; \, F(x,i) \geq K} \exp\left[\theta \left(F(x,i) - r \right)  -g_{\gamma}(\theta) t \right].$$
For all $K$, the r.h.s. of the above inequality is upper bounded by
$ \leq  f_{\theta}(w) \exp \left[  -g_{\gamma}(\theta) t \right].$
The conclusion follows.
\end{proof}

\begin{proof}[Proof of Proposition~\ref{p:indep-cond-init}]

Consider $0<b<v^{-1}$, and fix $\theta>0$. Choose $\gamma>0$ large enough so that
\begin{equation*}g_{\gamma}(\theta) b> J(b).\end{equation*}
Denote by $w=(F,r,A)$ the initial condition, and consider the set $B_{n} := \{ (x,i); \,  F(x,i) \leq -\ceil{\gamma b n} \}$.
Let $m_{n}:= \sum_{(x,i) \in B_{n}}  \exp( \theta (F(x,i) - \ceil{\gamma b n})$.
Now let $\Xi_{n} :=   \inf \{s \geq 0; \, \exists  (x,i) \in B_{n}, \,   F^{0}_{s}(x,i) = 0 \}$.
We see that  $\Xi_{n}\leq bn $ implies that    $\sup_{ (x,i) \in B_{n}} \sup_{0 \leq s \leq bn}   F^{0}_{s}(x,i)   \geq  0$.
Thanks to Lemma~\ref{l:neglig-temps-atteinte} and translation invariance of the model, we deduce that
\begin{equation}\label{e:temps-loin-grand}\P(\Xi_{n} \leq  bn )  \leq m_{n}  \exp(-g_{\gamma}(\theta) bn).\end{equation}
>From  the fact that $w$ satisfies~{\bf(G)}, we obtain that, for all $\varphi>0$,
$y \leq 0$,
$\# \{ (x,i) \in A ; \,           F(x,i)  =  y  \} \leq f_{\varphi}(w) \exp( -\varphi y).$
As a consequence, whenever $\varphi<\theta$, we have that
\begin{equation}\label{e:peu-de-particules}m_{n}  \leq  f_{\varphi}(w) (1 - \exp( \varphi-\theta) )^{-1}  \exp( \varphi \ceil{\gamma b n} ) .\end{equation}
Now consider $(x,i) \in A\setminus B_{n}$,  so that $F(x,i) > -\ceil{\gamma b n}$.  By an easy coupling argument, we see that, since $F(x,i) \leq 0$,
\begin{equation}\label{e:une-particule-a-droite}\P(T_{w(x,i)} \leq  bn  )  \leq  \P(T_{\delta_{0}} \leq  bn  ).\end{equation}
Moreover, according to~{\bf(G)},
\begin{equation}\label{e:petit-card}    \#   A \setminus B_{n}  \leq    f_{\varphi}(w) \exp(\varphi \ceil{\gamma b n}  ) .\end{equation}
Now, by $(ii)$ of Proposition~\ref{p:charact-T}, $$\{    \Xi_{n} > bn    \} \cap   \{      T_{w}(n) \leq bn      \}  \subset  \left\{    \inf_{   (x,i) \in  A \setminus B_{n}}  T_{w(x,i)} \leq  bn  \right\}.$$
We deduce from~(\ref{e:temps-loin-grand}), (\ref{e:peu-de-particules}), (\ref{e:une-particule-a-droite}), ~(\ref{e:petit-card}) and the 
union bound that
\begin{equation}\label{e:presquefini} \P( T_{w}(n) \leq bn  ) \leq  f_{\varphi}(w)  e^{\varphi \ceil{\gamma b n}  }  \left[(1 - \exp( \varphi-\theta) )^{-1} \exp(-g_{\gamma}(\theta) bn) +
    \P(T_{\delta_{0}} \leq  bn  ) \right].\end{equation}
Now, according to Proposition~\ref{p:ldp-partielle},
$$\lim_{n \to +\infty} n^{-1} \log \P( T^{0}_{\delta_{0}}(n) \leq b n          )  = -J(b).$$
Since we have chosen $\gamma$ so that $g_{\gamma}(\theta) b> J(b)$, we deduce from~(\ref{e:presquefini}) that
$$\limsup_{n \to +\infty} n^{-1} \log \P( T^{0}_{w}(n) \leq b n          )  \leq -J(b) +  \varphi \gamma b.$$
Since $\varphi>0$ is arbitrary, we deduce that
\begin{equation}\label{e:la-limite-sup}\limsup_{n \to +\infty} n^{-1} \log \P( T^{0}_{w}(n) \leq b n          )  \leq -J(b).\end{equation}
On the other hand, consider a given $(x,i) \in A$.
Clearly
$$\P( T_{w}(n) \leq bn  ) \geq   \P( T_{w(x,i)}(n) \leq bn  ).$$ Now consider $\tilde{\tau} = \inf \{  s \geq 0; \, F^{0}_{s}(x,i) = 0 \}$.
Clearly, $\tilde{\tau}$ is a.s. finite, and, conditional upon $\tilde{\tau}$, $ T_{w(x,i)}(n) - \tilde{\tau}$ has the (unconditional) distribution of $T_{\delta_{0}}(n)$.
Choosing any $M$ such that  $\P(\{  \tilde{\tau} \leq M  \}    )>0$,  one has that $ \P( T_{w(x,i)}(n) \leq bn  )    \geq     \P(\{  \tilde{\tau} \leq M  \}    )       \P( T_{\delta_{0}}(n) \leq bn-M)$.
Taking an arbitrary $c>b$, we deduce that $$\liminf_{n \to +\infty} n^{-1} \log \P( T^{0}_{w}(n) \leq b n          ) \geq -J(c).$$
By continuity of $J$, we conclude that
$$\liminf_{n \to +\infty} n^{-1} \log \P( T^{0}_{w}(n) \leq b n          ) \geq -J(b).$$
The above inequality, together with~(\ref{e:la-limite-sup}) concludes the proof.
\end{proof}

\begin{proof}[ Proof of Theorem~\ref{t:ldp}]

Consider a non-empty closed subset $F \subset [0,+\infty)$, and let $b:= \inf F$.
 Assume that $b \leq v$. We have that $\inf_{F} I = 0$, so the upper bound of the LDP for $F$ is always satisfied.
Assume now that $b>v$.
One has that $\P(  t^{-1}r^{0}_{t}(w)  \in F ) \leq \P( r^{0}_{t}(w) \geq \ceil{t b} )
= \P( T^{0}_{w} (\ceil{t b}) \leq t)$.
Proposition~\ref{p:indep-cond-init} entails that $\lim_{t \to +\infty}  t^{-1} \log \P( T^{0}_{w}(\ceil{t b} \leq t) \leq -I(b)$,
so that the upper bound of the LDP holds for $F$ since $I$ is non-decreasing.

Consider now an open set $G \subset (v,+\infty)$. For every $b \in G$, there exists an interval
 $[b,c) \subset G$.
  By the  large deviations upper bound, we know that $\limsup_{t \to +\infty} t^{-1} \log \P(r^{0}_{t}(w) \geq bt)
 \leq -I(b)$ and that  $\limsup_{t \to +\infty} t^{-1} \log \P(r^{0}_{t}(w) \geq ct)
 \leq -I(c)$.
By strict monotonicity of $I$ on $(v,+\infty)$, we have that $I(b)<I(c)$, so we can conclude that
$\liminf_{t \to +\infty} t^{-1} \log \P(bt \leq r^{0}_{t}(w) < ct) \geq -I(b)$.
As a consequence, $\liminf_{t\to\infty}\frac{1}{t} \P ( t^{-1} r^{0}_{t}(w) \in G  ) \geq -I(b)$.
 Since this holds for an arbitrary $b \in G$, the lower bound of the LDP for $G$ follows.

Consider now a non-empty open set $G \subset [0,+\infty)$ such that $G \cap [0,v] \neq \emptyset$.
Then $\inf_{G} I = 0$. On the other hand,
there is a non-empty interval of the form  $[c,b) \subset G \cap [0,v]$.
In Section~\ref{s:slowdown}, we prove that, under Assumption~{\bf(G)},
\begin{equation}\label{e:slowdown-lower-bound-gen}
\liminf_{t\to\infty}t^{-1}\log \P\left[ c \leq \frac{r^{0}_t(w)}{t}\leq b\right]=0.
\end{equation}
Applying Inequality~(\ref{e:slowdown-lower-bound-gen}), we see that
$ \liminf t^{-1} \log \P(  t^{-1} r^{0}_{t} \in  G  ) = 0$, so that the lower bound of the LDP holds.
\end{proof}

\section{Speedup probabilities}\label{s:speedup}

The main result in this section is Proposition~\ref{p:speedup}:
\begin{equation}\label{e:dev-exp}  \mbox{ for any $b > v$, } \limsup_{t \to +\infty} t^{-1} \log \P( r^{0}_{t}(\I_{0}) \geq bt ) < 0.\end{equation}
For the sake of readability, the reference to the initial condition $\I_{0}$ is often dropped in this section, so that $r^{\epsilon}_{t}$ should be read as
$r^{\epsilon}_{t}(\I_{0})$.

Our strategy for proving Proposition~\ref{p:speedup} is to exploit the renewal structure already used in~\cite{ComQuaRam2} to prove the CLT.
However, this renewal structure leads to random variables (renewal time, and displacement of the front at a renewal time) whose tails have polynomial decay
 (see~Appendix~\ref{a:poly-decay}, and asymptotic exponential bounds such as (\ref{e:dev-exp}) cannot be derived from such random variables.
 Whether it is possible to modify the definition of the renewal structure so as to obtain random variables enjoying an exponential decay of the tails, as required
 for a direct proof of Proposition~\ref{p:speedup} is unclear and instead we make use of a different idea. Indeed, we apply the renewal structure defined
 in~\cite{ComQuaRam2} to a perturbation of the original model, one in which the random walks have a small bias to the right.
  Again, a law of large numbers holds:
\begin{prop}\label{p:lgn-position}
For all small enough $\epsilon \geq 0$, there exists $0<v_{\epsilon}<+\infty$ such that
$$\lim_{t\to\infty} t^{-1} r^{\epsilon}_t=v_{\epsilon} , \,  \P-a.s. \mbox{ and in $L^{1}(\P)$}.$$
\end{prop}

 The interest of introducing a bias to the right is that, reworking the estimates of~\cite{ComQuaRam2} in this context, we can show
 that for any small value of the bias parameter $\epsilon>0$, exponential decay of the tail of the renewal times holds, so that the following result can be proved.

 \begin{prop}\label{p:speedup-biais}
There exists $\epsilon_{0}>0$ such that, for any $0 \leq \epsilon \leq \epsilon_{0}$, for any $b > v_{\epsilon}$,
\begin{equation*}\limsup_{t \to +\infty} t^{-1} \log \P( r^{\epsilon}_{t} \geq bt ) < 0.\end{equation*}
\end{prop}

 On the other hand, it is shown in Corollary~\ref{c:position-monotone} above that, as expected,
 biasing the random walks to the right cannot decrease the position of the front,
 so that at each time $t$, a comparison holds between the position of the front in the original model and in the model with a bias.
 We deduce that
  \begin{prop}\label{p:compare-loi}
 For any $0 \leq \epsilon < 1/2$ and $t \geq 0$, and all $x \in \{1,2,\ldots \}$,
 $$\P( r^{0}_{t} \geq  x) \leq \P( r^{\epsilon}_{t} \geq  x).$$
  \end{prop}

 As a consequence, we can prove that~(\ref{e:dev-exp}) holds for all $b$ such that there exists an $0 \leq \epsilon \leq \epsilon_{0}$
for which $v_{\epsilon} < b$. Noting that $v_{\epsilon}$ is a non-decreasing function of $\epsilon$, we see that  the
following result would make our strategy work for all $b>v$:
\begin{prop}\label{p:cv-vitesse}
 \begin{equation}\label{e:cv-vitesse}\lim_{\epsilon \to 0+} v_{\epsilon} = v.\end{equation}
  \end{prop}

It is indeed natural to expect such a continuity property to hold, but proving it seems to require substantial work.

Indeed, write
\begin{equation}\label{e:cv-en-t}v_{\epsilon} =    \lim_{t \to +\infty}  t^{-1}  \E(r^{\epsilon}_{t}).\end{equation}
\begin{equation*} v =    \lim_{t \to +\infty}  t^{-1}  \E(r^{0}_{t}).\end{equation*}
For fixed $t$, it is possible (using the dominated convergence theorem) to prove
that \begin{equation}\label{e:cv-en-epsilon}\lim_{\epsilon \to 0+} \E (r^{\epsilon}_{t}) =
\E(r^{0}_{t}).\end{equation}
Hence, to prove Identity~(\ref{e:cv-vitesse}), it is enough to prove that
$$\lim_{\epsilon \to 0+}  \lim_{t \to +\infty}  t^{-1}  \E(r^{\epsilon}_{t})
=  \lim_{t \to +\infty}  \lim_{\epsilon \to 0+}    t^{-1}  \E(r^{\epsilon}_{t}).$$

Our strategy for proving Proposition~\ref{p:cv-vitesse} is based on the observation that, if some sort of uniformity
with respect to $ \epsilon \in [0 , \epsilon_{0}]$ is achieved in~(\ref{e:cv-en-t}), then
the limits with respect to $\epsilon \to 0+$ and to $t\to +\infty$ in~(\ref{e:cv-en-t})-(\ref{e:cv-en-epsilon}) can be exchanged.
Reworking the estimates in~\cite{ComQuaRam2} to obtain uniform upper bounds
 (with respect to  $0 \leq \epsilon \leq \epsilon_{0}$) for the second moments of the random
 variables (renewal time, and displacement of the front at a renewal time) defined by the renewal structure, we can prove that the required
 uniformity in~(\ref{e:cv-en-t}) holds.

\subsection{Some random variables on $\D(\L_{\theta})$}\label{s:srvod}

It will be convenient in the sequel to work with random
variables defined on the canonical space of trajectories 
$\D(\L_{\theta})$ rather than on $(\Omega, \F, \P)$.
We use the $\hat{ }$ sign in order to make apparent the distinction between random variables defined on $\Omega$ and
their counterparts. on  $\D(\L_{\theta})$.

On $\D(\L_{\theta})$, we define the following random variables.
Let $w_{\cdot}=(w_{t})_{t \geq 0}=(\hat{F}_{t},\hat{r}_{t},\hat{A}_{t})_{t \geq 0} \in \D(\L_{\theta})$.
The random process $(\hat{r}_{t})_{t \geq 0}$ is defined through  $(w_{t})_{t \geq 0}=(\hat{F}_{t}, \hat{r}_{t},\hat{A}_{t})_{t \geq 0} $. 
Under the probability measure $\Q^{\epsilon,\theta}_w$ the process $(\hat F_t,\hat r_t,\hat A_t)_{t\geq 0}$ has the same
law as $(F^\epsilon_t,r^\epsilon_t,A^\epsilon_t)_{t\ge 0}$.

For all $s \geq 0$, let  $Z_{s,x,i}(w_{\cdot}) :=  x $ if $(x,i) \notin A_{s}$, and $Z_{s,x,i}(w_{\cdot}) =\hat{F}_{s}(x,i)$ otherwise.
For $y \in \Z$, let $\hat{T}(y) := \inf \{  s \geq 0; \, \hat{r}_{s}=y  \}$ if $y \geq \hat{r}_0+1$,
 and let $\hat{T}(y):=0$ otherwise.
Let also $G_{s,x,i}(w_{\cdot}) := Z_{\hat{T}(x)+s , x ,i}$.
With respect to $\Q^{\epsilon,\theta}_{w}$, the processes $(G_{s,x,i})_{s \geq 0}$ form a family independent
nearest-neighbor random walks on $\Z$ with jump rate 2 and step distribution $(1/2+\epsilon) \delta_{+1} + (1/2-\epsilon) \delta_{-1}$.

For $z \in \Z$, and $w=(F,r,A) \in \L_{\theta}$,  define $\phi_{z}(w)$ by
$$\phi_{z}(w) :=   \sum_{(x,i) \in A \cap \{  \ldots ,z-1,z \} \times \{1,\ldots, a \}} \exp( \theta (F(x,i)-r) ),$$
and for $z_{1} < z_{2} \in \Z$, let
$$
m_{z_1,z_2}(w):=\sum_{(x,i) \in A \cap \{z_1+1,\ldots, z_2\} \times \{1,\ldots, a \}}  \un(  z_{1}+1  \leq F(x,i) \leq z_{2}    ).
$$

We use the notation $\theta_{s}$ to denote the canonical time-shift on $\D(\L_{\theta})$ and  the notation $\varpi_{y}$ to denote the truncated space-shift on
$\D(\L_{\theta})$ defined by
$\varpi_{y}(F,r,A) = (F',r',A')$, with $A'= \{  (x-y,i); \, (x,i) \in A  , \, x \geq y  \}$,
$F'(x):=F(x+y)$, $r':=r-y$. In words, this corresponds to removing all the particles that are born at the left of $y$, and then shifting all
birth positions by $y$.
We denote by $(\F^{\epsilon,\theta}_{t})_{t \geq 0}$ the usual augmentation of the
natural filtration on $\D(\L_{\theta})$ with respect to the Markov family $(\Q^{\epsilon,\theta}_{w})_{w \in \L_{\theta}}$.

\subsection{An elementary speedup estimate}

The following lemma is stated in~\cite{ComQuaRam2} in the
  case $\epsilon=0$, and its adaptation to the more general case $0 \leq \epsilon <1/2$ is straightforward.
 
\begin{lemma}\label{l:lemme10}
Let $\lambda(\epsilon,\theta):=2 (\cosh \theta - 1) + 4 \epsilon \sinh \theta + a(1+2\epsilon) \exp \theta$ and
$c_{\gamma}(\epsilon,\theta):= \gamma \theta - \lambda(\epsilon,\theta)$. For all $0 \leq \epsilon <1/2$, $w \in \L_{\theta}$, and $t \geq 0$,
$$\Q^{\epsilon,\theta}_{w}(\hat{r}_{t} - \hat{r}_0 \geq \gamma t) \leq \phi_{\hat{r}_{0}}(w) \exp( - c_{\gamma}(\epsilon,\theta) t   ).$$
\end{lemma}

\subsection{Definition of the renewal structure}\label{s:def-renouvellement}

We follow the definition of the renewal structure in \cite{ComQuaRam2}.
Consider a parameter \begin{equation}\label{e:valeur-de-M}M:= 4(a+9).\end{equation}
 Let $\nu_{0}:=0$ and $\nu_{1}$
 be the first time one of the random walks
 $\{(G_{s,r_{0},i})_{s \geq 0}; \, 1 \leq i \leq a \} $, hits the site $\hat{r}_{0}+1$ (the random walks $(G_{s,x,i})$ are defined
 in section~\ref{s:srvod}).
Next, define $\nu_2$ as the first time one of the random
walks $\{(G_{s,z,i})_{s \geq 0} ; \, \hat{r}_{0}\le z\le \hat{r}_{0}+1, \, 1 \leq i \leq a \}$ hits the site $\hat{r}_{0}+2$.
In general, for $k\ge 2$, we define $\nu_k$ as the first time
one of the  random walks
 $\{(G_{s,z,i})_{s \geq 0}; \, \hat{r}_{0} \lor (\hat{r}_{0}+k-M)\le z\le \hat{r}_{0}+k-1, \, 1 \leq i \leq a \}$,
 hits the
site $\hat{r}_{0}+k$. For $n\in {\mathbb N}$, let
\begin{equation}\nonumber
\tilde r_t:=\hat{r}_{0}+n,\qquad{\rm if}\qquad
\sum_{k=0}^n\nu_k\le t<\sum_{k=0}^{n+1}\nu_k.
\end{equation}
The following proposition (see Lemma 1 from~\cite{ComQuaRam2}), shows that the so-called auxiliary front $\tilde r_t$
can be used to estimate the position of the front $\hat{r}_{t}$.
\begin{prop}
For every $0 \leq \epsilon <1/2$, $\theta>0$ and $w \in \L_{\theta}$,
the following holds $\Q^{\epsilon,\theta}_{w}-$almost surely:
$$ \mbox{ for every $t \geq 0$, } \tilde r_t \leq \hat{r}_{t}.$$
\end{prop}

Now, observe that for any $w=(F,r,A)$ such that $r \times \{1,\ldots,a\} \subset A$ and $F(r,i)= r$ for all $1 \leq i \leq a$, 
with respect to $\Q^{\theta,\epsilon}_{w}$,
for each $1\le j\le M-1$, the random variables $(\nu_{i})_{i \geq 1}$ are a.s. finite, and that the random variables
$\{\nu_{Mk+j}:k\ge 1\}$ are i.i.d. and have finite expectation since $M\ge 3$. We deduce that
a.s. (see also \cite{ComQuaRam}),
\begin{equation*}
\lim_{t\to\infty} \tilde r_t/t=:\alpha(\epsilon)>0.
\end{equation*}
First note that $\alpha(\epsilon)$ does not depend on $\theta$ nor on $w$ since the distribution of the random
walks $(G_{s,x,i})_{s \geq 0}$ with respect to $\Q^{\theta,\epsilon}_{w}$ does not.
Moreover, $\alpha(\epsilon)$ is a non-decreasing function of $\epsilon$ by an immediate coupling
argument.

Now consider $\epsilon_{0}<1/2, \theta>0, \alpha_{1},\alpha_{2}>0$ such that
\begin{equation}\left\{ \begin{array}{l}\label{e:alph1-et-2}
0<\alpha_{1}<\alpha_{2}<\alpha(0), \\
\theta^{-1}(   2 (\cosh \theta  -1) + 4 \epsilon_{0}  \sinh \theta  )   < \alpha_{1},\\
4 \epsilon_{0}<\alpha_{1}.\end{array}
\right.\end{equation}

In the sequel, we always assume that  $0 \leq \epsilon \leq \epsilon_{0}$.

Let us define the following random variables on $\D(\L_{\theta})$:
\begin{equation*}
\nonumber
\left\{ \begin{array}{l}
U(w_{\cdot}):=\inf \{t\ge  0; \, \tilde r_t -\hat{r}_{0}< \floor{\alpha_2 t} \},\\
V( w_{\cdot}):=\inf\{t\ge 0; \,  \max_{\hat{r}_{0}-L+1 \leq z \leq \hat{r}_{0}-1}   Z_{t,x,i}   >     \floor{\alpha_1 t} + \hat{r}_{0} \},\\
W(w_{\cdot}):=\inf\{t\ge 0; \, \phi_{\hat{r}_{0}-L}(w_{t})      \ge e^{\theta (\floor{\alpha_1 t} -(\hat{r}_t-\hat{r}_{0}))}\}.\end{array} \right.
\end{equation*}

Note that, for all $\epsilon$, $U,V,W$ are stopping times with respect to $(\F^{\epsilon,\theta}_{t})_{t \geq 0}$, and that they are mutually independent with respect
 to  $\Q^{\theta,\epsilon}_{w}$.

Let \begin{equation*}   D:= \min(U, V, W).\end{equation*}

Now let $p>0$ be such that
\begin{equation*}  p \exp(\theta) < 1, \end{equation*}
and $L$ such that
\begin{equation}
\label{e:L}
L^{1/4} \geq M+1 \quad {\rm and} \quad   a \exp(-L \theta)(1- \exp(-\theta))^{-1}<p.
\end{equation}
For $x \in \Z$, let
\begin{equation*}
J_x(w_{\cdot}) := \inf\{j\ge 1; \,  \phi_{x+(j-1)L}(w_{\hat{T}(x+jL)})\le p , \,
m_{x+jL-L^{1/4},x+jL}(w_{\hat{T}(x+jL)})\ge a \floor{L^{1/4}}/2\}.
\end{equation*}
 Let $S_0:=0$ and $R_0:=\hat{r}_{0}$.
Then define for $k\ge 0$,
\begin{equation} \nonumber S_{k+1}:=\hat{T}(R_k+J_{R_k}L),\qquad
 D_{k+1}:=D\circ\theta_{S_{k+1}}+S_{k+1},\qquad R_{k+1}=r_{D_{k+1}}
\end{equation}
\begin{equation} \nonumber
K:=\inf\{k\ge 1:S_k<\infty ,D_k=\infty\},
\end{equation}
and define the {\it regeneration time}
\begin{equation*}
\kappa:=S_{K},
\end{equation*}
Note that $\kappa$ is {\em not} a stopping time with
respect to  $(\F^{\epsilon,\theta}_t)_{t \geq 0}$.
Define ${\mathcal G}^{\epsilon,\theta}$, the information up to time $\kappa$,
as the smallest $\sigma$-algebra containing all sets of the
form $\{\kappa\le t\}\cap A, \, A\in \F^{\epsilon,\theta}_t, \, t \geq 0$.

\subsection{Properties of the renewal structure}

Throughout this section, we assume that $\theta,\alpha_{1}, \alpha_{2}, \epsilon_{0}$
satisfy the assumptions listed in part~\ref{s:def-renouvellement}.

\begin{prop}\label{p:moments-renouvellement}
The following properties hold:
\begin{itemize}
 \item[(i)] There exist $0<C,L^{*}<+\infty$
not depending on $\epsilon$ (but possibly depending on the choice of
$\theta,\alpha_{1}, \alpha_{2}, \epsilon_{0}$) such that, for $L:=L^{*}$, and all $0 \leq \epsilon \leq \epsilon_{0}$,
$$\E^{\epsilon,\theta}_{\I_{0}} ( \kappa^{2} ) \leq C , \,
\E^{\epsilon,\theta}_{a \delta_{0}} ( \kappa^{2}| U=+\infty ) \leq C,$$
$$\E^{\epsilon,\theta}_{\I_{0}} ( (\hat{r}_{\kappa})^{2} ) \leq C \quad {\rm and}\quad
\E^{\epsilon,\theta}_{a \delta_{0}} ( (\hat{r}_{\kappa})^{2}| U=+\infty ) \leq C.$$
\item[(ii)] For all $0<\epsilon \leq \epsilon_{0}$, there exist $0<C(\epsilon),L(\epsilon),t(\epsilon)<+\infty$ such that,
for $L:=L(\epsilon)$,
$$\E^{\epsilon,\theta}_{\I_{0}} ( \exp( t(\epsilon) \kappa ) ) \leq C(\epsilon) , \,
\E^{\epsilon,\theta}_{a \delta_{0}} (  \exp( t(\epsilon) \kappa )  | U=+\infty ) \leq C(\epsilon),$$
$$\E^{\epsilon,\theta}_{\I_{0}} ( \exp( t(\epsilon) \hat{r}_{\kappa} ) ) \leq C(\epsilon) \quad {\rm and}\quad
\E^{\epsilon,\theta}_{a \delta_{0}} (  \exp( t(\epsilon) \hat{r}_{\kappa} )  | U=+\infty ) \leq C(\epsilon).$$
\end{itemize}

\end{prop}

Proposition~\ref{p:moments-renouvellement} provides the key estimates needed for the proof of the main results in this section.
Most of the technical work needed to prove it consists in a
reworking of the estimates in~\cite{ComQuaRam2},
either proving that, for each positive value of the bias parameter $\epsilon$,  exponential estimates can be obtained instead of the polynomial
ones derived in~\cite{ComQuaRam2}, or that the polynomial estimates already obtained in~\cite{ComQuaRam2}
can be made uniform with respect to $0 \leq \epsilon \leq \epsilon_{0}$.
The proofs go  along the lines of~\cite{ComQuaRam2}, and are deferred to Appendix  \ref{a:renewal-estimates}.
In the sequel, we always assume that either $L:=L^{*}$ or $L:=L(\epsilon)$.
As a consequence of Proposition~\ref{p:moments-renouvellement} we see
that
for all $0 \leq \epsilon \leq \epsilon_{0}$,
$\Q^{\epsilon,\theta}_{\I_{0}} (0<\kappa<+\infty) = 1 \mbox{ and }
\Q^{\epsilon,\theta}_{a \delta_{0}}(0<\kappa<+\infty | U=+\infty ) = 1$.

As in~\cite{ComQuaRam2}, the following propositions and corollary can be proved.
\begin{prop}\label{p:loi-renouvellement}
Let $0 \leq \epsilon \leq \epsilon_{0}$. If $w=\I_{0}$ or $w=a \delta_{0}$, then
for any Borel subset $\Gamma$ of $\D(\L_{\theta})$,
$$\Q^{\epsilon,\theta}_{w}( \varpi_{\hat{r}_{\kappa}}(w_{\kappa+t})_{t \geq 0} \in \Gamma    | \G^{\epsilon,\theta}    ) =
\Q^{\epsilon,\theta}_{a \delta_{0}} (\Gamma | U=+\infty   ) \qquad \Q^{\epsilon,\theta}_{w}-a.s.$$
\end{prop}

Define
$\kappa_{1}:=\kappa$ and for $i\ge 1$, $\kappa_{i+1}:=\kappa_{i}+\kappa \circ \theta_{\kappa_{i}}$.
Now, for all $i \geq 1$, define $\G^{\epsilon,\theta}_{i}$ as the smallest $\sigma$-algebra containing all sets of the
form $\{\kappa_{i} \le t\}\cap A, \, A\in \F^{\epsilon,\theta}_t , \, t \geq 0$.The following general version of Proposition~\ref{p:loi-renouvellement} holds.
\begin{prop}\label{p:extension-loi-renouvellement}
Let $0 \leq \epsilon \leq \epsilon_{0}$ and $i \geq 1$.  If
 $w=\I_{0}$ and $w=a \delta_{0}$ then
for any Borel subset $\Gamma$ of $\D(\L_{\theta})$,
$$\Q^{\epsilon,\theta}_{w}( \varpi_{\hat{r}_{\kappa_{i}}}(w_{\kappa_{i}+t})_{t \geq 0} \in \Gamma    | \G^{\epsilon,\theta}_{i}    )
= \Q^{\epsilon,\theta}_{a \delta_{0}} (\Gamma | U=+\infty   ) \qquad \Q^{\epsilon,\theta}_{w}-a.s.$$
\end{prop}

\begin{coroll}
\label{c:iid}
The following properties hold:
\begin{itemize}
\item[(i)]  Under $\Q^{\epsilon,\theta}_{\I_{0}}$,
 $\kappa_1,\kappa_2-\kappa_1, \kappa_3-\kappa_2, \ldots$ are independent,
 and $\kappa_2-\kappa_1, \kappa_3-\kappa_2, \ldots$ are identically distributed with law identical to that of $\kappa$ under
 $\Q^{\epsilon,\theta}_{a \delta_{0}} (\cdot | U=+\infty   )$.
\item[(ii)]  Under $\Q^{\epsilon,\theta}_{\I_{0}}$,
 $\hat{r}_{\kappa_1},\hat{r}_{\kappa_2}-\hat{r}_{\kappa_1}, \hat{r}_{\kappa_3}-\hat{r}_{\kappa_2}, \ldots$ are independent,
 and $\hat{r}_{\kappa_2}-\hat{r}_{\kappa_1}, \hat{r}_{\kappa_3}-\hat{r}_{\kappa_2}, \ldots$
 are identically distributed with law identical to that of $r_{\kappa}$ under
 $\Q^{\epsilon,\theta}_{a \delta_{0}} (\cdot | U=+\infty   )$.
\end{itemize}
\end{coroll}

We now give the proofs of Propositions~\ref{p:lgn-position},~\ref{p:speedup-biais} and~\ref{p:cv-vitesse}.

\begin{proof}[Proof of Proposition~\ref{p:lgn-position}]

First, note that the $\P-$a.s. convergence stated in Proposition~\ref{p:lgn-position} follows from the integrability of renewal times by a standard argument.
 To prove that the convergence also takes place in $L^{1}(\P)$, we note that, from Lemma~\ref{l:lemme10} above, it stems 
 that $\E^{\epsilon,\theta}_{\I_{0}} (\hat{r}_{t})<+\infty$ for all $t$ and that the family of random variables
$(t^{-1} \hat{r}_{t})_{t \geq 1}$ is uniformly integrable with respect to $\Q^{\epsilon,\theta}_{\I_{0}}$.
The convergence in $L^{1}(\P)$ then follows from the $\P-$a.s. convergence.

\end{proof}

\begin{proof}[Proof of Proposition~\ref{p:speedup-biais}]
Fix $0 < \epsilon \leq \epsilon_{0}$, and let $L:=L(\epsilon)$.
For all $t \geq 0$, define $a(t):= \sup\{  n \geq 1; \, \kappa_{n} \leq t  \},$
with the convention that $\sup \emptyset = 0$. From Corollary~\ref{c:iid} and Proposition~\ref{p:moments-renouvellement},
we deduce that, $a(t)<+\infty$ a.s. for all $t \geq 0$ and
that
 $\lim_{t \to +\infty} a(t) = +\infty$ a.s.
 Using the fact that the map $t \mapsto \hat{r}_{t}$ is non-decreasing, we have that
 $ \hat{r}_{t} \leq \hat{r}_{\kappa_{a(t)+1}}$.
  Now observe that, for any $0 < \epsilon \leq \epsilon_{0}$, any $b> v_{\epsilon}$, and any $0<c<+\infty$, by the union bound,
 $$ \Q^{\epsilon,\theta}_{\I_{0}}(  \hat{r}_{t} \geq bt  ) \leq \Q^{\epsilon,\theta}_{\I_{0}}(a(t) \geq \floor{ct})+
  \Q^{\epsilon,\theta}_{\I_{0}}(\hat{r}_{\kappa_{\floor{ct}+1} }\geq bt).$$
 Note that  $\Q^{\epsilon,\theta}_{\I_{0}}( a(t) \geq \floor{ct}   )
 \leq \Q^{\epsilon,\theta}_{\I_{0}}( \kappa_{ \floor{ct} } \leq t )$,
and observe that, by a standard large deviations bound for the i.i.d non-negative sequence $(\kappa_{i+1}-\kappa_{i})_{i \geq 1}$
 and Proposition~\ref{p:moments-renouvellement}
 for $\kappa_{1}$, whenever $c^{-1}< \E^{\epsilon,\theta}_{a \delta_{0}}(  \kappa  | U=+\infty)$,
  $\limsup_{t \to +\infty} t^{-1} \log \Q^{\epsilon,\theta}_{\I_{0}}(  \kappa_{\floor{ct}} \leq t)<0$.
On the other hand,  writing $\hat{r}_{\kappa_{\floor{ct}+1}} = \hat{r}_{\kappa_{1}}+ \sum_{i=1}^{\floor{ct} +1} (\hat{r}_{\kappa_{i+1}}- \hat{r}_{\kappa_{i}})$,
and using  Proposition~\ref{p:moments-renouvellement}, together with a standard large deviations argument (see e.g.~\cite{DemZei}),
we have that, as soon as $b/c >\E^{\epsilon,\theta}_{a \delta_{0}}(  \hat{r}_{\kappa}  | U=+\infty)  $,
  $\limsup_{t \to +\infty} t^{-1} \log \Q^{\epsilon,\theta}_{\I_{0}}(\hat{r}_{\kappa_{\floor{ct}+1} } \geq bt)<0$.

 Note that we can deduce from the renewal structure that
 \begin{equation}\label{e:caract-vitesse}v_{\epsilon} = \frac{\E^{\epsilon,\theta}_{a \delta_{0}}
 (  \hat{r}_{\kappa}  | U=+\infty)}{\E^{\epsilon,\theta}_{a \delta_{0}}(  \kappa | U=+\infty)}.\end{equation}
As a consequence, if $b>v_{\epsilon}$, we see that we can choose a $c>0$ such that
 $c^{-1}< \E^{\epsilon,\theta}_{a \delta_{0}}(  \kappa  | U=+\infty)$ and
 $b/c >\E^{\epsilon,\theta}_{a \delta_{0}}(  \hat{r}_{\kappa}  | U=+\infty)  $. 
\end{proof}

\begin{lemma}\label{l:mino-esp}
There exists $0<c<+\infty$ such that, for all $0 \leq \epsilon \leq \epsilon_{0}$,
$$\E^{\epsilon,\theta}_{a \delta_{0}}(  \kappa | U=+\infty) \geq c.$$
\end{lemma}

\begin{proof}
Use the fact that, by definition, $\kappa \geq \hat{T}(1)$, so that
$\E^{\epsilon,\theta}_{a \delta_{0}}(  \kappa | U=+\infty) \geq \E^{\epsilon,\theta}_{a \delta_{0}}(  \hat{T}(1) \un(U=+\infty))$.
Now, by coupling, $\Q^{\epsilon,\theta}_{a \delta_{0}}(U=+\infty) \geq \Q^{0,\theta}_{a \delta_{0}}(U=+\infty)$
for all $0 \leq \epsilon \leq \epsilon_{0}$. By coupling again, for all $u>0$,
$\Q^{\epsilon,\theta}_{a \delta_{0}}(\hat{T}(1) \geq u) \geq \Q^{\epsilon_{0},\theta}_{a \delta_{0}}(\hat{T}(1) \geq u)$.
Now, since $\Q^{\epsilon_{0},\theta}_{a\delta_{0}}(\hat{T}(1)>0)=1$, we can find $u>0$ small enough so that
$\Q^{\epsilon_{0},\theta}_{a \delta_{0}}(\hat{T}(1) \geq u) \geq 1 - (1/2)\Q^{0,\theta}_{a \delta_{0}}(U=+\infty)$.
Putting the previous inequalities together, we see that, for all $0 \leq \epsilon \leq \epsilon_{0}$,
$\Q^{0,\theta}_{a \delta_{0}}(\hat{T}(1) \geq u, U=+\infty) \geq (1/2)\Q^{0,\theta}_{a \delta_{0}}(U=+\infty)$.
The conclusion follows.
 \end{proof}

The following proposition contains the uniform convergence estimate that is required for the proof of Proposition~\ref{p:cv-vitesse}.
Broadly speaking, the idea is to control the convergence speed with second moment estimates on the renewal structure, so that uniform
estimates on these moments yield uniform estimates on the convergence speed.

\begin{prop}\label{p:limite-inf-uniforme}
 For all $\zeta>0$, there exists $t_{\zeta} \geq 0$
such that, for all $t \geq t_{\zeta}$ and all $0 \leq \epsilon \leq \epsilon_{0}$,
$$v_{\epsilon} \leq \E^{\epsilon,\theta}_{\I_{0}} (t^{-1} \hat{r}_{t}) + \zeta.$$
\end{prop}

\begin{proof}

Let $0<\lambda<1$  be given, and let 
$$m(t,\epsilon) :=   \left\lfloor  (1-\lambda) t \left( \E^{\epsilon,\theta}_{a \delta_{0}}(  \kappa | U=+\infty) \right)^{-1}    \right\rfloor.$$
In the rest of the proof, we write $m$ instead of $m(t,\epsilon)$ for the sake of readability.
Note that, in view of Proposition~\ref{p:moments-renouvellement}, for all $0 \leq \epsilon \leq \epsilon_{0}$,
$\E^{\epsilon,\theta}_{a \delta_{0}}(  \kappa | U=+\infty)   \leq C^{1/2}$, so that $m \geq 1$ as soon as 
$t \geq C^{1/2}(1-\lambda)^{-1}$, which does not depend on $\epsilon$.

We now re-use the random variables $a(t)$ defined in the proof of Proposition~\ref{p:speedup-biais}
above.  Using the fact that $t \mapsto \hat{r}_{t}$ is non-decreasing,
we see that $\hat{r}_{t} \geq \hat{r}_{\kappa_{a(t)}}$.
Moreover,
 $\hat{r}_{\kappa_{a(t)}}
 \geq  \hat{r}_{\kappa_{a(t)}}
  \un( a(t) \geq m    )  $,
and
$\hat{r}_{\kappa_{a(t)}}     \un( a(t) \geq m    )    
\geq  \hat{r}_{\kappa_{m}}   \un( a(t) \geq m) $ when $m \geq 1$.
Taking expectations, we deduce that, when $m \geq 1$,
\begin{equation}\label{e:minoration-rkappa}  \E^{\epsilon,\theta}_{\I_{0}}( t^{-1}  \hat{r}_{t}   )   \geq
   \E^{\epsilon,\theta}_{\I_{0}}( t^{-1} \hat{r}_{\kappa_{m}}) -
  \E^{\epsilon,\theta}_{\I_{0}}( t^{-1} \hat{r}_{\kappa_{m}}   \un( a(t) < m) ) .\end{equation}

Consider the first term in the r.h.s. of~(\ref{e:minoration-rkappa}) above, and observe that
$$ \E^{\epsilon,\theta}_{\I_{0}} (\hat{r}_{\kappa_{m}}) = \E^{\epsilon,\theta}_{\I_{0}}(\hat{r}_{\kappa}) +
  (m-1) \E^{\epsilon,\theta}_{a \delta_{0}}(  \hat{r}_{\kappa} | U=+\infty).$$
>From Proposition~\ref{p:moments-renouvellement},
 $\E^{\epsilon,\theta}_{\I_{0}}(\hat{r}_{\kappa}) \leq C^{1/2} $ for all $0 \leq \epsilon \leq \epsilon_{0}$.
 Moreover, from Identity~(\ref{e:caract-vitesse}), $\left( \E^{\epsilon,\theta}_{a \delta_{0}}(  \hat{r}_{\kappa} | U=+\infty) \right) \left( \E^{\epsilon,\theta}_{a \delta_{0}}(  \kappa | U=+\infty) \right)^{-1} = v_{\epsilon}$.
We easily deduce that, as $t$ goes to infinity, uniformly with respect to $0 \leq \epsilon \leq \epsilon_{0}$,
\begin{equation}\label{e:mainterm}\E^{\epsilon,\theta}_{\I_{0}} (\hat{r}_{\kappa_{m}}) =  (1-\lambda) t v_{\epsilon}+O(1).\end{equation}

Consider now the second term in the r.h.s. of~(\ref{e:minoration-rkappa}).
By Schwarz's inequality, 
\begin{equation}\label{e:holder} \E^{\epsilon,\theta}_{\I_{0}}( t^{-1} \hat{r}_{\kappa_{m}}   \un( a(t) < m) )
\leq \left(   \E^{\epsilon,\theta}_{\I_{0}} \left[ (t^{-1} \hat{r}_{\kappa_{m}})^{2} \right]       \right)^{1/2}
 \Q^{\epsilon, \theta}_{\I_{0}} (  a(t) < m)^{1/2}.\end{equation}
>From Proposition~\ref{p:moments-renouvellement} and~Corollary~\ref{c:iid},
it is easily checked that
\begin{equation}\label{e:moments-2} \E^{\epsilon,\theta}_{\I_{0}}\left[ \left(  \hat{r}_{\kappa_{m}} \right)^{2} \right] \leq C m^{2}.\end{equation}

On the other hand, one has that $\Q^{\epsilon, \theta}_{\I_{0}} (  a(t) < m) \leq  \Q^{\epsilon, \theta}_{\I_{0}} (  \kappa_{m} \geq t)$.
>From Proposition~\ref{p:moments-renouvellement} and~Corollary~\ref{c:iid}, the variance of
$\kappa_{m}$ with respect to $\Q^{\epsilon, \theta}_{\I_{0}}$ is bounded above by $Cm$, so that we can use the 
Bienaym\'e-Chebyshev's inequality to prove that, whenever
$t > \E^{\epsilon,\theta}_{\I_{0}}(\kappa_{m})$,
 \begin{equation}\label{e:bienaime}\Q^{\epsilon, \theta}_{\I_{0}} (  a(t) < m)  \leq
   Cm (     t-    \E^{\epsilon,\theta}_{\I_{0}}(\kappa_{m}))^{-2}.\end{equation}

Now,  using Proposition~\ref{p:moments-renouvellement} as in the proof of~(\ref{e:mainterm}) above, 
we can easily prove that, as $t$ goes to infinity, uniformly with respect to $0 \leq \epsilon \leq \epsilon_{0}$, 
$$\E^{\epsilon,\theta}_{\I_{0}}(\kappa_{m})=(1-\lambda)t+O(1).$$ 
Putting the above identity together with~(\ref{e:bienaime}), (\ref{e:moments-2}) and~(\ref{e:holder}), we deduce that, as $t$ goes to infinity,
uniformly  with respect to $0 \leq \epsilon \leq \epsilon_{0}$, 
$$ \E^{\epsilon,\theta}_{\I_{0}}( t^{-1} \hat{r}_{\kappa_{m}}   \un( a(t) < m) )
\leq  C      m^{3/2}  (   \lambda t^{2} + O(t)  )^{-1}.$$
In view of Lemma~\ref{l:mino-esp}, we have that $m \leq c^{-1} t$ for all $0 \leq \epsilon \leq \epsilon_{0}$, so we can conclude that, 
uniformly with respect to $0 \leq \epsilon \leq \epsilon_{0}$, 
\begin{equation}\label{e:errorterm}
\lim_{t \to +\infty}  \E^{\epsilon,\theta}_{\I_{0}}( t^{-1} \hat{r}_{\kappa_{m}}   \un( a(t) < m) ) = 0.
\end{equation}

Plugging~(\ref{e:mainterm}) and~(\ref{e:errorterm}) in~(\ref{e:minoration-rkappa}), we finally deduce that, as
$t$ goes to infinity, uniformly with respect to $0 \leq \epsilon \leq \epsilon_{0}$, 
$$ \E^{\epsilon,\theta}_{\I_{0}}( t^{-1}  \hat{r}_{t}   )  \geq   (1-\lambda) v_{\epsilon}  + o(1).$$
The conclusion of the Proposition follows by noting that, since $v_{\epsilon} \leq v_{\epsilon_{0}}$,
$  (1-\lambda) v_{\epsilon} \geq v_{\epsilon} - \lambda v_{\epsilon_{0}}$.

\end{proof}

\begin{lemma}\label{l:cv-dom}
For all $t \geq 0$,
$$\lim_{\epsilon \to 0+} \E^{\epsilon, \theta}_{\I_{0}} (\hat{r}^{\epsilon}_{t}) = \E^{\epsilon,\theta}_{\I_{0}}(\hat{r}^{0}_{t}).$$
\end{lemma}

\begin{proof}
Consider a given $t \geq 0$. By Proposition~\ref{p:negligeabilite} in Appendix~\ref{s:negligeabilite},
with $\P$ probability one, we can find a (random) $K \leq 0$ such that
$\sup \{ F^{\epsilon_{0}}_{s}(x,i)  ; \, 0 \leq s \leq t, \, x < K, \, 1 \leq i \leq a \} \leq 0$,
so that
$\sup \{ F^{\epsilon}_{s}(x,i)  ; \, 0 \leq s \leq t, \, x < K, \, 1 \leq i \leq a \} \leq 0$
for all $0 \leq \epsilon \leq \epsilon_{0}$.
As a consequence, for all $0 \leq \epsilon \leq \epsilon_{0}$, with probability one,
$r_{t}^{\epsilon}(\I_{0}) = r_{t}^{\epsilon}(w(K))_{s \geq 0})$,
where $w(K)$ is the configuration defined by
$A=\{K, \ldots, 0 \} \times \{1, \ldots, a\}$, $r=0$ and $F(x,i)=x$ for all $(x,i) \in A$.

Since, for every $0 \leq \epsilon \leq \epsilon_{0}$,
with probability one $r^{\epsilon}_{t} \leq r^{\epsilon_{0}}_{t}$, we see that the value of
$r^{\epsilon}_{t}$ is entirely determined by the trajectories up to time $t$ of the random walks born
at sites $(x,i)$ with $K \leq x \leq r^{\epsilon_{0}}_{t}$.
With probability one again, we are dealing with a finite number of random walks, and a finite number of steps.
We now see that, for all $\epsilon$ small enough, these trajectories are identical to what they are for $\epsilon=0$,
so that $r^{\epsilon}_{t}=r^{0}_{t}$. Since $0 \leq r^{\epsilon}_{t} \leq r^{\epsilon_{0}}_{t}$ and $r^{\epsilon_{0}}_{t}$ is integrable w.r.t. $\P$,
we can use the dominated
convergence theorem to deduce the conclusion.
\end{proof}

\begin{proof}[Proof of Proposition~\ref{p:cv-vitesse}]

Let $\zeta>0$, and, following Proposition~\ref{p:limite-inf-uniforme},
consider a $t_{\zeta}$ such that, for
 all $t \geq t_{\zeta}$ and all $0 \leq \epsilon \leq \epsilon_{0}$,
$$v_{\epsilon} \leq \E^{\epsilon,\theta}_{\I_{0}} (t^{-1} \hat{r}_{t}) + \zeta.$$
Consider now, thanks to Proposition~\ref{p:lgn-position}, a $t \geq t_{\zeta}$ such that
$ \E^{0,\theta}_{\I_{0}} (t^{-1} \hat{r}_{t}) \leq v+ \zeta$.
Now, thanks to Lemma~\ref{l:cv-dom}, we know that, for all $\epsilon$ small enough,
$$\E^{\epsilon,\theta}_{\I_{0}} (t^{-1} \hat{r}_{t})  \leq \E^{0,\theta}_{\I_{0}} (t^{-1} \hat{r}_{t})+\zeta.$$
Putting together the above inequalities, we deduce that, for all $\epsilon$ small enough,
$v_{\epsilon} \leq v + 3\zeta.$
Since $v_{\epsilon} \geq v$, the conclusion follows.

\end{proof}

Now Proposition~\ref{p:speedup} follows from Proposition~\ref{p:speedup-biais}, Proposition~\ref{p:compare-loi} and
Proposition~\ref{p:cv-vitesse}, as explained in the beginning of this section.

\section{Slowdown large deviations}\label{s:slowdown}




For $x\geq 1$ and $t \geq 0$, let $(\zeta_{t})_{t \geq 0}$ denote a
continuous time  simple symmetric random walk  starting form $0$ of
total jump rate $2$.
Let
$$\bar G_{t}(x) := P \left(   \sup_{s \in [0,t)}  \zeta_{s}  < x  \right),    \   
G_{t}(x) :=   P \left( \zeta_{t}  \geq x  \right) .$$
In the sequel we will use the fact that for fixed $t\ge 0$, $G_t(\cdot)$ is non-decreasing and 
$\bar G_t(\cdot)$ is non-increasing, and that, thanks to the reflection principle, 
\begin{equation}
\label{reflection}
1-\bar G_{t}(x) = 2  G_{t}(x) - P\left(  \zeta_{t}  = x \right).
\end{equation}

\subsection{Proof of Theorem~\ref{t:slowdown} (a) and (c)}

We start with the proof of Theorem~\ref{t:slowdown} (c).
 The fact that $r_{t}=0$ means that no particle in the initial configuration hits $1$ before time $t$. Both the upper and lower bounds can then be understood 
heuristically as follows. Since we consider simple symmetric random walks,   for large $t$,
the constraint of not hitting $1$ before time $t$ has a cost only for particles
within a distance of order $t^{1/2}$ of the origin. Now these particles perform
 independent random walks, and their  number has an order of magnitude lying between 
$t^{u(\eta_{w})/2}$ and $t^{U(\eta_{w})/2}$.


We start with the lower bound.
When $U=+\infty$, the inequality holds trivially, so we assume in the sequel that
 $U<+\infty$.
The event $t^{-1} r_{t}(w) = 0$,
 implies that  none of the random walks
corresponding to particles in the initial condition $w$ hit $1$ before time $t$. By independence
of the random walks, the corresponding probability equals
$$ \prod_{x=0}^{-\infty}  \bar G_{t}(-x+1)^{\eta_{w}(x)}.$$
Now let $b_{1}>0$ be such that $1 - 2s \geq \exp\left(  - 4s  \right)$ for all $0 \leq s \leq b_{1}$.
From~(\ref{reflection}), we see that
for any $t\ge 0$ and $y\leq 0$, $\bar G_{t}(-y+1) \geq 1 - 2 G_{t}(-y+1)$.
By the central limit theorem,  we can find $t_{0}$ and $K>0$ such that, for all $t \geq t_{0}$ and
 $y \leq -Kt^{1/2}$, $G_{t}(-y+1) \leq b_{1}$. 
Let $k_{t} :=\ceil{Kt^{1/2}}$.
Then, for all $t \geq t_{0}$

$$ \prod_{x=-k_{t}}^{-\infty}  \bar G_{t}(-x+1)^{\eta_{w}(x)} \geq   \exp \left( - 4\sum_{x=-k_{t}}^{-\infty}  \eta_{w}(x)  G_{t}(-x+1)  \right).$$
Now, by definition of $G_{t}$,
 \begin{eqnarray*}\sum_{x=0}^{-\infty}  \eta_{w}(x)  G_{t}(-x+1)  =& E \left( \sum_{x=0}^{-\infty}  \eta_{w}(x) \un( \zeta_{t} \geq -x+1  )
    \right)\\
          =   E \left[  \un(  \zeta_{t} \geq 1   ) \left(        \sum_{x=0}^{- \zeta_{t} +1}  \eta_{w}(x)    \right)  \right]  
          =&
           E \left[ \un(  \zeta_{t} \geq 1   ) ( H_w(  - \zeta_{t}  +1  )  )  \right].  \end{eqnarray*}
         By assumption, $H_{w}(x) \leq |x|^{U+o(1)}$. H\"older's inequality  yields that
$$E \left[ \un(  \zeta_{t} \geq 1  ) ( H_w(  - \zeta_{t}  +1)  ) \right] \leq  t^{U/2+o(1)}.$$
 We deduce that for all $t \geq t_{0}$
 \begin{equation}\label{e:unbout} \prod_{x=-k_{t}}^{-\infty}  
\bar G_{t}(x)^{\eta_{w}(x)} \geq  \exp(   -  t^{U/2+o(1)} ).\end{equation}
Now, for $-k_{t} < y \leq 0$, observe that $\bar G_{t}(-y+1) \geq \bar G_{t}(1)$.
As a consequence,
  $$ \prod_{x=0}^{-k_{t}+1} \bar G_{t}(-x+1)^{\eta_{w}(x)}  \geq 
\bar G_{t}(1)^{H_{w}(-k_{t}+1)}.$$


\noindent But there exists $c_{4}>0$, such that, for large enough $t$, 
$\bar G_{t}(1)\geq c_{4} t^{-1/2}$.
Using again the fact that $H_{w}(x) \leq |x|^{U+o(1)}$, it is easy to deduce that
$\bar G_{t}(1)^{H_{w}(-k_{t}+1)} \geq    \exp(   -  t^{U/2+o(1)} )$, whence
 \begin{equation}\label{e:unautrebout} \prod_{x=0}^{-k_{t}+1} \bar G_{t}(-x+1)^{\eta_{w}(x)}  \geq   \exp(   -  t^{U/2+o(1)} ). \end{equation}
From~(\ref{e:unbout}) and~(\ref{e:unautrebout}), we deduce that
 \begin{equation*}\P(t^{-1} r^{0}_{t}(w)  \leq 0  ) \geq   \exp(   -  t^{U/2+o(1)} ).\end{equation*}
 Now, let us prove the upper bound when $u<+\infty$. Using an argument similar to the one used in the proof of the lower bound above, we easily obtain that 
 $$\P(t^{-1} r_{t}(w)  = 0  ) \leq   \exp\left( - E \left[ \un(  \zeta_{t} \geq 1  ) ( H_w(  - \zeta_{t}  +  1 ) ) \right]   \right).$$
It is then easy to deduce that 
$$E \left[ \un(  \zeta_{t} \geq 1  ) ( H_w(  - \zeta_{t}  +  1  )) \right] \geq  t^{u/2+o(1)},$$
and the upper bound follows.

We now turn to the proof of Theorem~\ref{t:slowdown} (a). 
The idea of the proof when $s(\eta)=1/2$ is to combine the following two arguments. First, for $b>0$, it costs nothing to prevent all the particles 
in the initial condition from hitting $\floor{bt}$ up to time $t$. Intuitively, this result comes from the fact that 
hitting $\floor{bt}$ before time $t$ has an exponential cost for any particle in the initial condition within distance $O(t)$ of the origin, and, due 
to~{\bf(G)}, there is a subexponentially large number of such particles. 

Second, in the worst case 
where all the particles attached to sites $1 \leq x \leq bt$ are turned into $X$ particles instantaneously at time zero, the cost of  preventing all 
these particles from hitting $bt$ up to time $t$ is of order $\exp(-t^{1/2+o(1)})$, due to the lower bound in~(\ref{e:slowdown-zero-bound}) proved above, 
The actual proof is in fact more complex since we want to
consider probabilities of the form $\P( c t \leq   r_{t}  \leq b t )$, and not only $\P(r_{t} \leq b t )$, and deal also with the case $s(\eta)<1/2$.

We state two lemmas  before giving the proof.

\begin{lemma}\label{l:atteinte-positif-fois-t}
Consider an initial condition  $w=(F,0,A)$  satisfying assumption {\bf (G)}. Then, for all $b>0$, and all $\varphi>0$,


$$\P \left[  \max_{(x,i) \in A} \sup_{0\le s\le t} F_{s}(x,i)    \geq bt  \right]  \leq 
  f_{\varphi}(w) \exp \left[ t (\cosh(2 \varphi) - 1 ) \right]  G_{t}(\floor{bt})^{1/2}.$$
\end{lemma}

\begin{proof}
The probability we are looking at is the probability that at least one of the random walks
corresponding to particles in $w$ exceeds $bt$ before time $t$.
By the union bound, this probability is smaller than
$$\sum_{x=0}^{-\infty}  \eta_{w}(x) (1-  \bar G_{t}(-x+\floor{bt}) ) 
 \leq  \sum_{x=0}^{-\infty}  2 \eta_{w}(x) G_{t}(-x+\floor{bt}) ).$$
Now observe that by definition of $G_{t}$,
 \begin{eqnarray*}\sum_{x=0}^{-\infty}  \eta_{w}(x)  G_{t}(-x+   \floor{bt} )  &=& E \left( \sum_{x=0}^{-\infty}  \eta_{w}(x) \un( \zeta_{t} \geq -x+ \floor{bt} )
    \right)\\
          =   E \left[  \un(  \zeta_{t} \geq  \floor{bt}   ) \left(        \sum_{x=0}^{- \zeta_{t} + \floor{bt}}  \eta_{w}(x)    \right)  \right]  
          &=&
           E \left[ \un(  \zeta_{t} \geq  \floor{bt}   ) ( H_w(  - \zeta_{t}  + \floor{bt}  )  )  \right].  \end{eqnarray*}
>From assumption {\bf (G)}, we deduce that, for all $\varphi>0$,  
$H_w(x) \leq f_{\varphi}(w) \exp(-\varphi x ) $ for all $x \leq 0$.
As a consequence, when $\zeta_{t} \geq  \floor{bt}$,  $H_{w}( -\zeta_{t} + \floor{bt}   )  \leq 
 H_{w}( -\zeta_{t}) \leq    f_{\varphi}(w) \exp( \varphi \zeta_{t} )$.
Applying Schwarz's inequality, we see that
$$  E \left[ \un(  \zeta_{t} \geq  \floor{bt}   ) ( H_w(  - \zeta_{t}  + \floor{bt}  )  )  \right] \leq
 \P\left( \zeta_{t} \geq  \floor{bt} \right)^{1/2}   f_{\varphi}(w)  E \left[    \exp( 2\varphi \zeta_{t} )    \right]^{1/2}.$$
Now note that
$ E \left[    \exp( 2\varphi \zeta_{t} )    \right] = \exp [2(\cosh (2 \varphi) - 1) t ]$.
\end{proof}

\begin{lemma}\label{l:croissance-cond-init}
Consider an initial condition  $w=(F,0,A)$ satisfying assumption {\bf (G)}.
 Then, for all $\varphi>0$
$$\E \left[ \sum_{(x,i) \in A_{t}}   \exp \left(  \varphi( F_{t}(x,i)- r_{t})  \right)  \right]  \leq 
    \exp [2(\cosh (\varphi) - 1) t ]  f_{\varphi}(w)   + a \E(r_{t}).$$
\end{lemma}


\begin{proof}
Write $\sum_{(x,i) \in A_{t}} = \sum_{(x,i) \in A} + \sum_{(x,i) \in  A_{t} \setminus A}$. 
For $(x,i) \in A$, observe that $\exp\left( \varphi (F_{t}(x,i) - r_{t}) \right)   \leq  \exp\left( \varphi F_{t}(x,i) \right)$
and that $ \E \left[    \exp( \varphi F_{t}(x,i) )    \right] = \exp [\varphi x+2(\cosh ( \varphi) - 1) t ]$.
As a consequence, 
\begin{equation}\label{e:gauchedezero}   \E \left[ \sum_{(x,i) \in A}   \exp \left(  \varphi( F_{t}(x,i)- r_{t})  \right)  \right]  \leq  
 \exp [2(\cosh (\varphi) - 1) t ] f_{\varphi}(w).\end{equation}
On the other hand, observe that  $A_{t} \setminus A  =  \{1, \ldots, r_{t} \} \times \{ 1,\ldots, a\}$. 
Since it is always true that $F_{t}(x,i) \leq r_{t}$,
\begin{equation}\label{e:droitedezero}\sum_{(x,i) \in A_{t} \setminus A}   \exp \left(  \varphi( F_{t}(x,i)- r_{t})  \right) \leq  a  r_{t}.\end{equation}
The result follows from putting together~(\ref{e:gauchedezero}) and~(\ref{e:droitedezero}). 
\end{proof}

\begin{proof}[Proof of Theorem~\ref{t:slowdown}] 


 Let $\alpha,\delta\in (0,1)$ be such that $c<v(1-\alpha)<b$,
$c<(1-\alpha) (1-\delta)v<(1-\alpha) (1+\delta)v<b$, and define
 $\gamma:=b-(1-\alpha) (1+\delta)v$.
For each $t>0$,  define the events 

\begin{eqnarray*}B_{t} &:=& \left\{ v(1-\alpha) (1-\delta)t   \leq  r_{(1-\alpha)t} \leq v(1-\alpha) (1+\delta)t \right\},\\
C_{t}&:=&    \left\{ \max_{  (x,i) \in A_{t(1-\alpha)}, \,  t(1-\alpha) \leq s  \leq t}   F_{s}(x,i)  \leq  r_{t(1-\alpha)} + \gamma t  \right\},\\
D_{t} &:=&  \left\{   \max_{r_{(1-\alpha)t} < x \leq bt, \, 1 \leq i \leq a , \, 0 \leq s  \leq  \alpha t}   x+Y_{x,i,s}      \leq bt \right\}. \end{eqnarray*}
Observe that 

\begin{equation}\label{e:inclusion-ralentissement}B_{t} \cap C_{t} \cap D_{t} \subset \{ct \leq r_{t} \leq bt \}.\end{equation}
Indeed, thanks to the choice of $\delta$, $B_{t}$ implies that
$r_{(1-\alpha)t} \geq ct$, so that $r_{t} \geq ct$.
 On the other hand, since $r_{(1-\alpha)t}<bt$ on $B_{t}$, the event $B_{t} \cap \{ r_{t} > bt \}$ 
implies that either a particle born before time $t(1-\alpha)$ at a position $x \leq r_{t(1-\alpha)}$, or 
a particle born  between time $(1-\alpha) t$ and $t$ at a position $r_{(1-\alpha)t}  < x < bt$, exceeds $bt$ at a time between $t(1-\alpha)$ and $t$. 
The former possibility is ruled out by $B_{t} \cap C_{t}$,  since on $B_{t} \cap C_{t}$, $r_{t} \leq r_{(1-\alpha)t}+\gamma t \leq bt$.
The latter possibility is ruled out by $D_{t}$. 
Now define
$$l(t):=\exp[2(\cosh (\varphi)-1) (1-\alpha) t ] f_{\varphi}(w)   + a \E(r_{(1-\alpha)t}),$$
and
$$H_{t} := \left\{     \sum_{(x,i) \in A_{(1-\alpha)t}}   \exp \left(  \varphi( F_{(1-\alpha)t}(x,i)- r_{(1-\alpha)t})  \right)  \leq 
     2 l(t)  \right\}.$$
By  Lemma~\ref{l:croissance-cond-init} and  Markov's inequality, for all $t \geq 0$,
$\P\left(  H_{t} \right)  \geq 1/2.$
Moreover, by the law of large numbers~(\ref{e:lgn-presque-sure}), $\lim_{t \to +\infty}\P( B_{t}) = 1.$     
 We deduce that there exists a $t_{0}$ such that, for all  $t \geq t_{0}$, $\P(  B_{t} \cap H_{t}  ) \geq 1/4$.
Let us call $\F_{t}$ the  $\sigma-$algebra generated by the history of the
particle system up to time $t$.
 Observe that $B_{t}$ and $H_{t}$ belong to $\F_{(1-\alpha)t}$, 
and by Lemma~\ref{l:atteinte-positif-fois-t}, 
 on $H_{t}$, 

$$\P(  C_{t}^{c} | \F_{(1-\alpha)t})  \leq 2 l(t) \exp \left[ \alpha t (\cosh(2 \varphi) - 1 ) \right] G_{\alpha t}(\floor{\gamma t})^{1/2}.$$
We deduce that

\begin{equation} \label{e:partie-negligeable} 
\P(  B_{t} \cap H_{t} \cap C_{t}^{c} ) \leq   
2l(t) \exp \left[ \alpha t (\cosh(2 \varphi) - 1 ) \right]  G_{\alpha t}(\floor{\gamma t})^{1/2}.\end{equation}
Moreover,
$$\P( D_{t} | \F_{(1-\alpha)t}  ) \geq \P(   r_{\alpha t}(\I_{0}) = 0 ),$$
so that
\begin{equation}\label{e:partie-principale} \P(  B_{t} \cap H_{t} \cap D_{t} ) \geq   (1/4)   \P(   r_{\alpha t}(\I_{0}) = 0 ) \geq \exp(-t^{1/2+o(1)}),\end{equation}
where the last inequality is due to the lower bound in~(\ref{e:slowdown-zero-bound}).
By standard large deviations bounds for the simple random walk, there exists $\zeta(\alpha,\gamma)>0$ depending only on $\gamma$ and $\alpha$
such that, as $t$ goes to infinity, $\liminf_{t \to +\infty} t^{-1} \log G_{\alpha t}(\floor{\gamma t})  = -\zeta(\alpha,\gamma)$.
Furthermore,
$\lim_{t \to +\infty}  t^{-1} \log( 2l(t) \exp \left[ \alpha t (\cosh(2 \varphi) - 1 ) \right] ) = \xi(\alpha,\varphi)$, 
where $\xi(\alpha,\varphi) :=  \alpha  (\cosh(2 \varphi) - 1 ) +
2(\cosh (\varphi) - 1) (1-\alpha)$.
We see that, choosing $\varphi$ small enough, $\xi(\alpha,\varphi)<\zeta(\alpha,\gamma)/2$.
For such a $\varphi$, (\ref{e:partie-negligeable}) and~(\ref{e:partie-principale})
show that  $  \P(  B_{t} \cap H_{t} \cap C_{t}^{c} ) = o(  \P(  B_{t} \cap H_{t} \cap D_{t} ) )$, 
so that 
$\P(    B_{t} \cap H_{t} \cap D_{t} \cap C_{t}   ) \geq \exp(-t^{1/2+o(1)})$.
It then follows from~(\ref{e:inclusion-ralentissement}) that $\P(  ct \leq r_{t} \leq bt) \geq \exp(-t^{1/2+o(1)})$,
so we are done when $s(\eta)=1/2$.

Now, let us choose 
an $(x,i)\in A$ for the initial condition $w=(F,0,A)$. Define also $\tau = \inf \{  s \geq 0; \, F_{s}(x,i) = 0 \}$.
Let 
\begin{eqnarray*} K_{t}&:=&\{  (1- (1-\alpha )(1+\delta))t   \leq \tau \leq   (1- (1-\alpha )(1-\delta))t  \},\\
                             L_{t}&:=&\{   ct  \leq  r_{(1-\alpha)(1-\delta)t+\tau}(w(x,i)) \leq   r_{(1-\alpha)(1+\delta)t+\tau}(w(x,i)) \leq  bt  \}, \\
                             L'_{t}&:=&  \{   ct  \leq  r_{(1-\alpha)(1-\delta)t}(\delta_{0}) \leq   r_{(1-\alpha)(1+\delta)t}(\delta_{0}) \leq  bt  \}, \\
                             M_{t}&:=&  \{         \mbox{ for all $(y,j) \in A\setminus \{ (x,i) \}$ and  all $0 \leq s \leq t$,  }F_{s}(y,j) \leq 0    \}.
\end{eqnarray*}
Observe that, on $M_{t}$, $r_{t}(w)=r_{t}(w(x,i))$. Moreover,  $K_{t} \cap L_{t} \subset \{ ct    \leq r_{t}(w(x,i))  \leq bt  \}$.  
As a consequence, \begin{equation}\label{e:inclusion-slowdown-bis}M_{t} \cap K_{t} \cap L_{t}  \subset  \{ ct    \leq r_{t}(w)  \leq bt  \}  .\end{equation}
But according to the lower bound of Theorem~\ref{t:slowdown} (b), 
$\P(  M_{t}  ) \geq \exp(-t^{U(\eta_{w})/2 + o(1)})$.
 On the other hand, conditional upon $\tau$, $r_{s+\tau}(w(x,i))$ has the (unconditional) distribution of
$r_{s}(\delta_{0})$, for all $s \geq 0$. As a consequence, $\P(K_{t} \cap L_{t})=\P(K_{t}) \P(L'_{t})$, and, by
the law of large numbers~(\ref{e:lgn-presque-sure}),  $\lim_{t \to +\infty} \P(L'_{t}) = 1$. 
Moreover, it is easily seen from elementary estimates on hitting times by a simple symmetric continuous time random walk that 
 $  \liminf_{t \to +\infty} t^{-1/2}  \P(  K_{t}  ) > 0$.
Finally, $M_{t}$ being defined in terms of random walks that do not enter the definition of $K_{t}$ and $L_{t}$, we deduce that
$M_{t}$ is independent from $K_{t} \cap L_{t}$. 
We finally deduce that $\P(M_{t} \cap K_{t} \cap L_{t}) \geq  \exp(-t^{U(\eta_{w})/2 + o(1)})$, 
and the result follows from~(\ref{e:inclusion-slowdown-bis}).

\end{proof}

\subsection{Proof of Theorem~\ref{t:optimal}}

As for the upper bound in~(\ref{e:slowdown-zero-bound}),
we easily obtain that 
 $$\P(t^{-1} r_{t}(w)  \leq bt  ) \leq   \exp\left( - E \left[ \un(  \zeta_{t} \geq \ceil{bt}  )  H_w(  - \zeta_{t}  +  \ceil{bt}  ) \right]   \right).$$
 It is easily checked that, for small enough $b>0$, 
\begin{equation}\label{e:explosion}\liminf_{t \to +\infty} t^{-1}\log E \left[ \un(  \zeta_{t} \geq \ceil{bt}  )  \exp \left(  \theta (\zeta_{t}  -  \ceil{bt}) \right)     \right] > 0.\end{equation}
This proves $(i)$.
We now prove $(ii)$. Again, it is easily checked that, for all $b>0$, there exists $\theta>0$ such that~(\ref{e:explosion}) holds.  
Choosing $b>v$, the result follows.

\subsection{Proof of Theorem~\ref{t:slowdown} (b)}

Note that by coupling, it is enough to prove the result with an initial condition consisting of exactly $a$ particles per site $x \leq 0$, that is, 
with $w=\I_{0}$.
Hence, we will establish that
for all $0 < b <v$, and all $\alpha>0$,
$$
\liminf_{t\to\infty}(\log t)^{-1} \left| \log \P\left[\frac{r^{0}_t(\I_{0})}{t}\leq b\right] \right| \geq 1/3.
$$
Using the fact that
$\P(T^{0}_{\I_{0}}(\floor{bt}) \geq t)  \leq   \P(r^{0}_{t}(\I_{0}) \leq bt) \leq \P(T^{0}_{\I_{0}}(\ceil{bt}) \geq t)$,
it is easy to see that~(\ref{e:slowdown-upper-bound}) is  equivalent to the  following
 result.
\begin{prop}\label{p:bornesubadd}
For every $c>v^{-1}$, as $n$ goes to infinity,
$$   \P(T^{0}_{\I_{0}}(n) \geq c n)   \leq \exp \left(-n^{1/3+o(1)}  \right).$$
\end{prop}

\noindent Our strategy for proving Proposition~\ref{p:bornesubadd} can be sketched as follows.
By subadditivity, for all  $m \geq 1$
$$T^{0}_{\I_{0}}(n) \leq \sum_{j=0}^{\floor{n/m}} T^{0}_{\I_{mj}}(m(j+1)),$$
so that
\begin{equation}\label{e:compare-par-sousadd}  \P( T^{0}_{\I_{0}}(n)  \geq cn ) \leq
\P\left( \sum_{j=0}^{\floor{n/m}} T^{0}_{\I_{mj}}(m(j+1)) \geq  (mc)  \floor{n/m} \right).\end{equation}
Now, by translation invariance, for all $j \geq 0$, $T^{0}_{\I_{mj}}(m(j+1))$ and
$T^{0}_{\I_{0}}(m)$ have the same distribution, and it can be shown that
$$\lim_{m \to +\infty} m^{-1} \E(T^{0}_{\I_{0}}(m)) = v^{-1}.$$
Hence, given $c>v^{-1}$ we can always find $m \geq 1$ such that $mc >  \E(T^{0}_{\I_{0}}(m))$,
so that the r.h.s. of~(\ref{e:compare-par-sousadd}) is the probability of a large deviation above the mean for the sum
$  \sum_{j=0}^{\floor{n/m}} T^{0}_{\I_{mj}}(m(j+1))$.
We then seek to apply large deviations bounds for i.i.d. variables in order to estimate this probability.
Of course, the random variables $\left\{ T^{0}_{\I_{mj}}(m(j+1)) \right\}_{j \geq 0}$ are {\it not} independent, but the dependency
between  $\left\{ T^{0}_{\I_{mj}}(m(j+1)) \right\}_{j \leq j_{1}}$ and  $\left\{ T^{0}_{\I_{mj}}(m(j+1)) \right\}_{j \geq j_{2}}$ is weak when $j_{2}-j_{1}$ is large.
Indeed, for given $j$, $T^{0}_{\I_{mj}}(m(j+1))$ mostly depends on the behavior of the random walks born at sites close to $mj$.
We implement this idea by using a technique already exploited in~\cite{RamSid} in a similar context. Given $\ell \geq 1$, we
define a family $\left\{T^{'}_{\I_{mj}}(m(j+1)) \right\}_{j \geq 0}$ of hitting times as follows: $T^{'}_{\I_{mj}}(m(j+1))$ uses the same
random walks as  $T^{0}_{\I_{mj}}(m(j+1))$ for particles born at sites $(x,i)$ with $mj-m\ell<x < m(j+1)$, but uses fresh independent random walks
for particles  born at sites $(x,i)$ with $x \leq mj - m\ell$. We can then prove that the following  properties hold:
\begin{itemize}
\item[(a)] For all $j \geq 0$, the family $\left\{   T^{'}_{\I_{mj + p m(\ell+1) }}(mj+ p m (\ell+1) + m) \right\}_{p \geq 0}$ is i.i.d.;
\item[(b)] when $\ell$ is large, the probability that $T^{'}_{\I_{mj}}(m(j+1))   = T^{0}_{\I_{mj}}(m(j+1))$ is close to $1$.
\end{itemize}
We can thus obtain estimates on the  r.h.s. of~(\ref{e:compare-par-sousadd}) by estimating separately the probability that
 $ T^{'}_{\I_{mj}}(m(j+1))   = T^{0}_{\I_{mj}}(m(j+1)) $
for all $0 \leq j \leq \floor{m/n}$, and the probability that
$\sum_{j=0}^{\floor{n/m}} T^{'}_{\I_{mj}}(m(j+1)) \geq  (mc)  \floor{n/m}$. Now, thanks to property (a) above,
this last sum can be split evenly into $\ell+1$ subsums of i.i.d. random variables distributed
 as $T^{0}_{\I_{0}}(m)$. Controlling the tail of $T^{0}_{\I_{0}}(m)$ then allows
us to apply large deviation bounds for i.i.d. variables separately to each of 
these subsums. In fact, the proof of~(\ref{e:slowdown-upper-bound}) is a bit more subtle, since it also makes use of
a positive association property, but we do not go into the details here (see Remark~\ref{r:pourquoi-assoc} below).

\subsection{Proof of Proposition~\ref{p:bornesubadd}}

Observe that, since for all $u \in \Z$ and $k \geq 0$, $\I_{u} \oplus k = \I_{u+k}$,
the subadditivity property  (part $(iii)$) of Proposition~\ref{p:charact-T}
reads as: $$\mbox{ for all $n, m \geq 0$, }T^{0}_{\I_{0}}(n+m) \leq  T^{0}_{\I_{0}}(n) +    T^{0}_{\I_{n}}(m).$$
Now, let $c>v^{-1}$.
Thanks to subadditivity,  for all $m \geq 1$ we have that
$$   \P(T^{0}_{\I_{0}}(n) \geq c n )  \leq    \P\left(  \sum_{j=0}^{\floor{n/m}} T^{0}_{\I_{mj}}(m(j+1)) \geq cn \right).$$

In Steps 1 and 2 below, $m$ and $\ell$ denote fixed positive integers, while $\alpha$ denotes a fixed real number $0 < \alpha < 1$.
For the sake of readability, the dependence with respect to these numbers is usually not mentioned explicitly in the notations.
Only in Step 3 have the values of $m, \ell$ and $\alpha$ to be specified.  

\subsubsection{Step 1: Comparison with a sum of i.i.d. random variables}

Assume that the probability space $(\Omega, \F, \P)$ is such that we have access to an i.i.d. family of random variables
$$\left[ (\tau^{j}_{k}(u,i), U^{j}_{k}(u,i))  ; \,  j \geq 0, \, k \geq 1, \, u \in \Z , \, 1 \leq i \leq a           \right],$$
independent from $$\left[ (\tau_{k}(u,i), U_{k}(u,i))  ; \, \, k \geq 1, \, u \in \Z , \, 1 \leq i \leq a           \right],$$ and
 such that,
 for all $(j,k,u,i)$,  $\tau^{j}_{k}(u,i)$ has an exponential(2) distribution, $U^{j}_{k}(u,i)$ has the uniform distribution on $(0,1)$,
 and   $\tau^{j}_{k}(u,i)$  and $U^{j}_{k}(u,i)$
are independent.

Let
 $$ \varepsilon^{j}_{n}(x,i):= 2 ( \un(  U^{j}_{n}(x,i)
 \leq 1/2 )  )  - 1. $$
Now, for all $\ell \geq 1$ and $j \in  \Z$, define, for all $u,v \in \Z$ such that $u<v$, and $1 \leq i \leq a$,
$$\B_{j}(u,i,v) := \inf \left\{ \sum_{k=1}^{m} \tau^{j}_{k}(u,i);  \,  u + \sum_{k=1}^{m} \varepsilon^{j}_{k}(u,i)  = v ,  \, m \geq 1\right\},$$
and let
$$ \C_{j}(u,i,v)   :=  \left\{   \begin{array}{l}   \B_{j}(u,i,v)  \mbox{ if }  u \leq  mj - m\ell \\
  \A^{0}(u,i,v) \mbox{ if } u > mj - m\ell
   \end{array} \right.  ,$$
   where $\A^0$ is defined in display (\ref{def-a}) in Section~\ref{s:construction}.
Then let
  \begin{equation*}  T^{'}_{\I_{mj}}(m(j+1)) := \inf \sum_{g=1}^{L-1}  \C_{j}(x_{g},i_{g},x_{g+1}) ,\end{equation*}
 where the infimum is taken over all finite sequences with $L \geq 2$, $x_{1},\ldots,x_{L} \in \mathbb{Z}$ and $i_{1},\ldots, i_{L-1}$ such that
 $x_{1} \leq mj$,  $ mj<x_{2} < \cdots < x_{L-1}<m(j+1), \, x_{L}=m(j+1)$,
 $i_{1},i_{2},\ldots, i_{L-1} \in \{1, \ldots, a \}$.  
Clearly, $T^{'}_{\I_{mj}}(m(j+1))$ and $T^{0}_{\I_{mj}}(m(j+1))$ have the same distribution. Moreover, we have the
 following lemma, whose proof  is immediate.

\begin{lemma}\label{l:tempsindep}
For every $j \in \Z$, the family of random variables
$$\left(   T^{'}_{\I{ mj+ pm(\ell+1)   }}( mj+ p m(\ell+1)+m ) \right)_{p \in \Z}  $$ is i.i.d.
\end{lemma}


We now study the event $ \left\{ T^{'}_{\I_{mj}}(m(j+1)) =  T^{0}_{\I_{mj}}(m(j+1)) \right\}$. To this end, 
let
\begin{equation*}  J_{j} := \inf \sum_{g=1}^{L-1}  \C_{j}(x_{g},i_{g},x_{g+1}) , \,
  K_{j} := \inf \sum_{g=1}^{L-1}  \A^{0}_{j}(x_{g},i_{g},x_{g+1}) ,\end{equation*}
 where in both cases the infimum is taken over all finite sequences with $L \geq 2$, $x_{1},\ldots,x_{L} \in \mathbb{Z}$ and $i_{1},\ldots, i_{L-1}$ such that
 $x_{1} \leq mj- m \ell$,  $ mj<x_{2} < \cdots < x_{L-1}<m(j+1), \, x_{L}=m(j+1)$,
 $i_{1},i_{2},\ldots, i_{L-1} \in \{1, \ldots, a \}$.
 Let also
\begin{equation*}  L_{j} := \inf \sum_{g=1}^{L-1}  \A^{0}_{j}(x_{g},i_{g},x_{g+1}) ,\end{equation*} where the infimum is taken over all finite sequences with $L \geq 2$, $x_{1},\ldots,x_{L} \in \mathbb{Z}$ and $i_{1},\ldots, i_{L-1}$ such that
$  mj- m \ell <  x_{1} \leq mj$,  $ mj<x_{2} < \cdots < x_{L-1}<m(j+1), \, x_{L}=m(j+1)$,
$i_{1},i_{2},\ldots, i_{L-1} \in \{1, \ldots, a \}$.

 Observe that, $ T^{'}_{\I_{mj}}(m(j+1))  = \min(J_{j},L_{j})$ and that $ T^{'}_{\I_{mj}}(m(j+1))  = \min(K_{j},L_{j})$. As a consequence, 
$$\{ \min( J_{j} , K_{j} )  \geq        L_{j}     \}   \subset    \left\{ T^{'}_{\I_{mj}}(m(j+1)) =  T^{0}_{\I_{mj}}(m(j+1))    \right\}.$$
For $\alpha>0$, we now define  $$D(j):= \left\{    \min( J_{j} , K_{j} )  <   \alpha (m\ell)^{2} \right\},$$
and  $$F(j) := \{   L_{j} \geq  \alpha  (m\ell)^{2} \},$$
so that
\begin{equation}\label{e:egalite} F(j)^{c}  \cap D( j)^{c} \subset \left\{ T^{'}_{\I_{mj}}(m(j+1)) =  T^{0}_{\I_{mj}}(m(j+1))    \right\}.\end{equation}


\begin{lemma}\label{l:controle1}
There exist $\lambda_{1}(a)$ and $\lambda_{2}>0$, not depending on $m,\ell,\alpha$,  such that 
$$\P(D(j) )    \leq    4a \alpha (m\ell)^{2}  G_{\alpha
  (m\ell)^{2}}(m\ell))  +
 \lambda_{1}(a) \exp \left(  - \lambda_{2} \alpha (m\ell)^{2}    \right)  =: \lambda.$$
\end{lemma}

\begin{proof}[Proof of Lemma~\ref{l:controle1}]

Consider the random walks born at sites $(x,i)$ for $x \leq mj- \alpha (m\ell)^{2}$. 
By Lemma~\ref{l:neglig-temps-atteinte} choosing $\gamma=1$ and $\theta>0$ small enough so that $g_{\gamma}(\theta)>0$, we obtain the existence of
$\lambda_{1}(a)>0$ and $\lambda_{2}>0$ such that the probability that any of the walks born at a site $(x,i)$ with $x \leq mj- \alpha (m\ell)^{2}$
hits $mj$ before time  $\alpha (m\ell)^{2}$ is $\leq  \lambda_{1}(a) \exp \left(  - \lambda_{2} \alpha (m\ell)^{2}    \right)$.
 On the other hand, for $mj - \alpha (m\ell)^{2} <   x \leq mj-m\ell$, the probability that a walk started at $x$
hits $mj$ before time $\alpha (m \ell)^{2}$ is less than the corresponding probability for the walk started at $mj-m\ell$,  that is, $1-\bar{G}_{\alpha (m \ell)^{2}}(m\ell)$.
In turn, this probability is less than $2 G_{\alpha (m\ell)^{2}}(m\ell)$.
A union bound over all the corresponding events yields the result.

\end{proof}

\begin{lemma}\label{l:controle2}
There exist $V_{1}(m),V_{2}(m)>0$, not depending on $\ell,\alpha$, such that
for all $j$, 
$$\P( F(j)) \leq  V_{1}(m) \exp\left( -V_{2}(m)  \alpha^{1/2} \ell\right).$$
\end{lemma}


\begin{proof}

By translation invariance, we can assume that $j=0$. Let $t = \alpha (m\ell)^{2}$. Since $F(0)$ implies that no random walk 
born at a site $  - m \ell +1 \leq x \leq 0$ hits $1$ before time $\alpha  (m\ell)^{2}$, one has that 
$\P(F(0))= \prod_{  - m \ell +1  \leq  x \leq 0 }  \bar{G}_{t}(1-x)^{a}$. 
Since $0 \leq \alpha \leq 1$, we see that $t^{1/2} \leq m \ell$, so that 
$\P(F(0)) \leq  \prod_{ - \floor{t^{1/2}} +1  \leq  x \leq 0 }  \bar{G}_{t}(1-x)^{a}.$
Using monotonicity of $\bar{G}_{t}$, we deduce that 
$\P(F(0)) \leq \bar{G}_{t}( \floor{t^{1/2}})^{a \floor{t^{1/2}}}$.
 
By the central limit theorem,  $\lim_{t \to +\infty} G_{t}  \left(\floor{t^{1/2}}\right) > 0$, 
so that, since $\bar{G}_{t} \leq 1 - G_{t}$,  $\limsup_{t \to +\infty} \bar{G}_{t}  \left(\floor{t^{1/2}}\right) < 1$.
As a consequence, we can find  $\rho > 0$, and $t_{0} \geq 0$ such that, for all $t \geq t_{0}$, 
$\bar{G}_{t} \left(\floor{t^{1/2}}\right) \leq 1-\rho$.
For $t \geq t_{0}$, we deduce that 
$P(F(0)) \leq (1-\rho)^{a \floor{t^{1/2}}}$. 
For $t \leq t_{0}$, we see that we can find a large enough $V_{1}$ such that 
$P(F(0)) \leq V_{1}(1-\rho)^{a \floor{t^{1/2}}}$, using only the trivial bound $\P(F(0)) \leq 1$. 

\end{proof}

\begin{lemma}\label{l:assoc}
For all $t \geq 0$, the  events  $\left\{ \sum_{j=0}^{\floor{n/m}} T^{0}_{\I_{mj}}(m(j+1)) \geq cn  \right\}$ and
$\bigcup_{j=0}^{ \floor{n/m}} D(j)$
are negatively associated.
\end{lemma}

\begin{proof}

%

For an integer $K \geq 1$, let
$$\A^{0}(u,i,v,K):=\inf \left\{ \sum_{k=1}^{m} \tau_{k}(u,i);
\,  u + \sum_{k=1}^{m} \varepsilon_{k}(u,i,\epsilon)  = v ,  \,
1 \leq m \leq K \right\}.$$
Similarly, let 
$$\B_{j}(u,i,v,K) := \inf \left\{ \sum_{k=1}^{m} \tau^{j}_{k}(u,i);  \,  u + \sum_{k=1}^{m} \varepsilon^{j}_{k}(u,i)  = v ,  \, 1 \leq m \leq K\right\},$$
and let
$$ \C_{j}(u,i,v,K)   :=  \left\{   \begin{array}{l}   \B_{j}(u,i,v,K)  \mbox{ if }  u \leq  mj - m\ell \\
  \A^{0}(u,i,v,K) \mbox{ if } u > mj - m\ell
   \end{array} \right. .$$

 Now let $$T^{\epsilon}_{w}(u,K) := \inf \sum_{j=1}^{L-1} \A^{\epsilon}(x_{j},i_{j},x_{j+1},K)$$
 where the infimum is taken over all finite sequences with $L \geq 2$, $x_{1},\ldots,x_{L} \in \mathbb{Z}$ and $i_{1},\ldots, i_{L-1}$ such that
 $(x_{1}, i_{1}) \in A, \, x_{1} \geq -K, \, r<x_{2} < \cdots < x_{L-1}<u, \, x_{L}=u$, $i_{2},\ldots, i_{L-1} \in \{1, \ldots, a \}$.

Similarly, let 
\begin{equation*}  J_{j,K} := \inf \sum_{g=1}^{L-1}  \C_{j}(x_{g},i_{g},x_{g+1},K) , \,
  K_{j,K} := \inf \sum_{g=1}^{L-1}  \A^{0}_{j}(x_{g},i_{g},x_{g+1},K) ,\end{equation*}
 where in both cases the infimum is taken over all finite sequences with $L \geq 2$, $x_{1},\ldots,x_{L} \in \mathbb{Z}$ and $i_{1},\ldots, i_{L-1}$ such that
 $-K \leq x_{1} \leq mj- m \ell$,  $ mj<x_{2} < \cdots < x_{L-1}<m(j+1), \, x_{L}=m(j+1)$,
 $i_{1},i_{2},\ldots, i_{L-1} \in \{1, \ldots, a \}$.

 Observe that $\P-$almost surely, for all $u<v$, the sequence $(T^{0}_{\I_{u}}(v,K))_{K \geq 1}$
is ultimately stationary, and
that its limiting value is $T^{0}_{\I_{u}}(v)$. Similarly,  $\P-$almost surely, the sequences $(J_{j,K})_{K \geq 1}$
and $(K_{j,K})_{K \geq 1}$ are ultimately stationary, and their respective limits are
 $J_{j}$
and $K_{j}$.

Then let $S_{q,K}:= \sum_{p = 0}^{q} T^{0}_{\I{p m(\ell+1)}}(p m(\ell+1)+m ,K ) $
and
 $$D(j,K):= \left\{    \min( J_{j,K} , K_{j,K} )  <   \alpha (m\ell)^{2} \right\}.$$
Now let $g_{1} := \un \left( S_{q} \geq t \right) $,
$ g_{2} := \un \left(  \bigcup_{p=0}^{q} D(p(\ell+1)) \right)$,
and
$g_{1,K}:=\un \left( S_{q,K} \geq t \right) $ and
$ g_{2,K} := \un \left(  \bigcup_{p=0}^{q} D(p(\ell+1),K) \right)$.

Note that $(g_{1,K})_{K \geq 1}$ is a bounded sequence of random variables that is $\P-$a.s. ultimately stationary and
converging to $g_{1}$ as $K$ goes to infinity. The same holds for
 $(g_{2,K})_{K \geq 1}$ and $g_{2}$.
Now, for every $K$, $g_{1,K}$ and $g_{2,K}$ are functions of a finite number of the random variables
$(U_{n}(x,i), U^{j}_{n}(x,i), \tau^{j}(x,i), \tau(x,i) ; \, n \geq 1, \, x \in \Z, \, 1 \leq i \leq a )$.
Moreover, it is easy to check from the definitions that, with respect to these random variables, $g_{1,K}$ is non-increasing, while $g_{2,K}$ is non-decreasing.
 Since these random variables are independent,
we deduce that $\E( -g_{1,K} g_{2,K} ) \geq \E( -g_{1,K}) \E( g_{2,K} )$ (see e.g.~\cite{EsaProWal}).
Taking the limit as $K \to +\infty$, and using the dominated convergence theorem, we obtain the result.

\end{proof}

Now consider the following inclusion.

\begin{equation}\label{e:encoreuneinclusion}   \left\{    \sum_{j=0}^{\floor{n/m}} T^{0}_{\I_{mj}}(m(j+1)) \geq cn      \right\}  \subseteq    X \cup Y \cup Z ,  \end{equation}
  where
  \begin{align*}
& X:=  \left\{   \sum_{j=0}^{\floor{n/m}} T^{0}_{\I_{mj}}(m(j+1)) \geq cn   \right \}  \cap \bigcup_{j=0}^{\floor{n/m}} D(j)  , \\
 & Y:=  \left\{  \sum_{j=0}^{\floor{n/m}} T^{0}_{\I_{mj}}(m(j+1)) \geq cn   \right \}  \cap  \bigcap_{j=0}^{\floor{n/m}}  \left( D(j)^{c} \cap
    F(j)^{c} \right)        , \\
 & Z:=  \bigcup_{j=0}^{\floor{n/m}}  F(j).
 \end{align*}
Let
$$ f(n,c) := \P\left(  \sum_{j=0}^{\floor{n/m}} T^{0}_{\I_{mj}}(m(j+1)) \geq  cn \right).$$
Then
$ f(n,c) \leq \P(X) + \P(Y) +\P(Z).$
Now, according to Lemmas~\ref{l:assoc}, \ref{l:controle1}  we see that
$$\P(X) \leq  f(n,c)\times ( \floor{n/m} +1) \lambda.$$
>From (\ref{e:egalite}), we see that
$$\P(Y)  \leq   \P \left( \sum_{j=0}^{\floor{n/m}} T^{'}_{\I_{mj}}(m(j+1)) \geq cn \right).$$
>From Lemma~\ref{l:controle2}, we see that,
$$\P(Z) \leq  ( \floor{n/m} +1)   V_{1}(m) \exp\left( -V_{2}(m)  \alpha^{1/2} \ell\right).$$
This leads to the following bound.

\begin{equation}\label{e:cequifaitmarcher}\delta(n) f(n,c)  \leq
 \P \left( \sum_{j=0}^{\floor{n/m}} T^{'}_{\I_{mj}}(m(j+1)) \geq cn \right) +  
 ( \floor{n/m} +1)   V_{1}(m) \exp\left( -V_{2}(m)  \alpha^{1/2} m \ell\right).
 \end{equation}
where
$ \delta(n) :=  1 -  ( \floor{n/m} +1) \lambda$.



Using the independence properties of the random variables
$T^{'}_{\I_{mj}}(m(j+1))$ (Lemma~\ref{l:tempsindep}), and the union bound, we see that
the following inequality holds
\begin{equation*} \P \left( \sum_{j=0}^{\floor{n/m}} T^{'}_{\I_{mj}}(m(j+1)) \geq cn \right) \leq (\ell+1)      \mathfrak{F}_{m}^{\otimes
   k(n) } \left(  [ cn  (\ell+1)^{-1}     , +\infty)          \right), \end{equation*}
where $\mathfrak{F}_{m}^{\otimes k}$ denotes the distribution of the sum of $k$ independent copies of $T^{0}_{\I_{0}}(m)$, 
and where $k(n):=1+\floor{ \textstyle{\frac{n/m - 1 }{\ell+1}  } }$.

\subsubsection{Step 2: Large deviations estimates for i.i.d. random variables}

We start with a general bound on the tail of $T^{0}_{w}(0,m)$. 

\begin{lemma}\label{l:queue}
 There exist $A_{m},c_{m}>0$ such that, for all
 $w=(F,0,A)$
and $t\ge 0$
$$ \P( T^{0}_{w}(0,m) \geq t ) \leq A_{m} \exp \left( -c_{m} H_{w}(-\floor{t^{1/2}})  \right).$$
\end{lemma}


\begin{proof}

Observe that the event $T^0_{w}(m) \geq t$ implies that none of the random walks born at sites $F(x,i)$, $(x,i) \in A$ has hit $m$ before time $t$.
As a consequence,  $\P(T^{0}_{w}(m) \geq t) \leq   \prod_{x=0}^{-\floor{t^{1/2}}}  \bar G_{t}(-x+m)^{\eta_{w}(x)}.$
Using monotonicity of  $\bar{G}_{t}$, we deduce that 
  $\P(T^{0}_{w}(m) \geq t) \leq  G_{t}(m+\floor{t^{1/2}})^{H_{w}(-\floor{t^{1/2}})}$.
Re-using the notations of the proof of Lemma~\ref{l:controle2}, we see that, for all $t \geq t_{0}$, 
  $\P(T^{0}_{w}(m) \geq t) \leq (1-\rho)^{H_{w}(-\floor{t^{1/2}})}$.
  Now, for $t \leq t_{0}$, we can find $A_{m}$ such that, using only the trivial bounds 
  $H_{w}(-\floor{t^{1/2}}) \geq 0$ and $\P(T^{0}_{w}(m) \geq t) \leq 1$, 
 $\P(T^{0}_{w}(m) \geq t) \leq
 A_{m}  (1-\rho)^{ H_w( -\floor{t^{1/2}})}$
 for all $0 \leq t \leq t_{0}$.
\end{proof}

\begin{remark}
The lower bound~(\ref{e:slowdown-lower-bound}) shows that
the upper bound of Lemma~\ref{l:queue} yields the right order of magnitude for the tail of $T_{w}(m)$,
 at least when  $U(\eta_{w}) < 2$.
\end{remark}

The probabilities of large deviations for $\mathfrak{F}^{\otimes k}_{m}$ are 
 described by
 the following lemma, whose proof is deferred to Appendix~\ref{a:ldev-racine}.

\begin{lemma}\label{l:devracine}

Let $(R_{j})_{j \geq 1}$ be a sequence of i.i.d. non-negative
random variables with common distribution $\mu$. Let $ M := \int x
d\mu(x)$. Assume that there exist $A,c>0$ such that for
every $x\ge 0$
\begin{equation}\label{e:queue} \mu( [x,+\infty) ) \leq A \exp \left(-c x^{1/2} \right).\end{equation} Then $M<+\infty$
and for all $f > M$, there exists $h(f)>0$
and $n(f)$ such that if $n\ge n(f)$
$$ P\left(  n^{-1}(R_{1}+\cdots+R_{n})  \geq f   \right)  \leq
\exp\left( -h(f)n^{1/2} \right).$$
\end{lemma}

\subsubsection{Step 3: Conclusion}

Lemma~\ref{l:devracine} above can be applied to probabilities of large deviations of the form
$\mathfrak{F}_{m}^{\otimes  k } \left(  [ bk     , +\infty)          \right)$, where $b> \E(T_{\I_{0}}(m))$, and our goal is to control probabilities
of the form   $\mathfrak{F}_{m}^{\otimes
   k(n) } \left(  [ cn  (\ell+1)^{-1}     , +\infty)          \right)$.
It is easily checked that
\begin{equation}\label{e:decoup-n-grand} cn   (\ell+1)^{-1}    \geq      k(n)   cm 
\left(    1 + \textstyle{\frac{m (\ell + 1)}{n}  }   \right)^{-1}.\end{equation}


Now observe that Kingman's subadditive ergodic theorem (see e.g.~\cite{Lig2})
can be applied to the sequence of random variables $(T^{0}_{\I_{u}}(v))_{u \leq v}$.
Indeed, these variables are non-negative, integrable (Lemma~\ref{l:queue}), 
and satisfy  the required distributional translation invariance properties. 
We deduce that
\begin{equation}\label{e:lgn-temps} \lim_{m \to +\infty}  m^{-1} \E (T_{\I_{0}}(m)) = v^{-1}.\end{equation}
As a consequence, for all $c> v^{-1}$, we can find $m \geq 1$ large enough so that
\begin{equation*}cm  >  \E (T_{\I_{0}}(m)).\end{equation*}
In the sequel, we assume that $m$ is chosen such that this inequality holds. Now let us choose $\ell:=\ell_{n}=n^{1/3}$.
Taking into account Lemmas~\ref{l:queue}, \ref{l:devracine}, and~(\ref{e:decoup-n-grand}), we see that,
as $n$ goes to infinity,  there exists a constant $h_{1}>0$ such that
 \begin{equation}\label{e:la-partie-iid}\mathfrak{F}_{m}^{\otimes
   k(n) } \left(  [ cn  (\ell_{n}+1)^{-1}     , +\infty)          \right) = O \left( \exp(-h_{1}n^{1/3}) \right)  .\end{equation}

Now, for  $0<\zeta<1/2$, let us choose $\alpha:=\alpha_{n}=n^{-\zeta}$, and consider Inequality~(\ref{e:cequifaitmarcher}).
With our definitions, $\alpha_{n}^{1/2}  (m \ell_{n})  = m n^{1/3- \zeta/2}$ while $m \ell_{n} = m n^{1/3}$.
As a consequence, a moderate deviations bound for the simple random walk (see e.g.~\cite{DemZei})  yields that
$G_{\alpha_{n} (m \ell_{n})^{2} }( m \ell_{n}+1) = O\left( \exp( -h_{2} n^{\zeta}) \right)$ for some constant $h_{2}>0$, whence the fact that
$ \delta(n)   = 1+ o(1)$.
Using~(\ref{e:la-partie-iid}),  we see that Inequality~(\ref{e:cequifaitmarcher}) entails that, for large $n$,
$$   f(n,c) \leq  O\left( \exp \left(  -h_{3} n^{1/3- \zeta/2}     \right)     \right).$$
Since $\zeta$ can be taken arbitrarily small, the conclusion of Proposition~\ref{p:bornesubadd} follows.


\begin{remark}

In view of~(\ref{e:slowdown-lower-bound}) and~(\ref{e:slowdown-upper-bound}), we see that our upper and lower bounds on slowdown probabilities do not match. One may wonder
whether it is possible to improve upon either of these bounds so as to find the exact order of magnitude of slowdown large deviations probabilities.
What we can prove (the details are not given here) is that the $\exp\left( - n^{1/3 + o(1)}  \right)$ bound in Proposition~\ref{p:bornesubadd}
gives the best order of magnitude that can be reached by following our proof strategy based on subadditivity.
Indeed, despite the fact that each $T^{0}_{\I_{mj}}(m(j+1))$ has a tail decaying roughly as $\exp\left( - t^{1/2}\right)$,
so that the probabilities of large deviations above the mean would be of order $\exp\left( -n^{-1/2}\right) $
if these random variables were independent,
the positive dependence between these variables makes such large deviations much more likely, with probabilities of order
$\exp\left( -n^{1/3}\right) $.







\end{remark}

\begin{remark}\label{r:pourquoi-assoc}

One may wonder whether the use of association (see Lemma~\ref{l:assoc}) is really needed in the proof.
Indeed, a simpler approach would be to bound the probability
 of the event $X$ in~(\ref{e:encoreuneinclusion}) above by
$\P \left( \bigcup_{j=0}^{\floor{n/m}} D(j) \right)$.
By properly choosing  $\alpha_{n}$ and $\ell_{n}$, we could make this  probability of the order
of $\exp\left(- n^{1/3+o(1)}\right)$,
  compared to the $\exp\left(-h_{2} n^{-\zeta}\right)$
  obtained in the proof of Proposition~\ref{p:bornesubadd}. However, such a choice
 interferes with the other bounds used in the proof, (making $\alpha_{n}$ smaller increases the
 probability of $F(j)$. The best order of magnitude we could obtain with that simpler method is
  $\exp\left( - n^{2/7 + o(1)}  \right)$.
\end{remark}

\section{Appendix: large deviations of i.i.d. random variables with $\exp(-t^{1/2})$ tails}\label{a:ldev-racine}

Neither the result stated in Lemma~\ref{l:devracine} nor the idea of its proof are new, but we failed in finding a reference providing both a
 statement suited to our purposes and a short proof,
so we chose to give a detailed exposition.

We refer to the paper~\cite{MikNag} for a review of results concerning large deviations of random variables with subexponential tails, and to Theorem 4.1 in~\cite{BalDalKlu} for an example
of a result from which Lemma~\ref{l:devracine} may be derived. See also the recent preprint~\cite{DenDieShn}.

\begin{proof}[Proof of Lemma~\ref{l:devracine}]

Let $A$ and $c$ be as in the statement of the lemma. And let $G$ be defined by $G(x):= \mu( [x,+\infty))$.

Let $A_{n}$ be the following event: $A_{n}:= \bigcap_{1 \leq i \leq n} \{ R_{i} \leq n  \}$.
 By the union bound, $P(A_{n}^{c}) \leq n \mu( [n,+\infty) )$, so that, by Assumption~(\ref{e:queue}) above and Lemma~\ref{l:calcul} below,
  \begin{equation}\label{e:complement}P(A_{n}^{c}) =  O \left[   n  \exp \left( -(c/2)n^{1/2} \right)   \right].\end{equation}

 We now apply the Cram\'er bound for i.i.d. random variables possessing finite exponential moments (see e.g.~\cite{DemZei})
 to the i.i.d. bounded random variables $R_{i,n}$ defined by $R_{i,n}:=\min\left(R_{i},n \right).$ For every $\lambda>0$, the following inequality holds.
  \begin{equation}\label{e:cramer}P \left(    n^{-1}(R_{1,n}+\cdots+R_{n,n})  \geq f         \right)
  \leq \exp \left[- n \lambda f \right] \left[  E\exp\left(\lambda R_{1,n}\right)      \right]^{n}.\end{equation}
Let $\lambda_{n} :=  (c/3) n^{-1/2}$ and $K_{n}:=n^{1/4}$. By definition $E\exp\left(\lambda_{n} R_{1,n}\right) = \int_{[0,n)} \exp( \lambda_{n} x  ) d\mu(x)
   + \exp(\lambda_{n} n ) \mu( [n,+\infty) ).$
Let us split the above integral into $\int_{[0,n)} = \int_{[0,K_{n})} +\int_{[K_{n},n)}$. Fix a real number $\alpha>0$.
Since $\lambda_{n} K_{n}$ goes to zero as $n$ goes to infinity, we have, for all large enough $n$ (depending on $\alpha$), an inequality of the following form:
for every $x \in [0,K_{n})$, $\exp(\lambda_{n} x) \leq  1 + (1+\alpha) \lambda_{n} x$.
Taking the integral in this inequality, we obtain that, for all large enough $n$,
$$ \int_{[0,K_{n})} \exp(\lambda_{n} x) d \mu(x) \leq  \mu  ( [0,K_{n}) )  + (1+\alpha) \lambda_{n}  \int_{[0,K_{n})} x d\mu(x).$$
Since $\alpha$ is arbitrary in the above argument, we see that
\begin{equation}\label{e:decoup1} \int_{[0,K_{n})} \exp(\lambda_{n} x)  d \mu(x) \leq  \mu  ( [0,K_{n}) )  + (1+o(1)) \lambda_{n}  \int_{[0,K_{n})} x d\mu(x).\end{equation}

By definition, $ M =     \int_{[0,K_{n})} x d\mu(x) + \int_{[K_{n},+\infty)} x d\mu(x)$.
Integration by parts yields that
$  \int_{[K_{n},+\infty)} x d\mu(x) =   - \left[   x  G(x) \right]_{K_{n}}^{+\infty}  + \int_{[K_{n},+\infty)} G(x) dx$.
Assumption~(\ref{e:queue}) above says that $G(x) \leq A \exp(-c x^{1/2})$. As a consequence,
$- \left[   x  G(x) \right]_{K_{n}}^{+\infty} \leq A K_{n} \exp \left( -c K_{n}^{1/2} \right) $.
Moreover,  Lemma~\ref{l:calcul} yields that
 $\int_{[K_{n},+\infty)} G(x) dx =  O \left[     \exp \left( -(c/2)K_{n}^{1/2} \right)   \right]$.

Putting the above estimates together, and
using the definitions of $\lambda_{n}$ and $K_{n}$, the above estimates clearly imply that $  \int_{[K_{n},+\infty)} x d\mu(x) = o(\lambda_{n})$.
Similarly,   $\mu  ( [K_{n},+\infty) ) = o(\lambda_{n})$.
As a consequence, Inequality~(\ref{e:decoup1}) above yields that
\begin{equation*} \int_{[0,K_{n})} \exp(\lambda_{n} x)  d\mu(x)\leq  1 + (1+o(1))  M \lambda_{n}.\end{equation*}

We now study $\int_{[K_{n},n)} \exp(\lambda_{n} x)  d\mu(x)$. Integration by parts says that
 $\int_{[K_{n},n)} \exp(\lambda_{n} x)  d\mu(x)    =     - \left[ \exp\left( \lambda_{n} x \right) G(x) \right]_{K_{n}}^{n} + \int_{K_{n}}^{n}  \lambda_{n} \exp\left( \lambda_{n} x \right)  G(x) dx.$
Observe that, with our definitions of $\lambda_{n}$ and $K_{n}$, for every $0 \leq x \leq n$, $\lambda_{n} x \leq (c/3) x^{1/2}$.
As a consequence,
$ \exp\left( \lambda_{n} x \right)  G(x)  \leq    A \exp\left( -(2c/3) x^{1/2} \right)$.
This estimate, together with Lemma~(\ref{l:calcul}), yields that, as $n$ goes to infinity,
$\int_{K_{n}}^{n} \exp\left( \lambda_{n} x \right)  G(x) dx = o(1)$.
Similarly,
 $\left[ \exp\left( \lambda_{n} x \right) G(x) \right]_{K_{n}}^{n} = o(\lambda_{n})$. As a consequence, as $n$ goes to infinity,
 $\int_{[K_{n},n)} \exp(\lambda_{n} x)  d\mu(x)    = o(\lambda_{n})$. Similarly, $\exp(\lambda_{n} n ) \mu( [n,+\infty) ) = o(\lambda_{n})$.

  Finally, we obtain the following estimate:
$E\exp\left(\lambda_{n} R_{1,n}\right) = 1 + \lambda_{n} m (1+o(1))$. As $n$ goes to infinity, an expansion yields that
$\left[E\exp\left(\lambda_{n} R_{1,n}\right) \right]^{n}  = \exp\left(    n M \lambda_{n} (1+o(1))  \right).$
>From Cram\'er's inequality~(\ref{e:cramer}), we obtain that
\begin{equation}\label{e:enfin} P \left(    n^{-1}(R_{1,n}+\cdots+R_{n,n})  \geq f   \right)
 \leq \exp \left(- n \lambda_{n} (f-M) (1+o(1))) \right).
\end{equation}

Now, on the event $A_{n}$, $R_{i}=R_{i,n}$ for all $1 \leq i \leq n$.

As a consequence, $P\left(  n^{-1}(R_{1}+\cdots+R_{n})  \geq f   \right) \leq
 P \left(    n^{-1}(R_{1,n}+\cdots+R_{n,n})  \geq f         \right)  +
P(A_{n}^{c}) $.

The statement of the Lemma now follows from the bound~(\ref{e:complement}) on $P(A_{n}^{c})$ and the large deviations bound~(\ref{e:enfin}) for $R_{1,n}+\cdots+R_{n,n}$.

\end{proof}

\begin{lemma}\label{l:calcul}
For every $\nu>0$,
as $x \to +\infty$,
$$\int_{x}^{+\infty} \exp\left(-\nu u^{1/2}\right)   du = O \left[     \exp \left( -(\nu/2)x^{1/2} \right)   \right].$$

\end{lemma}

\begin{proof}[Proof of Lemma~\ref{l:calcul}]

Observe that there exists $d_{1}>0$ such that, for every $u\geq 1$, $u^{1/2} \exp\left( -(\nu/2)u^{1/2}\right) \leq  d_{1}$.
As a consequence, $ \exp\left( -\nu u^{1/2}\right) \leq d_{1}u^{-1/2}  \exp\left( -(\nu/2)u^{1/2}\right)$, so that
$$  \int_{x}^{+\infty}  \exp\left( - \nu u^{1/2}\right)  du \leq d_{1}  \int_{x}^{+\infty} u^{-1/2}    \exp\left( -(\nu/2) u^{1/2}\right)  du.$$
The r.h.s. of the above inequality is then equal to $ d_{1} (4/ \nu)  \exp\left( -(\nu/2)x^{1/2}\right).$

\end{proof}

\section{Appendix: Polynomial tail of renewal variables when $\epsilon=0$}\label{a:poly-decay}

That $\kappa$ may have a polynomial tail is not necessarily an obstacle for proving an exponential bound for speedup large deviations, as long as $\hat{r}_{\kappa}$
has an exponential tail. However, as we now prove, both $\kappa$ and
$\hat{r}_{\kappa}$ have a polynomial tail under $\Q^{0,\theta}_{w}$ when $w$ satisfies {\bf (G)}.

Let $w$ be such that  $r \times \{1,\ldots, a\} \subset A$, $F(r,i)=r$ for all $1 \leq i \leq a$ and $\phi_{r-L}(w) \leq p$.
Let $A_{t}:=\{  U \geq t \}$, so that   $A_{t}= \{ \tilde{r}_{u}  - \hat{r}_{0} \geq \floor{\alpha_{2} u} \mbox{ for all } u<t \}$.
  Let $B_{t} := \{ \tilde{r}_{t}- \hat{r}_{0}  \leq \floor{2 \alpha(0) t} \}$.
Now choose $K>0$ such that $K \alpha_{2}>2 \alpha(0)$ and consider the event $C_{t}$ that the (at most) $aM$ random walks involved in the definition
 of $\nu_{\tilde{r}_{t}+1}$ remain below their position at time $t$
during the time interval $[t,Kt]$. On $B_{t} \cap C_{t}$, $\tilde{r}_{Kt}- \hat{r}_{0} \leq 2 \alpha(0)t$, so that, with our choice of $K$, for $t$ large enough (non-random), 
$B_{t} \cap C_{t} \subset\{ U  \leq Kt \}$. 

Now, we know that $\Q^{0,\theta}_{w}(A_{t})\geq \Q^{0,\theta}_{w}(U = +\infty) \geq \delta_{2}>0$.
Moreover, it is easily checked that, by the law of large numbers, $\lim_{t \to +\infty} \Q^{0,\theta}_{w}(   B_{t}^{c}  )=0$ uniformly with respect to all $w$ such 
that $F(r,i)=r$ for all $1 \leq i \leq a$. As a consequence, for large enough $t$ (not depending on $w$),   $\Q^{0,\theta}_{w}( A_{t} \cap B_{t}   ) \geq \delta_{2}/2$. 
 Moreover, conditional on $A_{t} \cap B_{t}$, $C_{t}$ has a probability larger that $c t^{-aM/2}$ for some $c>0$. 
  As a consequence, there exists $d>0$ such that, for large enough $t$ (not depending on $w$), $\Q^{0,\theta}_{w} (A_{t} \cap B_{t} \cap C_{t}) \geq d t^{-aM/2}$, so that
  $\Q^{0,\theta}_{w}(  t  \leq  U \leq Kt  ) \geq d t^{-aM/2}$.
 Since $U$, $V$ and $W$ are independent and $ \Q^{0,\theta}_{w}(V=+\infty) \geq \delta_{1}>0$ and $ \Q^{0,\theta}_{w}(W=+\infty)\geq \delta_{3}>0$,
we deduce that  $ \Q^{0,\theta}_{w}(  t  \leq U < +\infty, \, V=+\infty, \, W=+\infty) > d \delta_{1} \delta_{3} t^{-aM/2}$.
Then observe that,  on the event $\{ t  \leq U < +\infty, \, V=+\infty, \, W=+\infty \}$, one has that $D \geq t$ and $\hat{r}_{D} \geq \hat{r}_{t} \geq \floor{\alpha_{2} t}$.
Since $\kappa \geq D \circ \theta_{S_{1}}+S_{1}$ and $\hat{r}_{\kappa} \geq \hat{r}_{D \circ \theta_{S_{1}}+S_{1}}$, this ends the proof.

\section{Appendix: negligibility of remote particles}\label{s:negligeabilite}

\begin{prop}\label{p:negligeabilite}
For any $w \in \L_{\theta}$, $0 \leq \epsilon < 1/2$, and any $t \geq 0$,  with $\P$ probability one,
$$\lim_{K \to -\infty} \sup_{0 \leq s \leq t} \sum_{(x,i) \in A; \, x \leq r+K} \exp(\theta(F^{\epsilon}_{s}(x,i) -r))   = 0.$$
\end{prop}

\begin{proof}

For all $x,i,t$, let $  C_{x,i,t} := \exp(\theta(F^{\epsilon}_{t}(x,i) -r))$
and $\gamma:=[2(\cosh\theta -1)  + 4 \epsilon \sinh \theta ]$.
Let also $$H_{K,k}(s):=\sum_{(x,i) \in A; \,  r+K+k \leq x \leq r+K}   C_{x,i,s} \mbox{ and } H_{K,-\infty}(s):=\sum_{(x,i) \in A; \,  x \leq r+K}  
 C_{x,i,s}.$$

Since for every $(x,i)$,
$  (  C_{x,i,s} \exp ( -  \gamma s ))_{s \geq 0}$ is a c\`adl\`ag martingale,
so is $(H_{K,k}(t) \exp ( -  \gamma t ))_{t \geq 0}    $, and we have the following inequality,
valid for all $\lambda>0$:
\begin{equation*}\P \left(  \sup_{0 \leq s \leq t}    H_{K,k}(s) \exp(-\gamma s)   > \lambda       \right)    \leq 
\lambda^{-1} \E \left(  H_{K,k}(0)  \right).\end{equation*}
Now
$\E \left(   H_{K,k}(0)  \right) =  \sum_{(x,i) \in A; \,  r+K+k \leq x \leq r+K} \exp(\theta(F^{\epsilon}(x,i) -r)).$
We deduce that
\begin{equation}\label{e:borne-deviations-martingale}\P \left(  \sup_{0 \leq s \leq t}  H_{K,k}(s)   > \lambda       \right)    \leq 
\lambda^{-1} \exp ( \gamma t )   \sum_{(x,i) \in A; \,  r+K+k \leq x \leq r+K} \exp(\theta(F^{\epsilon}(x,i) -r)) .\end{equation}
Now observe that, for every $s$, the sequence $(H_{K,k}(s))_{k = 0,-1,\cdots}$ is non-decreasing since we are summing
non-negative terms.
As a consequence,
$\P \left(   \sup_{0 \leq s \leq t}  H_{K,-\infty}(s) > \lambda \right) $ equals $\P\left( \bigcup_{k=0}^{-\infty}   \sup_{0 \leq s \leq t}  H_{K,k}(s) > \lambda  \right),$
which is the probability of the union of a non-decreasing sequence of events, and so is equal to
$\lim_{k \to -\infty} \P \left(  \sup_{0 \leq s \leq t} H_{K,k}(s) > \lambda \right)$.
As a consequence, by~(\ref{e:borne-deviations-martingale}),
\begin{equation}\label{e:borne-proba-deviations}\P \left(   \sup_{0 \leq s \leq t}  H_{K,-\infty}(s) > \lambda \right) 
\leq \lambda^{-1} \exp ( \gamma t )    \sum_{(x,i) \in A; \,  x \leq r+K} \exp(\theta(F^{\epsilon}(x,i) -r)).\end{equation}

Now observe that, for every $s$, the sequence
$   ( \sum_{(x,i) \in A; \,  x \leq r+K} C_{x,i,s})_{K = 0,-1,\cdots}$ is non-increasing, since we are summing non-negative terms.
As a consequence, $ \lim_{K \to -\infty}  \sup_{0 \leq s \leq t} H_{K,-\infty}(s)$ exists, and
$\P \left(  \lim_{K \to -\infty}  \sup_{0 \leq s \leq t }  H_{K,-\infty}(s) > \lambda   \right)$ equals   
$\P \left( \bigcap_{K \leq 0}  \sup_{0 \leq s \leq t}  H_{K,-\infty}(s)  > \lambda   \right),$
which is the probability of the intersection of a non-increasing sequence of events, and so is equal to
$\lim_{K \to -\infty}  \P (   \sup_{0 \leq s \leq t}   H_{K,-\infty}(s)    )  > \lambda ).$
>From Inequality~(\ref{e:borne-proba-deviations}), we see that this last expression equals zero.
\end{proof}

\section{Appendix: estimates on the renewal structure}\label{a:renewal-estimates}

In the sequel
every constant $C_{i}$ or $\delta_{i}$ appearing in the estimates is assumed to depend on
$a,\theta,\epsilon_{0}, \alpha_{1}, \alpha_{2}, p, L, \epsilon$, unless there is a special mention that dependence with respect to some of these parameters is absent.
The notation
$(\xi^{\epsilon}_{s})_{s \geq 0}$ stands for a nearest-neighbor random walk on $\Z$
with jump rate 2 and step
distribution $(1/2+\epsilon) \delta_{+1} + (1/2-\epsilon) \delta_{-1}$, started at zero.
 The probability
measure governing $(\xi^{\epsilon}_{s})_{s \geq 0}$ is denoted by $P$.
We use the shorthand $M':=M/4-1$, which is an integer number according to~(\ref{e:valeur-de-M}).
We also use the notation
$$\L'_{\theta} := \{ w=(F,r,A) \in \L_{\theta}; \, r \times \{1,\ldots, a\} \subset A, \, F(r,i)=r \mbox{ for all }1 \leq i \leq a \}.$$
For every $x \in \Z$, let $M_{t,x,i}:=\sup_{0\le s\le t} Z_{s,x,i}$.
Let also, for $z \in \Z$,
\begin{equation}\label{d:psi}
\psi_{z,t}:=\sum_{(x,i); \, x \leq z, \, (x,i) \in A_{t}} \exp (\theta (M_{t,x,i} - \hat{r}_{t})).
\end{equation}
Let $\mu_{\epsilon}:=  \theta \alpha_{1} - 2 (\cosh \theta  -1) - 4 \epsilon  \sinh \theta  $.
Now, for all $0 \leq \epsilon \leq \epsilon_{0}$,
$\mu_{\epsilon}\geq  \mu_{\epsilon_{0}}$, and,
according to~(\ref{e:alph1-et-2}), $\mu_{\epsilon}\geq  \mu_{\epsilon_{0}} >0$ for all $0 \leq \epsilon \leq \epsilon_{0}$.

\begin{lemma}\label{l:lemme2} There exists $C_{1}<+\infty$ not depending on $\epsilon$ or $L$
such that, for all $0 \leq \epsilon \leq \epsilon_{0}$ and
all $w =(F,r,A)\in \L_{\theta}$,
 $$\Q^{\epsilon,\theta}_{w} ( t < W < +\infty  )  \leq C_{1} \phi_{r-L}(w) \exp(- \mu_{\epsilon} t).$$
\end{lemma}

\begin{proof}
Without loss of generality we assume $r=0$.
Let us first note that

\begin{equation} \nonumber
 \Q^{\epsilon,\theta}_{w}\left[t< W<\infty\right]
\le  \Q^{\epsilon,\theta}_{w}\left[\cup_{s\ge t}\left\{\phi_{-L}(w_s)\ge
e^{\theta(\lfloor\alpha_1 s\rfloor-\hat{r}_s)}\right\}\right].
\end{equation}

By the fact that $s \mapsto M_{s,x,i}$ is nondecreasing, and the union bound,
we deduce that
$$
\Q^{\epsilon,\theta}_{w}\left[t< W<\infty\right]
\leq \sum_{n= \floor{t}}^{+\infty}   \Q^{\epsilon,\theta}_{w}  \left[
\sum_{(x,i) \in A\cap \{ \ldots,-(L-1),-L \} \times \{1, \ldots, a \} }e^{\theta M_{n+1,x,i}}
\ge
 e^{\theta \floor{\alpha_1 n} } \right].$$
Using the Markov inequality, we obtain that
\begin{equation}\label{e:borneW}
\Q^{\epsilon,\theta}_{w}\left[t< W<\infty\right]
\leq \sum_{n= \floor{t}}^{+\infty}  \exp( -   \theta \floor{\alpha_1 n}   )
\sum_{(x,i) \in A\cap \{ \ldots,-(L-1),-L \}  \times \{1,\ldots, a \} }  \E^{\epsilon,\theta}_{w} \left( e^{\theta M_{n+1,x,i}} \right).
\end{equation}

For $(x,i) \in A$, write $(Z_{s,x,i})_{s}$
as the independent sum of a symmetric nearest neighbor random
walk on $\Z$ with rate $2-4\epsilon$ and a Poisson process
with rate $4 \epsilon$.
Using the reflection principle to treat the symmetric part, and the fact that the Poisson process part is non-decreasing,
we deduce that
$$\E^{\epsilon,\theta}_{w} \left( e^{\theta M_{s,x,i}} \right) \leq 2 \exp(\theta F(x,i))
\exp \left( s  \left[ 2 (\cosh \theta  -1) + 4 \epsilon  \sinh \theta  \right]  \right).$$

 Plugging the last identity into~(\ref{e:borneW}) and summing,
we finish the
proof of the Lemma.

\end{proof}

 Define for $t\ge 0$,  and $z\le \hat{r}_{0}$,

\begin{equation} \label{d:Ntz}
N_{t,z}(w_{\cdot}):=e^{\theta \hat{r}_{t}- [2(\cosh\theta -1)  - 4 \epsilon \sinh \theta ] t }\phi_{z}(w_{t}).
\end{equation}

\begin{lemma}
\label{l:lemme3} For all $0 \leq \epsilon \leq \epsilon_{0}$, and all
$w=(F,r,A) \in \L_{\theta}$, the family
  $(N_{t,z})_{t\geq 0}$ is a c\`adl\`ag
  $(\F_t^{\epsilon, \theta})_{t \geq 0}$-martingale with respect to $\Q^{\epsilon,\theta}_{w}$,.
\end{lemma}

\begin{proof} Let us remark that,

\begin{equation} \nonumber
N_{t,z}=\sum_{(x,i) \in A, x\le z} e^{\theta
Z_{t,x,i}-[2(\cosh\theta -1)  - 4 \epsilon \sinh \theta ] t}.
\end{equation}
Now, each one of the terms in the above sum is an $(\F_t^{\epsilon, \theta})_{t \geq 0}$-martingale
Furthermore, since $\phi_z(0)<+\infty$, the martingales
 $\sum_{(x,i) \in A, -n\le x\le z} e^{\theta
Z_{t,x,i}- [2(\cosh\theta -1)  - 4 \epsilon \sinh \theta ] t   }$, converge in $L^1( \Q^{\epsilon,\theta}_{w})$ 
to $N_{t,z}$ as $n\to\infty$. Thus, $(N_{t,z})_{t\ge 0}$ is an
$(\F_t^{\epsilon, \theta})_{t \geq 0}$-martingale.
That the paths are  c\`adl\`ag is an easy consequence of $(w_{s})_{s \geq 0}$ being c\`adl\`ag.

\end{proof}

\begin{lemma}\label{l:lemme4}
For every  $0 \leq  \epsilon \leq \epsilon_{0}$ and $w = (F,r,A)\in \L_{\theta}$,
$$ \Q^{\epsilon,\theta}_{w}\left[W<\infty\right]
\leq \exp(\theta)  \phi_{r-L}(w).$$
\end{lemma}

\begin{proof}
See~\cite{ComQuaRam2}.
\end{proof}

\begin{lemma}\label{l:lemme5}
There exist $0<C_{2},C_{3}<+\infty$ not depending on $\epsilon$ or $L$
 such that, for all $0 \leq \epsilon \leq \epsilon_{0}$,
 $w=(F,r,A) \in \L_{\theta}$ and $t \geq 0$,
$$ \Q^{\epsilon,\theta}_{w}\left[t<V<\infty \right]
\leq L C_{2} \exp(-C_{3} t).$$

\end{lemma}

\begin{proof}

Without loss of generality, assume that $r=0$. Then $\Q^{\epsilon,\theta}_{w}(t < V < +\infty)$ is bounded above
by the probability that one of the random walks born at a site between $-L+1$ and $-1$
is at the right of $\floor{\alpha_{1} s}$ at some
time $s \geq t$.
By coupling, we see that the worst case is when all the walks start at zero, in which case,
by the union bound, the probability is less than $aL$ times the probability
for a single random walk started at zero to exceed $\floor{\alpha_{1} s}$ at some
time $s \geq t$.
Let $\tau:=\inf \{ s \geq t; \, \xi^{\epsilon}_{s} \geq \floor{\alpha_{1} s} \}$.

Using the fact that
$(\exp( \theta \xi^{\epsilon}_{s} - [2(\cosh\theta -1)  - 4 \epsilon \sinh \theta ] s  )  )_{s \geq 0}$
is a martingale, and applying Doob's stopping
theorem,
we obtain the bound
$P( \tau < +\infty) \leq \exp(\theta)  \exp(    - \mu_{\epsilon}t)$.
The result follows.

\end{proof}

\begin{lemma}\label{l:lemme6}
There exists $\delta_{1}>0$ not depending on $\epsilon$ such that, for
all $0 \leq \epsilon \leq \epsilon_{0}$ and
$w=(F,r,A) \in \L_{\theta}$,
$$ \Q^{\epsilon,\theta}_{w}\left[V<\infty \right] \leq 1 - \delta_{1}.$$
\end{lemma}

\begin{proof}

 Without loss of generality
we can assume that $r=0$. Note that the probability $\Q^{\epsilon,\theta}_w[V<\infty]$
is upper bounded by the probability that a random walk within a group
of $aL$ independent ones all initially at site $x=0$, at some
time $t\ge 0$ is at the right of $\floor{\alpha_1 t}$.
 But this probability is $1-f(\epsilon)^{aL}$,
 where
$f(\epsilon):=P(\mbox{for all $s \geq 0$, }\xi^{\epsilon}_{s} \leq \floor{\alpha_{1} s})$.
By coupling, observe that $f$ is a non-increasing function of $\epsilon$.
For $\epsilon=\epsilon_{0}$, the asymptotic speed of the walk is $4 \epsilon_{0}$.
Since, from~(\ref{e:alph1-et-2}) $\alpha_{1}>4 \epsilon_{0}$, an easy consequence of the law of large numbers is that
 $f(\epsilon_{0}) >0$.
This ends the proof.
\end{proof}

\begin{lemma}\label{l:lemme7-unif}
There exists $0<C_{4}<+\infty$ not depending on $\epsilon$ or $L$
such that  for all $\epsilon \leq \epsilon_{0}$  and
$w=(F,r,A) \in \L'_{\theta}$,
 and all $t > 0$,
$$ \Q^{\epsilon,\theta}_{w}\left[t<U<\infty \right] \leq C_{4} t^{-M'}.$$
\end{lemma}

\begin{proof}

The proof given in~\cite{ComQuaRam} for $\epsilon=0$ is based on tail estimates on
the random variables $(\nu_{k})_{k \geq 0}$.
By coupling, for all $0 \leq \epsilon<1/2$, and every $s \geq 0$,
$\Q^{\epsilon,\theta}_{w}( \nu_{k} \geq s  ) \leq \Q^{0,\theta}_{w}( \nu_{k} \geq s  )$.
Thus, the estimate in~\cite{ComQuaRam} is in fact uniform over $\epsilon$.

\end{proof}

\begin{lemma}\label{l:lemme7-addendum-unif}
There exists $0<C_{45}<+\infty$ not depending on $\epsilon$ or $L$
such that  for all $\epsilon \leq \epsilon_{0}$  and all $t > 0$,
$$ \Q^{\epsilon,\theta}_{\I_{0}}\left[    \cup_{s \geq t}     \hat{r}_{s} < \floor{\alpha_{1}s} \right] \leq  C_{45} t^{-M'}.$$
\end{lemma}

\begin{proof}

Since we start with the initial condition $\I_{0}$, we can define a modified auxiliary front 
$(\tilde{r}'_{s})_{s \geq 0}$ by replacing the random variables $(\nu_{k})_{k \geq 0}$ used in the definition of 
$(\tilde{r}_{s})_{s \geq 0}$ by the random variables $(\nu'_{k})_{k \geq 0}$ defined as follows.
Let $\nu'_{0}:=0$ and, for $k \geq 1$, $\nu'_k$ is the first time
one of the  random walks
 $\{(G_{s,z,i})_{s \geq 0}; \, (\hat{r}_{0}+k-M)\le z\le \hat{r}_{0}+k-1, \, 1 \leq i \leq a \}$,
 hits the
site $\hat{r}_{0}+k$. 
With this definition, $ \tilde{r}'_{s} \leq \hat{r}_{s}$ for all $s \geq 0$, and, 
for each $1\le j\le M-1$, the random variables
$\{\nu'_{Mk+j}:k\ge 0\}$ are i.i.d. with finite moment of order $M/2$, whereas this is only true for $\{\nu_{Mk+j}:k\ge 1\}$.
The argument of~\cite{ComQuaRam} used to  prove  Lemma~\ref{l:lemme7-unif} can then be easily adapted to prove the present result.
Alternatively, one can invoke Lemma~\ref{l:deviations-moments}.
 \end{proof}

\begin{lemma}\label{l:lemme7-exp}
For every $0<\epsilon<1/2$,there exist $0<C_{5}(\epsilon), C_{6}(\epsilon)<+\infty$ not depending on $L$ such that, for every
$w=(F,r,A) \in \L'_{\theta}$,
 and every $t > 0$,
$$ \Q^{\epsilon,\theta}_{w}\left[t<U<\infty \right] \leq C_{5}(\epsilon) \exp( - C_{6}(\epsilon) t).$$
\end{lemma}

\begin{proof}

We observe that, for a given $\epsilon>0$, $\nu_{k}$ has an exponentially decaying tail due to the positive bias of the random walks $(G_{s,x,i})_{s \geq 0}$.
Using standard large deviations estimates rather than moment estimates
in the proof of Lemma~\ref{l:lemme7-unif}, we get the result.

\end{proof}

Using a similar argument, we can prove the following Lemma.

\begin{lemma}\label{l:lemme7-addendum-exp}
For all $0<\epsilon \leq \epsilon_{0}$, there exists $0<C_{53}(\epsilon), C_{54}(\epsilon)<+\infty$ 
not depending on $L$ such that  for all $\epsilon \leq \epsilon_{0}$  and all $t > 0$,
$$ \Q^{\epsilon,\theta}_{\I_{0}}\left[    \cup_{s \geq t}     \hat{r}_{s} < \floor{\alpha_{1}s} \right] \leq C_{53}(\epsilon) \exp(-C_{54}(\epsilon) t).$$
\end{lemma}

\begin{lemma}\label{l:lemme7-bis}
There exists $\delta_{2}>0$ not depending on $\epsilon$ such that, for all
$0 \leq \epsilon \leq \epsilon_{0}$,
$w=(F,r,A) \in \L'_{\theta}$,
 and $t > 0$,
$$ \Q^{\epsilon,\theta}_{w}\left[U<\infty \right] \leq 1 - \delta_{2}.$$
\end{lemma}

\begin{proof}

By coupling, we see that $\Q^{\epsilon,\theta}_{w}\left[U<\infty \right]$ is a non-increasing function of $\epsilon$.
Thus the estimate for $\epsilon=0$ proved in~\cite{ComQuaRam} is enough.

\end{proof}

\begin{lemma}\label{l:lemme9-unif}

 Let  $\beta$  be such that $0<\beta<\alpha(0)$.
Then there exists $0<C_{7}(\beta)<\infty$ not depending on $\epsilon$ or $L$ such
 that, for all $0 \leq \epsilon \leq \epsilon_{0}$, the following properties hold for all 
 $w=(F,r,A) \in \L_{\theta}$.

\begin{itemize}

\item[a)] If $r=0$ and $w \in \L'_{\theta}$, and $n \geq 1$,
\begin{equation} \nonumber
\Q^{\epsilon,\theta}_{w}
\left[ \hat{T}(n)>  n/\beta  \right]\le C_{7}(\beta) n^{-a/2}.
\end{equation}

\item[b)] Assume that $r=0$,
$m_{-\floor{L^{1/4}},0}(w)\geq a \floor{L^{1/4}}/2$ and $n \geq 1$. Then,

\begin{equation} \nonumber
\Q^{\epsilon,\theta}_{w}
\left[ \hat{T}(n)>  n/\beta  \right]\le   (C_{7}(\beta) \floor{L^{1/4}} n^{-1/2} )^{a \floor{L^{1/4}}/2} +  C_{7}(\beta) n^{-M'}.
\end{equation}

\item[c)] Assume that $r=0$.  For all $k\ge M$ and $n \geq 1$, we have,

\begin{equation} \nonumber
\Q^{\epsilon,\theta}_{w}
\left[ \hat{T}(n+k)-\hat{T}(k)>  n/\beta  \right]\le C_{7}(\beta) n^{-M'}.
\end{equation}

\end{itemize}
\end{lemma}

\begin{proof}
The proof given in~\cite{ComQuaRam2} for $\epsilon=0$ is based on tail estimates for
the random variables $(\nu_{k})_{k \geq 0}$ and
for hitting times of symmetric random walks,
so that, by coupling,  the estimates
proved in~\cite{ComQuaRam2} are in fact uniform over $\epsilon$.

\end{proof}

\begin{lemma}\label{l:lemme9-unif-complement}

 Let  $\beta$  be such that $0<\beta<\alpha(0)$.
Then there exists $0<C_{37}(\beta)<\infty$ not depending on $\epsilon$ or $L$ such
 that, for all $0 \leq \epsilon \leq \epsilon_{0}$, for all $w = (F,r,A)\in \L_{\theta}$ such that
 $r=0$ and
$m_{-\floor{L^{1/4}},0}(w)\geq a \floor{L^{1/4}}/2$, for all $n \geq 1$,

\begin{equation} \nonumber
\Q^{\epsilon,\theta}_{w}
\left[ \hat{T}(nL)>  nL/\beta  \right]\le   C_{37}(\beta)   (nL^{1/2})^{-M'}.
\end{equation}

\end{lemma}

\begin{proof}
Easy consequence of Lemma~\ref{l:lemme9-unif} b), using the first inequality in~(\ref{e:L}).
\end{proof}

\begin{lemma}\label{l:lemme9-exp}

For all  $0<\epsilon<1/2$ and $\beta$ such that $0<\beta<\alpha(0)$, there exist $0<C_{8}(\epsilon,\beta),C_{9}(\epsilon,\beta)<\infty$ not depending on $L$ such that:
for every $w = (F,r,A)\in \L'_{\theta}$, and $n \geq 1$,
\begin{equation} \nonumber
\Q^{\epsilon,\theta}_{w}
\left[ \hat{T}(n)>  n/\beta  \right]\le C_{8}(\beta,\epsilon) \exp( -C_{9}(\beta,\epsilon) n).
\end{equation}

\end{lemma}

\begin{proof}

Stems from the exponential decay of the tail of  $\nu_{k}$, as in Lemma~\ref{l:lemme7-exp}.

\end{proof}

\begin{coroll}\label{c:corollaire2-unif}
There exists  $0<C_{10},C_{11}<\infty$ not depending on $\epsilon$ or $L$
such that, for all $\epsilon \leq \epsilon_{0}$,
all $w = (F,r,A)\in \L'_{\theta}$ such that  $\phi_{r-L}(w) \leq p$,
and all $t>0$,
\begin{equation} \nonumber
\Q^{\epsilon,\theta}_{w}(t<D<\infty) \leq C_{10}(t^{-M'} + L \exp(-C_{11}t)).
\end{equation}
\end{coroll}

\begin{coroll}\label{c:corollaire2-exp}
For every $0<\epsilon < 1/2$, there exist
$0<C_{12}(\beta,\epsilon), C_{13}(\beta,\epsilon)<\infty$ not depending on $L$ such that,
for all $w = (F,r,A) \in \L'_{\theta}$ such that  $\phi_{r-L}(w) \leq p$, and for all $t>0$,
\begin{equation} \nonumber
\Q^{\epsilon,\theta}_{w}(t<D<\infty) \leq L C_{12}(\beta,\epsilon) \exp( -C_{13}(\beta,\epsilon) t).
\end{equation}
\end{coroll}

\begin{coroll}\label{c:corollaire2-bis}
There exists  $0<\delta_{3}<\infty$  such that, for all $0 \leq \epsilon \leq \epsilon_{0}$,
and all $w = (F,r,A) \in \L'_{\theta}$ such that  $\phi_{r-L}(w) \leq p$,
\begin{equation} \nonumber
\Q^{\epsilon,\theta}_{w}(D<\infty) \leq 1 - \delta_{3}.
\end{equation}
\end{coroll}

\begin{proof}[Proof of the corollaries~\ref{c:corollaire2-unif}, \ref{c:corollaire2-exp} and~\ref{c:corollaire2-bis}]
See~\cite{ComQuaRam2}.
\end{proof}

\begin{lemma}\label{l:lemme11-unif}

There exists $0<C_{14},C_{15}<+\infty$ not depending on $\epsilon$ or $L$
such that, for all $0 \leq \epsilon \leq \epsilon_{0}$, all $w=(F,r,A) \in \L'_{\theta}$ such that $\phi_{r-L}(w) \leq p$, and
all $t >0$, $$\Q^{\epsilon,\theta}_{w}(\hat{r}_{D} - r >  t, D < +\infty) \leq C_{14} \left(t^{-M'} +   L \exp(-C_{15} t) \right).$$

\end{lemma}

\begin{lemma}\label{l:lemme11-exp}
For every $0< \epsilon \leq \epsilon_{0}$, there exists $0<C_{16}(\epsilon), C_{17}(\epsilon) <+\infty$ not depending on $L$
such that, for all $w=(F,r,A) \in \L'_{\theta}$ such that $\phi_{r-L}(w) \leq p$, and for all $t >0$,
$$\Q^{\epsilon,\theta}_{w}(\hat{r}_{D} - r >  t, D < +\infty) \leq L C_{16}(\epsilon) \exp(   - C_{17}(\epsilon) t).$$

\end{lemma}

\begin{proof}[Proof of Lemmas~\ref{l:lemme11-unif} and~\ref{l:lemme11-exp}]

Consider $\gamma_{0}>0$ large enough so that \begin{equation}\label{e:gamma0}c_{\gamma_{0}}(\epsilon_{0}, \theta)>0.\end{equation}
Observe that then $c_{\gamma_{0}}(\epsilon,\theta) \geq c_{\gamma_{0}}(\epsilon_{0},\theta)$ for all $0 \leq \epsilon \leq \epsilon_{0}$.

Observe that by the union bound and the fact that $(\hat{r}_{s})_{s}$ is non-decreasing,
$\Q^{\epsilon,\theta}_{w}(\hat{r}_{D} - r >  t, D < +\infty) \leq
\Q^{\epsilon,\theta}_{w}(\hat{r}_{ t \gamma_{0}^{-1}} - r >  t, D \leq t \gamma_{0}^{-1}) +
\Q^{\epsilon,\theta}_{w}(t \gamma_{0}^{-1} < D < +\infty)$. Moreover, note that, by definition, 
$\phi_{r}(0) \leq \phi_{r-L}(0)+aL$.
Then apply Lemma~\ref{l:lemme10} and Corollaries~\ref{c:corollaire2-unif} and~\ref{c:corollaire2-exp}.

\end{proof}

\begin{lemma}\label{l:lemme12}
Consider $w=(F,r,A) \in \L'_{\theta}$ such that $\phi_{r-L}(w)\leq p$. Then, for all $0 \leq \epsilon \leq \epsilon_{0}$,
$\Q^{\epsilon,\theta}_{w}$-a.s. on the event $\{D<\infty\}$ we have,
\begin{equation} \nonumber
\phi_{r-L}(D)\leq e^{\theta}.
\end{equation}
\end{lemma}

\begin{proof}
See~\cite{ComQuaRam2}.
\end{proof}

\begin{coroll}\label{c:corollaire3}
There exists $0<C_{18}<+\infty$ not depending on $\epsilon$ or $L$,
such that, for all $0 \leq \epsilon \leq \epsilon_{0}$, and all $w=(F,r,A) \in \L'_{\theta}$ satisfying $\phi_{r-L}(w) \leq p$,

\begin{equation} \nonumber
\E^{\epsilon,\theta}_w[\phi_{\hat{r}_D}(D), D<\infty] \leq C_{18} L.
\end{equation}

\end{coroll}

\begin{proof}
See~\cite{ComQuaRam2}.
\end{proof}

\begin{lemma}\label{l:lemme13-unif}
There is a constant $0<C_{19}<+\infty$ not depending on $\epsilon$ or $L$,
such that, for all $0 \leq \epsilon \leq \epsilon_{0}$, and all $w=(F,r,A) \in \L'_{\theta}$:
\begin{itemize}

\item[a)]
$\Q^{\epsilon,\theta}_{w} \left(m_{r,r+\floor{L^{1/4}}}(w_{  \hat{T}(r+ \floor{L^{1/4}} )  })   < a \floor{L^{1/4}}/2  \right) \leq  C_{19} L^{-a/8};$

\item[b)]
$\Q^{\epsilon,\theta}_{w} \left(m_{\hat{r}_{D}+L-\floor{L^{1/4} },\hat{r}_{D}+L}(w_{\hat{r}_{D}+L})
 < a \floor{L^{1/4}} /2  \right) \leq C_{19} L^{-aM'/8(M'+1)}.$

\end{itemize}

\end{lemma}

\begin{proof}

Without loss of generality, assume that $r=0$. For the sake of readability, let $n:=\floor{L^{1/4}}$.
We start with the proof of a).

 Choose $4 \epsilon_{0} <\beta<\alpha(0)$. Then,

\begin{equation}
\label{nn2}
\Q^{\epsilon,\theta}_{w}\left[m_{0,n}(w_{\hat{T}(n)})< \frac{a n}{2}\right]\le
\Q^{\epsilon,\theta}_{w}\left[m_{0,n}(w_{\hat{T}(n)})<\frac{a n}{2},  \hat{T}(n)\le\frac{n}{\beta}\right]
+ \Q^{\epsilon,\theta}_{w} \left[\hat{T}(n)>\frac{n}{\beta}\right].
\end{equation}
Note that the event $\{m_{0,n}(w_{\hat{T})(n)}<a n /2, \hat{T}(n)\le n/\beta\}$ is
contained in the event that at least one particle
born at any of the sites $\floor{n/2},\floor{n/2}+1,\ldots,n$
hits some site $x\le 0$ in
a time shorter than or equal to $n/\beta$. Hence, we can conclude that,

\begin{equation}
\label{nn}
\Q^{\epsilon,\theta}_{w}
\left[m_{0,n}(w_{\hat{T}(n)})< \frac{an}{2}, \hat{T}(n)\le\frac{n}{\beta}\right]
\le a(n+1-\floor{n/2}) P[ \Lambda^{\epsilon}_{n/\beta} \leq -\floor{n/2}],
\end{equation}
where and $\Lambda^{\epsilon}_t:=\inf_{0\le s\le t} \xi^{\epsilon}_{s}$.

Noting that, by coupling, $P[ \Lambda^{\epsilon}_{n/\beta} \leq -n/2]$ is non-increasing
as a function of $\epsilon$, we can assume that $\epsilon=0$.

Now, by the reflection principle,  $P[ \Lambda^{0}_{n/\beta} \leq -n/2] \le 2 P[\xi^{0}_{n/\beta}\le -n/2]$.
Hence, from inequality (\ref{nn}), we see that
$ \Q^{\epsilon,\theta}_{w} \left[m_{0,n}(\hat{T}(n))<a n/2, \hat{T}(n)\le\frac{n}{\beta}\right]$
is bounded above by $a(n+1)P[\xi^{0}_{n/\beta}\le -n/2]$. By a standard large deviations argument, for every $t\ge 0$ and
positive integer $x$, $P[\xi^{0}_t\ge x] \le e^{-t g(x/t)}$, where $g(u)>0$ for all $u>0$.
Hence,
$a(n+1) P[\xi^{0}_{n/\beta}\le -n/2] \le a(n+1)\exp \left\{-\frac{n}{\beta} g(\beta/2)\right\}$.
Finally, using part a) of Lemma~\ref{l:lemme9-unif} to bound the second term of inequality
(\ref{nn2}) and using the fact that $a(n+1)\exp \left\{
-\frac{n}{\beta} g(\beta/2)\right\}\le 1/n^{a/2}$ for $n$ large enough,
we conclude the proof of a).

Now for b),
${\mathbb P}_w[m_{\hat{r}_D+L-n,\hat{r}_D+L}(\hat{T}(\hat{r}_D+L))(w_{\hat{T}(\hat{r}_{D}+L)})<an/2]$ is upper bounded by,
\begin{eqnarray}
\nonumber
&\sum_{k:1\le k\le m} \Q^{\epsilon,\theta}_{w}   [m_{k+L-n,k+L}(w_{\hat{T}(k+L)})<an/2]
+  \Q^{\epsilon,\theta}_{w}   [\hat{r}_D>m,D<\infty]
\end{eqnarray}

Letting $m:=L^{a/(8(M'+1))}$, and using part a) and Lemma~\ref{l:lemme11-unif},
we obtain the result.

\end{proof}

Throughout the sequel, to simplify notation, we will
define on the event $\{D<\infty\}$ for each $n\ge 1$,
$$F_n:=\hat{T}_(\hat{r}_D+nL)-D, \,   F'_{n}:=\hat{T}(\hat{r}_{D}+nL).$$

\begin{lemma}
\label{l:lemme14-unif}
For every $0<\beta<\alpha(0)$, there exist $0<C_{20}(\beta), C_{21}(\beta)<\infty$ not depending on $\epsilon,L$, such that for
all $0 \leq \epsilon \leq \epsilon_{0}$, and all $w=(F,r,A) \in \L_{\theta}$ such that  $m_{r-\floor{L^{1/4}} , r}(w) \geq a \floor{L^{1/4}}/2$,
and $\phi_{r-L}(w)\le p$, and
for all natural $n \ge 1$,

\begin{equation} \nonumber
  \Q^{\epsilon,\theta}_{w}\left[ F_n>
\frac{nL}{\beta}, D<\infty \right]
\le C_{20}(\beta)  (nL^{1/2})^{-M'+1} .
\end{equation}

\end{lemma}

\begin{proof} Without loss of generality
we can assume that initially $r=0$. Note that
$\Q^{\epsilon,\theta}_{w}
\left[ F_n> \frac{nL}{\beta}, D<\infty\right]$
 is upper-bounded by
\begin{equation}
\label{f1}
\sum_{k:1\le k\le  \floor{L^{1/2}}n}
   \Q^{\epsilon,\theta}_{w}\!\!\left[ F_n> \frac{nL}{\beta}
,\hat{r}_D=k, D<\infty \right]
+  \Q^{\epsilon,\theta}_{w}\left[ \hat{r}_D>  n\floor{L^{1/2}}, D<\infty \right].
\end{equation}
Now, on the event $\{D<\infty\}$ we have that
$\hat{T}(\hat{r}_D)\le D$ so that $F_n\le \hat{T}(\hat{r}_D+nL)-\hat{T}(\hat{r}_D)$.
Hence,
\begin{equation} \nonumber
 \Q^{\epsilon,\theta}_{w}\left[ F_n>
\frac{nL}{\beta},\hat{r}_D=k, D<\infty \right]
\le
  \Q^{\epsilon,\theta}_{w}\left[ \hat{T}(k+nL)-\hat{T}(k)>\frac{nL}{\beta}\right].
\end{equation}
Now, by part c) of Lemma~\ref{l:lemme9-unif}, for all $k \geq M$ we have
 $\Q^{\epsilon,\theta}_{w}\left[ \hat{T}(k+nL)-\hat{T}(k)>\frac{nL}{\beta}\right]\le \frac{C_{7}(\beta)}{(nL)^{M'}}$.
 On the other hand for $1\leq k \leq M-1$,
$    \Q^{\epsilon,\theta}_{w}  \left[ \hat{T}(k+nL)-\hat{T}(k)>\frac{nL}{\beta}\right]\le
 \Q^{\epsilon,\theta}_{w} \left[\hat{T}(M+nL)>\frac{nL}{\beta}\right]$.
 
 Now let $\beta<\beta'<\alpha(0)$. Observe that, when $nL^{1/2} \geq M (\beta'/\beta - 1)^{-1}$, 
 $(nL+M)/\beta' \leq nL/\beta$, so that 
$ \Q^{\epsilon,\theta}_{w} \left[\hat{T}(M+nL)>\frac{nL}{\beta}\right] \leq 
\Q^{\epsilon,\theta}_{w} \left[\hat{T}(M+nL)>\frac{nL+M}{\beta'}\right] $   . 
  
Thus,  by Lemma~\ref{l:lemme9-unif-complement}, since  $m_{r-\floor{L^{1/4}} , r}(w) \geq a \floor{L^{1/4}}/2$, we know that
  \begin{equation}\label{e:bornef}\Q^{\epsilon,\theta}_{w} 
   \left[\hat{T}(M+nL)>\frac{nL+M}{\beta'}  \right]\le (C_{37}(\beta)  (nL^{1/2})^{-M'}).\end{equation}
When $nL^{1/2} \leq  M (\beta'/\beta - 1)^{-1}$, the same bound holds, with a possibly larger constant, using only the trivial inequality
$\Q^{\epsilon,\theta}_{w}(\cdot) \leq 1$.
  Using Lemma \ref{l:lemme11-unif} to estimate the second term of display~(\ref{f1}), and 
  combining with~(\ref{e:bornef}), we finish the proof.
\end{proof}

\begin{lemma}
\label{l:lemme14-exp}
For every $0<\epsilon \leq \epsilon_{0}$ and $0<\beta<\alpha(0)$, there exist $0<C_{22}(\beta,\epsilon), C_{23}(\beta,\epsilon)<\infty$
not depending on $L$, such that for
all $w=(F,r,A) \in \L'_{\theta}$ such that $\phi_{r-L}(w)\le p$,
for all natural $n \ge 1$,

\begin{equation} \nonumber
  \Q^{\epsilon,\theta}_{w}\left[ F_n>
\frac{nL}{\beta}, D<\infty \right]
\le C_{22}(\beta,\epsilon) \exp(-C_{23}(\beta,\epsilon) n L ).
\end{equation}

\end{lemma}

\begin{proof}

Consider $\ell>0$ such that $\beta (1+\ell) < \alpha(0)$.

As in the proof of the previous lemma,
$\Q^{\epsilon,\theta}_{w}
\left[ F_n> \frac{nL}{\beta}, D<\infty\right]$
 is upper-bounded by

\begin{equation}
\label{f1-exp}
\sum_{k:1\le k\le  \floor{ \ell nL}}
   \Q^{\epsilon,\theta}_{w}\!\!\left[ F_n> \frac{nL}{\beta}
,\hat{r}_D=k, D<\infty \right]
+  \Q^{\epsilon,\theta}_{w}\left[ \hat{r}_D> \floor{\ell nL}, D<\infty \right].
\end{equation}

By Lemma~\ref{l:lemme11-exp},  $\Q^{\epsilon,\theta}_{w}\left[ \hat{r}_D> \floor{\ell nL}, D<\infty \right]
 \leq  L C_{16}(\epsilon) \exp(   - C_{17}(\epsilon)
\floor{\ell nL}  ) $.
On the other hand, for $1 \leq k \leq \floor{\ell  nL}$,
$ \Q^{\epsilon,\theta}_{w}  \left[ \hat{T}(k+nL)-\hat{T}(k)>\frac{nL}{\beta}\right]\le
 \Q^{\epsilon,\theta}_{w} \left[\hat{T}( \floor{nL (1+\ell) })>\frac{nL}{\beta}\right]$.
By Lemma~\ref{l:lemme9-exp},
$\Q^{\epsilon,\theta}_{w} \left[\hat{T}( \floor{n L (1+\ell) })>\frac{nL}{\beta}\right]
\leq
C_{8}(\beta(1+\ell),\epsilon) \exp( -C_{9}(\beta(1+\ell),\epsilon)   \floor{nL(1+\ell)})$.

\end{proof}

\begin{lemma}
\label{l:lemme16}
Consider $w=(F,r,A) \in \L_{\theta}$ such that $r=0$. Then for all $0 \leq \epsilon \leq \epsilon_{0}$, the following properties hold.

\begin{itemize}

\item[a)] For every $h>0, s>0$ and $n\ge 1$
we have

 \begin{equation}
\Q^{\epsilon,\theta}_{w}\left[\psi_{0,\hat{T}(n)}>
h, \hat{T}(n)< s\right]
\le 2 \frac{\phi_{0}(w)}{ h}     \exp( s  ( 2 (\cosh \theta  -1) + 4 \epsilon  \sinh \theta  ) -\theta n).
\end{equation}

\item[b)]  For every $h>0, s>0$, $k\ge 1$ and $n\ge k$
we have a.s.
\begin{eqnarray*}
\Q^{\epsilon,\theta}_{w}\left[\left.\psi_{k,\hat{T}(n)}-\psi_{k-L,\hat{T}(n)}>
h, \hat{T}(n)-\hat{T}(k)< s\right|{\mathcal F}^{\epsilon, \theta}_{\hat{T}(k)}\right]
\\
 \le 2\frac{aL}{h}
 \exp( s  ( 2 (\cosh \theta  -1) + 4 \epsilon  \sinh \theta  ) -\theta (n-k)).
\end{eqnarray*}
\end{itemize}
\end{lemma}

\begin{proof}
See~\cite{ComQuaRam2}.
\end{proof}

\begin{coroll}\label{c:corollaire4}
There exists $0<C_{24}<+\infty$ not depending on $\epsilon$ or $L$ such that, for all $w=(F,r,A) \in \L_{\theta}$, for all $0 \leq \epsilon \leq \epsilon_{0}$,
$\lambda>0$, $n \geq 1$,
$$\Q^{\epsilon,\theta}_{w}\left[  \psi_{\hat{r}_{D},F'_{n}}  > \lambda, \, F_{n} \leq \alpha_{1}^{-1} nL, \, D<+\infty  \right]
\leq \lambda^{-1} C_{24} L \exp(  -\alpha_{1}^{-1}nL \mu_{\epsilon_{0}}).$$
\end{coroll}

\begin{proof}
See~\cite{ComQuaRam2}.
\end{proof}

\begin{coroll}\label{c:corollaire5}
There exists $0<C_{25}<+\infty$ not depending on $\epsilon$ or $L$ such that, for all $0 \leq \epsilon \leq \epsilon_{0}$,
and all $w=(F,r,A) \in \L'_{\theta}$ such that $m_{r-\floor{L^{1/4}},r}(w) \geq a\floor{L^{1/4}}/2$,
\begin{eqnarray*}\Q^{\epsilon,\theta}_{w}\left[  \{ \psi_{\hat{r}_{D}} > p , \hat{T}(\hat{r}_{D}+L) \}
\cup \{ m_{\hat{r}_{D}+L-\floor{L^{1/4}},\hat{r}_{D}+L}(w_{\hat{r}_{D}+L}) < a\floor{L^{1/4}}/2         \} , \, D<+\infty  \right] \\
\leq    C_{25}  L^{-aM'/(8(M'+1))}.\end{eqnarray*}
\end{coroll}

\begin{proof}
See~\cite{ComQuaRam2}.
\end{proof}

\begin{lemma}\label{l:lemme15}
 Let $q\ge 1$ be an integer.
Consider two sequences $(a_k)_{k\ge 1}$ and $(c_k)_{k\ge 1}$
of non-negative
real numbers such that $\sum_{k=1}^\infty a_k<1$ and such that

\begin{equation}
\label{rec-1}
c_1\le a_1,
\end{equation}
and  for every $m\ge 2$ we have that,
\begin{equation}
\label{recurrence}
c_m\le a_{m}+\sum_{k=1}^{m-1}a_{m-k}c_k.
\end{equation}
 For all integers $q \geq 0$, let $A_{q}:=\sum_{k=1}^{+\infty} a_{k} k^{q}$ and
 $C_{q}:=\sum_{k=1}^{+\infty} c_{k} k^{q}$. For $t \geq 0$,
 let $\mathcal{A}(t):=\sum_{k=1}^{+\infty} a_{k} \exp(tk)$ and  $\mathcal{C}(t):=\sum_{k=1}^{+\infty} c_{k} \exp(tk)$.
The following properties hold:

\begin{itemize}

\item[a)] Assume that $q \geq 1$ is such that $A_{q}<+\infty$. Then  $C_{k}<+\infty$ for all $1 \leq k \leq q$, and
$$  C_{q} \leq (1-A_{0})^{-1} \left( A_{q} + \sum_{k=1}^{q}  \binom{q}{k}  C_{q-k} A_{k}\right).$$

\item[b)] Assume that $\mathcal{A}(t_{0})<+\infty$ for some $t_{0}>0$. Then $\mathcal{A}(t)<1$ for all small enough $t>0$ and, for all such $t$,
$$\mathcal{C}(t)  \leq  (1-\mathcal{A}(t))^{-1} \mathcal{A}(t).$$

\end{itemize}

\end{lemma}

\begin{proof}

Part a) is proved in~\cite{ComQuaRam2}.
As for part b), observe that the power series $a(z):=\sum_{k=1}^{+\infty} a_{k} z^{k}$ has a convergence radius $\geq \exp(t_{0})$.
As a consequence, the map $t \mapsto a(\exp(t))$ is well-defined and continuous for $t \leq t_{0}$. For $t=0$, $a(\exp(t))=\sum_{k=1}^{+\infty}a_{k}<1$ by assumption.
By continuity, $a(\exp(t))<1$ for all $t>0$ small enough.

Summing~(\ref{rec-1}) and~(\ref{recurrence}), we see that, for all $m \geq 1$ and $t \geq 0$,
$\sum_{i=1}^{m} c_{i} \exp(ti) \leq a_{1} \exp(t) + \sum_{i=2}^{m}\left( a_{i} \exp(ti) + \sum_{k=1}^{i-1} a_{i-k} c_{k} \exp(ti) \right)$,
so that
$\sum_{i=1}^{m} c_{i} \exp(ti) \leq  \sum_{i=1}^{m} a_{i} \exp(ti)  + \sum_{k=1}^{m-1}   c_{k} \exp(tk)  \left( \sum_{i=k+1}^{m}  a_{i-k} \exp(t(i-k))\right)$.
As a consequence,  $\sum_{i=1}^{m-1} c_{i} \exp(ti) \leq  \mathcal{A}(t)  +  \mathcal{A}(t) \sum_{k=1}^{m-1} c_{k} \exp(tk))$.
 When $\mathcal{A}(t)<1$, we deduce that $\sum_{i=1}^{m-1} c_{i} \exp(ti) \leq  (1- \mathcal{A}(t))^{-1}\mathcal{A}(t)$.
Letting $m$ go to infinity, we conclude the proof.

 \end{proof}

\begin{lemma}\label{l:technique-suite}

Let $(O,\H,\T)$ be a probability space, and $(\H_{n})_{n \geq 1}$ be a non-decreasing sequence of sub-$\sigma-$algebras of
$\H$. Let $(B_{n})_{n \geq 1}$, $(A^{n}_{k})_{n \geq 2 , \, 0 \leq k \leq n-1}$ and $(B'_{n})_{n \geq 2}$ be sequences of
events in $\H$ such that the following properties hold:

\begin{itemize}

\item[(i)] for all $n \geq 1$, $B_{n} \in \H_{n}$

\item[(ii)] for all $n \geq 2$,
$B_n\subset B_{n-1}\cap \left(B'_n\cup A^n_0\cup A^n_1\cup
\cdots\cup A^n_{n-1}\right)$.

\end{itemize}

Now assume that we have defined a sequence $(a_{n})_{n \geq 1}$ of non-negative real numbers enjoying the following properties:
  \begin{enumerate}
 \item $\T(B_{1}) \leq a_{1}$;
 \item for all $n \geq 2$, $\T(B'_{n} | \H_{n-1}) \leq a_{1}$ a.s.;
 \item for all $n \geq 3$,  $\T(A^{n}_{n-1} | \H_{n-2}) \leq a_{2}$ a.s.;
 \item for all $n \geq 2$, $\T(A^{n}_{0}) \leq a_{n}/2$ a.s.;
 \item for all $n \geq 2$, $\T(A^{n}_{1}) \leq a_{n}/2$ a.s.;
 \item for all $n \geq 4$ and all $2 \leq k \leq  n-2$, $\T(A^{n}_{k} | \H_{k-1}) \leq a_{n-k+1}$ a.s.;
 \end{enumerate}
 then, letting $c_{n}:=\T(B_{n})$ for all $n \geq 1$, the inequalities~(\ref{rec-1}) and~(\ref{recurrence})
are satisfied by the two sequences   $(a_{n})_{n \geq 1}$ and  $(c_{n})_{n \geq 1}$.
\end{lemma}

\begin{proof}

First, observe that Inequality~(\ref{rec-1}) is a mere consequence of assumption (1).
Assume now that $n \geq 2$.
By the union bound,
\begin{eqnarray}\label{inc}
 \T\left(B_n\right)
&\le&
\sum_{k=0 }^{n-1 }\T\left(
A^n_k,B_{n-1}\right)
+\T\left(B'_n, B_{n-1} \right).
\end{eqnarray}

Now, since $B_{n-1} \in \H_{n-1}$, assumption (2) entails that
 $\T\left(B'_n, B_{n-1} \right) \leq a_{1} \T(B_{n-1})$.

 On the other hand, (4) and (5) imply that
 $\T(A^{n}_{0})+\T(A^{n}_{1}) \leq a_{n}$.

 When $n=2$, we deduce from~(\ref{inc}) that
  $\T(B_{n}) \leq  \T(A^{n}_{0})+\T(A^{n}_{1}) +  \T\left(B'_n, B_{n-1} \right)$, so that
  $\T(B_{n}) \leq a_{n} + a_{1} \T(B_{n-1})$, and so~(\ref{recurrence}) is proved for $n=2$.

 Assume now that $n \geq 3$. Since by assumption $B_{n-1} \subset B_{n-2}$,
 $\T(A^{n}_{n-1}, B_{n-1}) \leq \T(A^{n}_{n-1}, B_{n-2})$. Now, thanks to assumption (3) and to the fact that
 $B_{n-2} \in \H_{n-2}$, $\T(A^{n}_{n-1}, B_{n-2}) \leq a_{2} \T(B_{n-2})$.

 For $n=3$,  we deduce from~(\ref{inc}) that
  $\T(B_{n}) \leq  \T(A^{n}_{0})+\T(A^{n}_{1}) + \T(A^{n}_{n-1}, B_{n-1}) + \T\left(B'_n, B_{n-1} \right)$,
  so that
  $    \T(B_{n}) \leq a_{n} + a_{2} \T(B_{n-2}) +    a_{1} \T(B_{n-1})   $,
 and so~(\ref{recurrence}) is proved for $n=3$.

 Assume now that $n \geq 4$.
 For $2 \leq k \leq n-2$, the fact that $B_{n-1} \subset B_{k-1}$ implies that
 $\T(A^{n}_{k}, B_{n-1}) \leq \T(A^{n}_{k}, B_{k-2})$. Since $B_{k-1} \in \H_{k-1}$,
 assumption (6) entails that
 $    \T(A^{n}_{k}, B_{k-1}) \leq a_{n-k+1}    \T(B_{k-1}).$

As a consequence, plugging the previous estimates into Inequality~(\ref{inc}), we obtain that
 $$  \T(B_{n}) \leq  a_{n} + a_{2} \T(B_{n-2}) +    a_{1} \T(B_{n-1}) +  \sum_{k=2 }^{n-2 }
a_{n-k+1}    \T(B_{k-1}),$$
which is exactly~(\ref{recurrence}).

\end{proof}

\begin{lemma}\label{l:lemme17-unif} There exists $0<L_{0}<+\infty$ not depending on $\epsilon$ such that,
for all $L \geq L_{0}$ there exists $0<C_{26}<+\infty$ not depending on $\epsilon$, 
such that for all $0 \leq \epsilon \leq \epsilon_{0}$, the following properties hold.

\begin{itemize}

\item[a)] For all $n \geq 1$, $\Q^{\epsilon,\theta}_{\I_{0}}(J_{0} \geq n ) \leq C_{26}n^{3-M'}$.

\item[b)] For all $w=(F,r,A) \in \L'_{\theta}$ such that
$m_{r-\floor{L^{1/4}},r}(w) \geq a\floor{L^{1/4}}/2$, and $\phi_{r-L}(w) \leq p$,
we have that,
for all $n \geq 1$, $\Q^{\epsilon,\theta}_{w}(J_{\hat{r}_{D}} \geq n , D<+\infty) \leq C_{26}n^{3-M'}$

\item[c)] For all $n \geq 1$, $\Q^{\epsilon,\theta}_{a\delta_{0}}(J_{0} \geq n, U>\hat{T}_{nL} ) \leq C_{26} n^{3-M'}$.

\end{itemize}

\end{lemma}

In the sequel, we use the notation $\F_{t}$ instead of $\F^{\epsilon,\theta}_{t}$ to alleviate notations.

\begin{proof}[Proof of part a)]

For all $n \geq 1$, let
\begin{eqnarray}&
\nonumber
B_n:=\cap_{i=1}^n \left\{\psi_{(i-1)L, \hat{T}(iL)}>p\right\}\cup B'_i,&\\
\nonumber &
B'_i:=\left\{
m_{iL-\floor{L^{1/4}},iL}(w_{\hat{T}(iL)})<a  \floor{L^{1/4}} /2\right\}.
&\end{eqnarray}

Since $\phi_z(w_{t})\le \psi_{z,t}$, the following inequality holds:
\begin{equation}
\Q^{\epsilon,\theta}_{\I_{0}}(J_{0} > n )
\le  \Q^{\epsilon,\theta}_{\I_{0}}( B_n).
\end{equation}

For $n \geq 2$ and $1 \leq k \leq n-1$, let
\begin{equation} \nonumber
\Delta^{n}_k:=\psi_{kL,\hat{T}(nL)}- \psi_{(k-1)L, \hat{T}(nL)},
\end{equation}
and let
\begin{equation} \nonumber
A^n_0:=\left\{\psi_{0,\hat{T}(nL)}>
{p}/{2^{n-1}}\right\},
\quad A^n_k:=\left\{
\Delta^{n}_k>
{p}/{2^{n-k}}\right\}.
\end{equation}

We now prove that the assumptions (i)-(ii) of Lemma~\ref{l:technique-suite}
are satisfied, with $(O,\H)$ being the space
$\D(\L_{\theta})$ equipped with the cylindrical $\sigma$-algebra, and
probability $\Q^{\epsilon,\theta}_{w}$,
and $\H_{n}:=\F_{\hat{T}(nL)}$ for all $n \geq 1$.

Assumption (i) is  immediate.
 Note that, for $n \geq 2$,
$\psi_{(n-1)L,\hat{T}(nL)}=\psi_{0,\hat{T}(nL)}+\sum_{k=1}^{n-1} \Delta^{n}_k$.
Since $\frac{1}{2^{n-1}}+\sum_{k=1}^{n-1}\frac{1}{2^{n-k}}=1$, we have that
\begin{equation}
\left\{\psi_{(n-1)L,\hat{T}(nL)}>p\right\}\subset
\left\{  \psi_{0,\hat{T}(nL)} > {p}/{2^{n-1}}\right\}
\cup\left[\cup_{k=1}^{n-1}\left\{ \Delta^{n}_k>{p}/{2^{n-k}}\right\}\right],
\end{equation}
so that (ii) is established.

We now look for a sequence $(a_{n})_{n \geq 1}$ such that assumptions (1)-(6) of Lemma~\ref{l:technique-suite}
are satisfied.
Assume that $n \geq 2$.
By the strong Markov property and Lemma~\ref{l:lemme9-unif} c), using the fact that, by~(\ref{e:L}), $L \geq M$, we have for any $1\le k\le n-1$, a.s.
$$
 \Q^{\epsilon,\theta}_{\I_{0}}\left( \hat{T}(nL) - \hat{T}(kL) \geq (n-k)L/\alpha_{1}   |     \F_{\hat{T}((k-1)L)} \right) \leq  C_{7}(\alpha_{1}) ((n-k)L)^{-M'}.
$$
By the strong Markov property again,
 and Lemma~\ref{l:lemme16} b), using the fact that $\mu_{\epsilon} \geq \mu_{\epsilon_{0}}$, we have that a.s.
 \begin{eqnarray*}
\Q^{\epsilon,\theta}_{\I_{0}}\left[  \Delta^{n}_{k}  >   p/2^{n-k}, \hat{T}(nL) - \hat{T}(kL) \leq (n-k)L/\alpha_{1}   | \F_{\hat{T}((k-1)L)} \right]
 \\
  \le  2aLp^{-1} 2^{n-k}
 \exp(-\mu_{\epsilon_{0}}   (n-k)L/\alpha_{1}).
 \end{eqnarray*}
We deduce that, for $n \geq 2$,  and  $1\le k\le n-1$, a.s.
\begin{equation}\label{e:ineg-1}\Q^{\epsilon,\theta}_{\I_{0}}(A_{k}^{n} | \F_{\hat{T}((k-1)L)}  ) \leq   C_{7}(\alpha_{1}) ((n-k)L)^{-M'} + 2aLp^{-1} 2^{n-k}
 \exp(-\mu_{\epsilon_{0}}   (n-k)L/\alpha_{1}).\end{equation}
Similarly, using Lemma~\ref{l:lemme9-unif-complement}, which is possible since $m_{-\floor{L^{1/4}},0}(\I_{0})\geq a \floor{L^{1/4}}/2$,
we have that
$$
 \Q^{\epsilon,\theta}_{\I_{0}}\left( \hat{T}(nL) \geq nL/\alpha_{1} \right) \leq   C_{37}(\alpha_{1}) (nL^{1/2})^{-M'}.
$$

On the other hand, by  Lemma~\ref{l:lemme16} a), we have that
  $$
\Q^{\epsilon,\theta}_{\I_{0}}\left[  \psi_{0,\hat{T}_{nL}}  >   p/2^{n}, \hat{T}(nL)  \leq n L/\alpha_{1} \right]
 \le  2p^{-1} 2^{n-1} \phi_{0}(\I_{0})
 \exp(-\mu_{\epsilon_{0}}   nL/\alpha_{1}).
$$
We deduce that
\begin{equation}\label{e:ineg-2}\Q^{\epsilon,\theta}_{\I_{0}}(A_{0}^{n}) \leq    C_{37}(\alpha_{1}) (nL^{1/2})^{-M'} + 
 2p^{-1} 2^{n-1} \phi_{0}(\I_{0})
 \exp(-\mu_{\epsilon_{0}}   nL/\alpha_{1}).\end{equation}

Now, for $n\ge 2$, by part $a)$ of  Lemma~\ref{l:lemme13-unif}, the strong Markov property, the fact that 
$(n-1)L \leq nL-\floor{L^{1/4}}$
and that there are at least $a$ particles at the rightmost visited site
at time $\hat{T}(nL-\floor{L^{1/4}})$, a.s.
\begin{equation}\label{e:ineg-3}\Q^{\epsilon,\theta}_{\I_{0}}(B'_{n} | \F_{\hat{T}((n-1)L)}) \leq C_{19} L^{-a/8}.\end{equation}
Finally, observe that, by the union bound,
$\Q^{\epsilon,\theta}_{\I_{0}}(B_{1})$ is upper bounded by
$\Q^{\epsilon,\theta}_{\I_{0}}(  \psi_{0,\hat{T}(L)} > p , \hat{T}(L) \leq L/\alpha_{1}   )
+\Q^{\epsilon,\theta}_{\I_{0}}( \hat{T}(L) > L/\alpha_{1}   )
+   \Q^{\epsilon,\theta}_{\I_{0}} ( m_{L-\floor{L^{1/4}},L}(w_{\hat{T}(L)})<a  \floor{L^{1/4}} /2).$

Thanks to Lemma~\ref{l:lemme9-unif} a), Lemma~\ref{l:lemme16} a) and Lemma~\ref{l:lemme13-unif} a),
we obtain that
\begin{equation}\label{e:ineg-4}\Q^{\epsilon,\theta}_{\I_{0}}(B_{1}) \leq   2p^{-1} \phi_{0}(\I_{0})
 \exp(-\mu_{\epsilon_{0}}  L/\alpha_{1}) +  C_{7}(\alpha_{1}) L^{-a/2} +  C_{19} L^{-a/8}.\end{equation}

Now we see that, by Inequalities~(\ref{e:ineg-3}) and~(\ref{e:ineg-4}), (1) and (2) of Lemma~\ref{l:technique-suite}
are satisfied if we let
$$a_{1} :=  2p^{-1} \phi_{0}(\I_{0})
 \exp(-\mu_{\epsilon_{0}}  L/\alpha_{1}) +  C_{7}(\alpha_{1}) L^{-a/2} +  C_{19} L^{-a/8}.$$
Now, for $m \geq 2$, let
\begin{eqnarray*}a_{m}:=   2\left[   C_{7}(\alpha_{1}) ((m-1)L)^{-M'} + 2aLp^{-1} 2^{m-1}
 \exp(-\mu_{\epsilon_{0}}   (m-1)L/\alpha_{1}) \right]    \\
   +   2\left[      C_{37}(\alpha_{1}) (mL^{1/2})^{-M'} +  2p^{-1} 2^{m-1} \phi_{0}(\I_{0})
 \exp(-\mu_{\epsilon_{0}}   mL/\alpha_{1})                         \right].\end{eqnarray*}
Inequalities~(\ref{e:ineg-1}) and~(\ref{e:ineg-2}) entail assumptions (3)-(4)-(5)-(6) of Lemma~\ref{l:technique-suite}.
Note that the sequence $(a_{m})_{m \geq 1}$ depends on $\epsilon_{0}$ but not on $\epsilon$.
Moreover, observe that, for large enough $L$ (not depending on $\epsilon$),  $\sum_{m=1}^{+\infty} a_{m} m^{M'-3} < +\infty$.
 On the other hand, as $L$ goes to infinity,
$\sum_{m=1}^{+\infty} a_{m}$ goes to zero, as can be checked by studying each term in the definition of $(a_{m})_{m \geq 1}$.
Part a) of Lemma~\ref{l:lemme17-unif} then follows from applying Lemma~\ref{l:lemme15}.

\end{proof}

\begin{proof}[Proof of part b)]

We use exactly the same strategy as for part a).

For all $n \geq 1$, let
\begin{eqnarray}&
\nonumber
B_n:=  \cap_{i=1}^n \left\{\psi_{\hat{r}_{D}+(i-1)L, \hat{T}(\hat{r}_{D}+iL)}>p, \, D < +\infty \right\}\cup B'_i,&\\
\nonumber &
B'_i:=\left\{
m_{\hat{r}_{D}+iL-\floor{L^{1/4}},\hat{r}_{D}+iL}(w_{\hat{T}(\hat{r}_{D}+iL)})<a  \floor{L^{1/4}} /2, \, D< +\infty\right\}.
&\end{eqnarray}

Since $\phi_z(w_{t})\le \psi_{z,t}$, the following inequality holds:
\begin{equation}
\Q^{\epsilon,\theta}_{w}(J_{\hat{r}_{D}} > n , \, D<+\infty)
\le  \Q^{\epsilon,\theta}_{w}(B_n).
\end{equation}

For $n \geq 2$ and $1 \leq k \leq n-1$, on $\{D < +\infty \}$, let
\begin{equation} \nonumber
\Delta^{n}_k:=\psi_{\hat{r}_{D}+kL,\hat{T}(\hat{r}_{D}+nL)}- \psi_{\hat{r}_{D}+(k-1)L, \hat{T}(\hat{r}_{D}+nL)},
\end{equation}
and let
\begin{equation} \nonumber
A^n_0:=\left\{\psi_{\hat{r}_{D},\hat{T}(\hat{r}_{D}+nL)}>
{p}/{2^{n-1}}, \, D<+\infty \right\},
\quad A^n_k:=\left\{
\Delta^{n}_k>
{p}/{2^{n-k}}, \, D<+\infty \right\},
\end{equation} for $ 1\le k\le n-1$.

We now prove that the assumptions (i)-(ii) of Lemma~\ref{l:technique-suite}
are satisfied, with $(O,\H, \T)$ being the space $\D(\L_{\theta})$ equipped with the
cylindrical $\sigma-$algebra, and
probability $\Q^{\epsilon,\theta}_{w}$,
and $\H_{n}:=\F_{\hat{T}(\hat{r}_{D}+nL)}$ for all $n \geq 1$.

Assumptions (i) is immediate.
 Note that, for $n \geq 2$, on $\{ D<+\infty \}$
$\psi_{\hat{r}_{D}+(n-1)L,\hat{T}(\hat{r}_{D}+nL)}=\psi_{\hat{r}_{D},\hat{T}(\hat{r}_{D}+nL)}+\sum_{k=1}^{n-1} \Delta^{n}_k$.
Since $\frac{1}{2^{n-1}}+\sum_{k=1}^{n-1}\frac{1}{2^{n-k}}=1$, we have that
\begin{eqnarray*}
\left\{\psi_{\hat{r}_{D}+(n-1)L,\hat{T}(\hat{r}_{D}+nL)}>p, \, D<+\infty \right\}\subset & \\
\left\{  \psi_{\hat{r}_{D},\hat{T}(\hat{r}_{D}+nL)} > {p}/{2^{n-1}}, \, D<+\infty \right\}&
\cup\left[\cup_{k=1}^{n-1}\left\{ \Delta^{n}_k>{p}/{2^{n-k}}, \, D<+\infty \right\}\right],
\end{eqnarray*}
so that (ii) is established.

We now look for a sequence $(a_{n})_{n \geq 1}$ such that assumptions (1)-(6) of Lemma~\ref{l:technique-suite}
are satisfied. Assume that $n \geq 2$.
By the strong Markov property and Lemma~\ref{l:lemme9-unif} c), using the fact that,
by~(\ref{e:L}), $L \geq M$, we have for any $1\le k\le n-1$, on the event $\{ D<+\infty \}$, a.s.
\begin{eqnarray*}
 \Q^{\epsilon,\theta}_{w}\left( \hat{T}(\hat{r}_{D}+nL) - \hat{T}(\hat{r}_{D}+kL) \geq (n-k)L/\alpha_{1}   |     \F_{\hat{T}(\hat{r}_{D}+(k-1)L)} \right)
 \\
 \leq  C_{7}(\alpha_{1}) ((n-k)L)^{-M'}.
\end{eqnarray*}
By the strong Markov property again,
 and Lemma~\ref{l:lemme16} b), using the fact that $\mu_{\epsilon} \geq \mu_{\epsilon_{0}}$, we have that, on $\{  D<+\infty \}$, a.s.
 \begin{eqnarray*}
\Q^{\epsilon,\theta}_{w}\left[  \Delta^{n}_{k}  >   p/2^{n-k}, \hat{T}(\hat{r}_{D}+nL) - \hat{T}(\hat{r}_{D}+kL) \leq (n-k)L/\alpha_{1}   | \F_{\hat{T}(\hat{r}_{D}+(k-1)L)} \right]
 \\
 \le  2aLp^{-1} 2^{n-k}
 \exp(-\mu_{\epsilon_{0}}   (n-k)L/\alpha_{1}).
\end{eqnarray*}
  We deduce that, for $n \geq 2$,  and  $1\le k\le n-1$, on $\{ D<+\infty\}$, a.s.
\begin{equation}\label{e:ineg-1-2}\Q^{\epsilon,\theta}_{w}(A_{k}^{n} | \F_{\hat{T}(\hat{r}_{D}+(k-1)L)}  ) \leq   C_{7}(\alpha_{1}) ((n-k)L)^{-M'} + 2aLp^{-1} 2^{n-k}
 \exp(-\mu_{\epsilon_{0}}   (n-k)L/\alpha_{1}).\end{equation}
Similarly, using Lemma~\ref{l:lemme14-unif}, which is possible since $m_{-\floor{L^{1/4}},0}(w)\geq a \floor{L^{1/4}}/2$
and $\phi_{r-L}(w) \leq p$,
we have that
$$
 \Q^{\epsilon,\theta}_{w}\left( \hat{T}(\hat{r}_{D}+nL) - D \geq nL/\alpha_{1}, \, D<+\infty \right) \leq
   C_{20}(\alpha_{1})  (nL^{1/2})^{-M'+1}.
$$
On the other hand, by  Corollary~\ref{c:corollaire4}, we have that
  $$
\Q^{\epsilon,\theta}_{w}\left[   \psi_{\hat{r}_{D}, \hat{T}_{\hat{r}_{D}}+nL}    >   p/2^{n-1}, \hat{T}(nL) - D \leq n L/\alpha_{1} \right]
 \le  p^{-1} 2^{n-1}  C_{24} L \exp(  -\alpha^{-1}nL \mu_{\epsilon_{0}}).
$$
We deduce that
\begin{equation}\label{e:ineg-2-2}\Q^{\epsilon,\theta}_{w}(A_{0}^{n}) \leq
  C_{20}(\alpha_{1})  (nL^{1/2})^{-M'+1}
 + p^{-1} 2^{n-1}  C_{24} L \exp(  -\alpha_{1}^{-1}nL \mu_{\epsilon_{0}}).\end{equation}

Now, for $n\ge 2$, by part $a)$ of  Lemma~\ref{l:lemme13-unif}, the strong Markov property,
 the fact that 
$(n-1)L \leq nL-\floor{L^{1/4}}$, 
and that there are at least $a$ particles at the rightmost visited site
at time $\hat{T}(\hat{r}_{D}+nL-\floor{L^{1/4}})$, on $\{ D<+\infty \}$, a.s.
\begin{equation}\label{e:ineg-3-2}\Q^{\epsilon,\theta}_{w}(B'_{n} | \F_{\hat{T}(\hat{r}_{D}+(n-1)L)}) \leq C_{19} L^{-a/8}.\end{equation}

Finally, observe that, by Corollary~\ref{c:corollaire5},
\begin{equation}\label{e:ineg-4-2}\Q^{\epsilon,\theta}_{w}(B_{1}) \leq
C_{25}  L^{-aM'/(8(M'+1))}.
\end{equation}

Now we see that, by Inequalities~(\ref{e:ineg-3-2}) and~(\ref{e:ineg-4-2}), (1) and (2) of Lemma~\ref{l:technique-suite} 
are satisfied if we let
$$a_{1} :=   C_{19} L^{-a/8}     +   C_{25}  L^{-aM'/(8(M'+1))} .$$

Now, for $m \geq 2$, let
\begin{eqnarray*}a_{m}:=   2\left[   C_{7}(\alpha_{1}) ((m-1)L)^{-M'} + 2aLp^{-1} 2^{m-1}
 \exp(-\mu_{\epsilon_{0}}   (m-1)L/\alpha_{1})  \right]    \\
   +   2\left[     C_{20}(\alpha_{1}) (mL^{1/2})^{-M'+1}
 + p^{-1} 2^{m-1}  C_{24} L \exp(  -\alpha_{1}^{-1}mL \mu_{\epsilon_{0}})                  \right].\end{eqnarray*}

Inequalities~(\ref{e:ineg-1-2}) and~(\ref{e:ineg-2-2}) entail assumptions (3)-(4)-(5)-(6) of Lemma~\ref{l:technique-suite}.

Note that the sequence $(a_{m})_{m \geq 1}$ depends on $\epsilon_{0}$ but not on $\epsilon$.
Moreover, observe that, for large enough $L$ (not depending on $\epsilon$),  $\sum_{m=1}^{+\infty} a_{m} m^{M'-3} < +\infty$.
 On the other hand, as $L$ goes to infinity,
$\sum_{m=1}^{+\infty} a_{m}$ goes to zero, as can be checked by studying each term in the definition of $(a_{m})_{m \geq 1}$.
Part b) of Lemma~\ref{l:lemme17-unif} then follows from applying Lemma~\ref{l:lemme15}.
\end{proof}

\begin{proof}[Proof of part c)]

For all $n \geq 1$, let
\begin{eqnarray}&
\nonumber
B_n:= \cap_{i=1}^n \left\{\psi_{(i-1)L, \hat{T}(iL)}>p, \, U>\hat{T}(iL)\right\}\cup B'_i,&\\
\nonumber &
B'_i:=\left\{
m_{iL-\floor{L^{1/4}},iL}(w_{\hat{T}(iL)})<a  \floor{L^{1/4}} /2, \,  U>\hat{T}(iL)   \right\}.
&\end{eqnarray}

Since $\phi_z(w_{t})\le \psi_{z,t}$, the following inequality holds:
\begin{equation}
\Q^{\epsilon,\theta}_{w}(J_{0} > n, U>\hat{T}_{nL} )
\le  \Q^{\epsilon,\theta}_{\I_{0}}( B_n).
\end{equation}
For $n \geq 2$ and $1 \leq k \leq n-1$, let
\begin{equation} \nonumber
\Delta^{n}_k:=\psi_{kL,\hat{T}(nL)}- \psi_{(k-1)L, \hat{T}(nL)},
\end{equation}
and let
\begin{equation} \nonumber
A^n_0:=\left\{\psi_{0,\hat{T}(nL)}>
{p}/{2^{n-1}}, \, U>\hat{T}(nL)\right\},
\quad A^n_k:=\left\{
\Delta^{n}_k>
{p}/{2^{n-k}}, \, U>\hat{T}(nL) \right\},
\end{equation} for $ 1\le k\le n-1$.

We now prove that the assumptions (i)-(ii) of Lemma~\ref{l:technique-suite}
are satisfied, with $(O,\H)$ being the space $\D(\L_{\theta})$ equipped with the cylindrical $\sigma-$algebra and
probability $\Q^{\epsilon,\theta}_{a \delta_{0}}$,
and $\H_{n}:=\F_{\hat{T}(nL)}$.
Assumption (i) is immediate. Assumption (ii) is proved as in a).

We now look for a sequence $(a_{n})_{n \geq 1}$ such that assumptions (1)-(6) of Lemma~\ref{l:technique-suite}
are satisfied.

Assume that $n \geq 2$. Exactly as in part a), we can prove that,
for $n \geq 2$,  and  $1\le k\le n-1$, a.s.
\begin{equation}\label{e:ineg-1-3}\Q^{\epsilon,\theta}_{a \delta_{0}}(A_{k}^{n} | \F_{\hat{T}((k-1)L)}  )
\leq   C_{7}(\alpha_{1}) ((n-k)L)^{-M'} + 2aLp^{-1} 2^{n-k}
 \exp(-\mu_{\epsilon_{0}}   (n-k)L/\alpha_{1}).\end{equation}

Now, note that, on $A_{0}^{n}$, one has $\hat{T}(nL) \leq (nL+1)/\alpha_{2}$ since $U>\hat{T}(nL)$, whence
$\hat{T}(nL) \leq nL/\alpha_{1}$ when $L \geq \alpha_{1}/(\alpha_{2}- \alpha_{1})$.

On the other hand, by  Lemma~\ref{l:lemme16} a), we have that
  $$
\Q^{\epsilon,\theta}_{a \delta_{0}}\left[  \psi_{0,\hat{T}_{nL}}  >   p/2^{n-1}, \hat{T}(nL)  \leq n L/\alpha_{1} \right]
 \le  2p^{-1} 2^{n-1} \phi_{0}(a \delta_{0})
 \exp(-\mu_{\epsilon_{0}}   nL/\alpha_{1}).
$$
We deduce that
\begin{equation}\label{e:ineg-2-3}\Q^{\epsilon,\theta}_{a \delta_{0}}(A_{0}^{n}) \leq     2p^{-1} 2^{n-1}\phi_{0}(a \delta_{0})
 \exp(-\mu_{\epsilon_{0}}   nL/\alpha_{1}).\end{equation}

Exactly as in a), a.s.
\begin{equation}\label{e:ineg-3-3}\Q^{\epsilon,\theta}_{a \delta_{0}}(B'_{n} | \F_{\hat{T}((n-1)L)}) \leq C_{19} L^{-a/8}.\end{equation}

Finally, observe that, by the union bound,
$\Q^{\epsilon,\theta}_{a \delta_{0}}(B_{1})$ is upper bounded by
$\Q^{\epsilon,\theta}_{a \delta_{0}}(  \psi_{0,\hat{T}(L)} > p , \hat{T}(L) \leq L/\alpha_{1}   ) +\Q^{\epsilon,\theta}_{a \delta_{0}}( \hat{T}(L) > L/\alpha_{1}   )  +   \Q^{\epsilon,\theta}_{\I_{0}} ( m_{L-\floor{L^{1/4}},L}(w_{\hat{T}(L)})<a  \floor{L^{1/4}} /2).$

Thanks to Lemma~\ref{l:lemme16} a) and Lemma~\ref{l:lemme13-unif} a) and Lemma~\ref{l:lemme9-unif},
we obtain that
\begin{equation}\label{e:ineg-4-3}\Q^{\epsilon,\theta}_{a \delta_{0}}(B_{1}) \leq   2p^{-1} \phi_{0}(a \delta_{0})
 \exp(-\mu_{\epsilon_{0}}  L/\alpha_{1}) +  C_{7}(\alpha_{1}) L^{-a/2} +  C_{19} L^{-a/8}.\end{equation}

Now we see that, by Inequalities~(\ref{e:ineg-3-3}) and~(\ref{e:ineg-4-3}), (1) and (2) of Lemma~\ref{l:technique-suite}
are satisfied if we let
$$a_{1} := 2p^{-1} \phi_{0}(a \delta_{0})
 \exp(-\mu_{\epsilon_{0}}  L/\alpha_{1}) +  C_{7}(\alpha_{1}) L^{-a/2} +  C_{19} L^{-a/8}.$$

Now, for $m \geq 2$, let
\begin{eqnarray*}a_{m}:=   2\left[     C_{7}(\alpha_{1}) ((m-1)L)^{-M'} + 2aLp^{-1} 2^{m-1}
 \exp(-\mu_{\epsilon_{0}}   (m-1)L/\alpha_{1})\right]    \\
   +   2\left[       2p^{-1} \phi_{0}(a \delta_{0}) 2^{m-1}
 \exp(-\mu_{\epsilon_{0}}   mL/\alpha_{1})     \right].\end{eqnarray*}
Inequalities~(\ref{e:ineg-1-3}) and~(\ref{e:ineg-2-3}) entail assumptions (3)-(4)-(5)-(6) of Lemma~\ref{l:technique-suite}.

Note that the sequence $(a_{m})_{m \geq 1}$ depends on $\epsilon_{0}$ but not on $\epsilon$.
Moreover, observe that, for large enough $L$ (not depending on $\epsilon$),  $\sum_{m=1}^{+\infty} a_{m} m^{M'-3} < +\infty$.
 On the other hand, as $L$ goes to infinity,
$\sum_{m=1}^{+\infty} a_{m}$ goes to zero, as can be checked by studying each term in the definition of $(a_{m})_{m \geq 1}$.
Part c) of Lemma~\ref{l:lemme17-unif} then follows from applying Lemma~\ref{l:lemme15}.

\end{proof}

\begin{lemma}\label{l:lemme17-exp} For every $\epsilon>0$, there exists $L_{1}(\epsilon)<+\infty$ such that, for all
$L \geq L_{1}(\epsilon)$, there exists
$0<C_{27}(\epsilon), C_{28}(\epsilon)<+\infty$  such that the following properties hold.

\begin{itemize}

\item[a)] For all $n \geq 1$, $\Q^{\epsilon,\theta}_{\I_{0}}(J_{0} \geq n ) \leq C_{27}(\epsilon) \exp(-C_{28}(\epsilon) n)$.

\item[b)] For all $w \in \L'_{\theta}$ such that
$m_{r-\floor{L^{1/4}},r}(w) \geq a\floor{L^{1/4}}/2$, and $\phi_{r-L}(w) \leq p$,
we have that,
for all $n \geq 1$, $\Q^{\epsilon,\theta}_{w}(J_{\hat{r}_{D}} \geq n , D<+\infty) \leq C_{27}(\epsilon) \exp(-C_{28}(\epsilon) n)$.

\item[c)] For all $n \geq 1$, $\Q^{\epsilon,\theta}_{a\delta_{0}}(J_{0} \geq n, U>\hat{T}_{nL} )
\leq C_{27}(\epsilon) \exp(-C_{28}(\epsilon) n)$.

\end{itemize}

\end{lemma}

\begin{proof}[Proof of part a)]

We use exactly the same definitions as in the proof of part a) of Lemma~\ref{l:lemme17-unif}, except that we
 look for a different sequence $(a_{n})_{n \geq 1}$ such that assumptions (1)-(6) of Lemma~\ref{l:technique-suite}
are satisfied.
Assume that $n \geq 2$.
By the strong Markov property and Lemma~\ref{l:lemme9-exp},
we have that, for any $1\le k\le n-1$, a.s.
$$
 \Q^{\epsilon,\theta}_{\I_{0}}\left( \hat{T}(nL) - \hat{T}(kL) \geq (n-k)L/\alpha_{1}   |     \F_{\hat{T}((k-1)L)} \right) \leq
 C_{8}(\alpha_{1},\epsilon) \exp( -C_{9}(\alpha_{1},\epsilon) (n-k)L).
$$
As in the proof of Lemma~\ref{l:lemme17-unif}, a.s.
 \begin{eqnarray*}
\Q^{\epsilon,\theta}_{\I_{0}}\left[  \Delta^{n}_{k}  >   p/2^{n-k}, \hat{T}(nL) - \hat{T}(kL) \leq (n-k)L/\alpha_{1}   | \F_{\hat{T}((k-1)L)} \right]
  \\
   \le  2aLp^{-1} 2^{n-k}
 \exp(-\mu_{\epsilon_{0}}   (n-k)L/\alpha_{1}).
\end{eqnarray*}
We deduce that, for $n \geq 2$,  and  $1\le k\le n-1$, a.s.
\begin{eqnarray}\label{e:ineg-1-4}\Q^{\epsilon,\theta}_{\I_{0}}(A_{k}^{n} | \F_{\hat{T}((k-1)L)}  ) 
 \leq
C_{8}(\alpha_{1},\epsilon) \exp( -C_{9}(\alpha_{1},\epsilon) (n-k)L) \nonumber
\\
 + 2aLp^{-1} 2^{n-k} \exp(-\mu_{\epsilon_{0}}   (n-k)L/\alpha_{1}).\end{eqnarray}

By Lemma~\ref{l:lemme9-exp} again,
$$
 \Q^{\epsilon,\theta}_{\I_{0}}\left( \hat{T}(nL) \geq nL/\alpha_{1} \right) \leq    C_{8}(\alpha_{1},\epsilon) \exp( -C_{9}(\alpha_{1},\epsilon) nL).
$$

On the other hand, as in the proof of Lemma~\ref{l:lemme17-unif},
  $$
\Q^{\epsilon,\theta}_{\I_{0}}\left[  \psi_{0,\hat{T}_{nL}}  >   p/2^{n-1}, \hat{T}(nL)  \leq n L/\alpha_{1} \right]
 \le  2p^{-1} 2^{n-1} \phi_{0}(\I_{0})
 \exp(-\mu_{\epsilon_{0}}   nL/\alpha_{1}).
$$
We deduce that
\begin{equation}\label{e:ineg-2-4}\Q^{\epsilon,\theta}_{\I_{0}}(A_{0}^{n}) \leq   C_{8}(\alpha_{1},\epsilon) 
\exp( -C_{9}(\alpha_{1},\epsilon) nL) +  2p^{-1} 2^{n-1} \phi_{0}(\I_{0})
 \exp(-\mu_{\epsilon_{0}}   nL/\alpha_{1}).\end{equation}

Now, for $n\ge 2$, as in the proof of Lemma~\ref{l:lemme17-unif}, a.s.
\begin{equation}\label{e:ineg-3-4}\Q^{\epsilon,\theta}_{\I_{0}}(B'_{n} | \F_{\hat{T}((n-1)L)}) \leq C_{19} L^{-a/8}.\end{equation}
Similarly,
\begin{equation}\label{e:ineg-4-4}\Q^{\epsilon,\theta}_{\I_{0}}(B_{1}) \leq   2p^{-1} \phi_{0}(\I_{0})
 \exp(-\mu_{\epsilon_{0}}  L/\alpha_{1}) +  C_{7}(\alpha_{1}) L^{-a/2} +  C_{19} L^{-a/8}.\end{equation}

Now we see that, by Inequalities~(\ref{e:ineg-3-4}) and~(\ref{e:ineg-4-4}), (1) and (2) of Lemma~\ref{l:technique-suite}
are satisfied if we let
$$a_{1} :=  2p^{-1} \phi_{0}(\I_{0})
 \exp(-\mu_{\epsilon_{0}}  L/\alpha_{1}) +  C_{7}(\alpha_{1}) L^{-a/2} +  C_{19} L^{-a/8}.$$
Now, for $m \geq 2$, let
\begin{eqnarray*}a_{m}:=   2\left[ C_{8}(\alpha_{1},\epsilon) \exp( -C_{9}(\alpha_{1},\epsilon) (m-1)L)
+ 2aLp^{-1} 2^{m-1} \exp(-\mu_{\epsilon_{0}}   (m-1) L/\alpha_{1})   \right]    \\
   +   2\left[   C_{8}(\alpha_{1},\epsilon) \exp( -C_{9}(\alpha_{1},\epsilon) mL) +  2p^{-1} 2^{m-1} \phi_{0}(\I_{0})
 \exp(-\mu_{\epsilon_{0}}   mL/\alpha_{1})  \right].\end{eqnarray*}

Inequalities~(\ref{e:ineg-1-4}) and~(\ref{e:ineg-2-4}) entail assumptions (3)-(4)-(5)-(6) of Lemma~\ref{l:technique-suite}.
Now observe that, for $L$ large enough, $\sum_{n=1}^{+\infty} a_{n} \exp(tn) < +\infty$ for $t>0$ small enough.
As $L$ goes to infinity,
$\sum_{n=1}^{+\infty} a_{n}$ goes to zero, as can be checked by studying each term in the definition of $(a_{m})_{m \geq 1}$.
Part a) then follows from applying Lemma~\ref{l:lemme15}.

\end{proof}

\begin{proof}[Proof of part b)]

We re-use exactly the same definitions as in the proof of part b) of Lemma~\ref{l:lemme17-unif}, except that we
 look for a different sequence $(a_{n})_{n \geq 1}$ such that assumptions (1)-(6) of Lemma~\ref{l:technique-suite}
are satisfied.
Assume that $n \geq 2$.
By the strong Markov property and Lemma~\ref{l:lemme9-exp}, we have for any $1\le k\le n-1$, on $\{ D<+\infty \}$ a.s.
\begin{eqnarray*}
 \Q^{\epsilon,\theta}_{w}\left( \hat{T}(\hat{r}_{D}+nL) - \hat{T}(\hat{r}_{D}+kL) \geq (n-k)L/\alpha_{1}   |     \F_{\hat{T}(\hat{r}_{D}+(k-1)L)} \right)
\\
 \leq  C_{8}(\alpha_{1},\epsilon) \exp( -C_{9}(\alpha_{1},\epsilon) (n-k)L).
\end{eqnarray*}
As in Lemma~\ref{l:lemme17-unif}, we have that, on $\{  D<+\infty \}$ a.s.
 \begin{eqnarray*}
\Q^{\epsilon,\theta}_{w}\left[  \Delta^{n}_{k}  >   p/2^{n-k}, \hat{T}(\hat{r}_{D}+nL) - \hat{T}(\hat{r}_{D}+kL) \leq (n-k)L/\alpha_{1}   | \F_{\hat{T}(\hat{r}_{D}+(k-1)L)} \right]
 \\
  \le  2aLp^{-1} 2^{n-k}
 \exp(-\mu_{\epsilon_{0}}   (n-k)L/\alpha_{1}).
\end{eqnarray*}

  We deduce that, for $n \geq 2$,  and  $1\le k\le n-1$, on $\{ D<+\infty\}$, a.s.
\begin{eqnarray}\label{e:ineg-1-5} \Q^{\epsilon,\theta}_{w}(A_{k}^{n} | \F_{\hat{T}(\hat{r}_{D}+(k-1)L)}  ) 
\leq C_{8}(\alpha_{1},\epsilon) \exp( -C_{9}(\alpha_{1},\epsilon) (n-k)L) \nonumber \\
 + 2aLp^{-1} 2^{n-k}
 \exp(-\mu_{\epsilon_{0}}   (n-k)L/\alpha_{1}).\end{eqnarray}

Similarly, using Lemma~\ref{l:lemme14-exp}, which is possible since $\phi_{r-L}(w) \leq p$,
we have that
$$
 \Q^{\epsilon,\theta}_{w}\left( \hat{T}(\hat{r}_{D}+nL) \geq nL/\alpha_{1}, \, D<+\infty \right) \leq
    C_{22}(\alpha_{1},\epsilon)  L \exp(-C_{23}(\alpha_{1},\epsilon) n L ).
$$

As in the proof of Lemma~\ref{l:lemme17-unif}, we have that
  $$
\Q^{\epsilon,\theta}_{w}\left[  \psi_{\hat{r}_{D}, \hat{T}_{\hat{r}_{D}+nL}}  >   p/2^{n-1}, \hat{T}(nL)  \leq n L/\alpha_{1} \right]
 \le  p^{-1} 2^{n-1}  C_{24} L \exp(  -\alpha^{-1}nL \mu_{\epsilon_{0}}).
$$
We deduce that
\begin{equation}\label{e:ineg-2-5}\Q^{\epsilon,\theta}_{w}(A_{0}^{n}) \leq
   C_{22}(\alpha_{1},\epsilon)  L \exp(-C_{23}(\alpha_{1},\epsilon) n L )
 + p^{-1} 2^{n-1}  C_{24} L \exp(  -\alpha_{1}^{-1}nL \mu_{\epsilon_{0}}).\end{equation}

Now, for $n\ge 2$, as in Lemma~\ref{l:lemme17-unif} a.s.
\begin{equation}\label{e:ineg-3-5}\Q^{\epsilon,\theta}_{w}(B'_{n} | \F_{\hat{T}(\hat{r}_{D}+(n-1)L)}) \leq C_{19} L^{-a/8},\end{equation}
and
\begin{equation}\label{e:ineg-4-5}\Q^{\epsilon,\theta}_{w}(B_{1}) \leq
C_{25}  L^{-aM/(16(M+1))}.
\end{equation}

Now we see that, by Inequalities~(\ref{e:ineg-3-5}) and~(\ref{e:ineg-4-5}), (1) and (2) of Lemma~\ref{l:technique-suite}
are satisfied if we let
$$a_{1} :=   C_{19} L^{-a/8}     +   C_{25}  L^{-aM/(16(M+1))} .$$
Now, for $m \geq 2$, let
\begin{eqnarray*}a_{m}:=   2\left[ C_{8}(\alpha_{1},\epsilon) \exp( -C_{9}(\alpha_{1},\epsilon) (m-1)L)  + 2aLp^{-1} 2^{m-1}
 \exp(-\mu_{\epsilon_{0}}   (m-1)L/\alpha_{1})     \right]    \\
   +   2\left[   C_{22}(\alpha_{1},\epsilon)  L \exp(-C_{23}(\alpha_{1},\epsilon) n L )
 + p^{-1} 2^{m-1}  C_{24} L \exp(  -\alpha_{1}^{-1}mL \mu_{\epsilon_{0}})         \right].\end{eqnarray*}
Inequalities~(\ref{e:ineg-1-5}) and~(\ref{e:ineg-2-5}) entail assumptions (3)-(4)-(5)-(6) of Lemma~\ref{l:technique-suite}.
Now observe that, for $L$ large enough, $\sum_{n=1}^{+\infty} a_{n} \exp(tn) < +\infty$ for $t>0$ small enough.
As $L$ goes to infinity,
$\sum_{n=1}^{+\infty} a_{n}$ goes to zero, as can be checked by studying each term in the definition of $(a_{m})_{m \geq 1}$.
Part b) then follows from applying Lemma~\ref{l:lemme15}.
\end{proof}

\begin{proof}[Proof of part c)]
We use exactly the same definitions as in the proof~\ref{l:lemme17-unif} c), except that we look for a different 
sequence $(a_{n})_{n \geq 1}$ such that assumptions (1)-(6) of Lemma~\ref{l:technique-suite}
are satisfied.

Assume that $n \geq 2$. Exactly as in the proof of part a) of the present lemma, we can prove that,
for $n \geq 2$,  and  $1\le k\le n-1$, a.s.
\begin{eqnarray}\label{e:ineg-1-6}
\Q^{\epsilon,\theta}_{\I_{0}}(A_{k}^{n} | \F_{\hat{T}((k-1)L)}  ) \leq
C_{8}(\alpha_{1},\epsilon) \exp( -C_{9}(\alpha_{1},\epsilon) (n-k)L) \nonumber
\\
+ 2aLp^{-1} 2^{n-k} \exp(-\mu_{\epsilon_{0}}   (n-k)L/\alpha_{1}). 
\end{eqnarray}

As in the proof of Lemma~\ref{l:lemme17-unif} c),
 \begin{equation}\label{e:ineg-2-6}\Q^{\epsilon,\theta}_{a \delta_{0}}(A_{0}^{n}) 
  \leq     2p^{-1} 2^{n-1} \phi_{0}(a \delta_{0})
 \exp(-\mu_{\epsilon_{0}}   nL/\alpha_{1}).\end{equation}

Similarly, a.s.
\begin{equation}\label{e:ineg-3-6}\Q^{\epsilon,\theta}_{a \delta_{0}}(B'_{n} | \F_{\hat{T}((n-1)L)}) \leq C_{19} L^{-a/8},\end{equation}
and
\begin{equation}\label{e:ineg-4-6}\Q^{\epsilon,\theta}_{a \delta_{0}}(B_{1}) \leq   2p^{-1} \phi_{0}(a \delta_{0})
 \exp(-\mu_{\epsilon_{0}}  L/\alpha_{1}) +  C_{7}(\alpha_{1}) L^{-a/2} +  C_{19} L^{-a/8}.\end{equation}
Now we see that, by Inequalities~(\ref{e:ineg-3-6}) and~(\ref{e:ineg-4-6}), (1) and (2) of Lemma~\ref{l:technique-suite}
are satisfied if we let
$$a_{1} := 2p^{-1} \phi_{0}(a \delta_{0})
 \exp(-\mu_{\epsilon_{0}}  L/\alpha_{1}) +  C_{7}(\alpha_{1}) L^{-a/2} +  C_{19} L^{-a/8}.$$
Now, for $m \geq 2$, let
\begin{eqnarray*}a_{m}:=   2\left[C_{8}(\alpha_{1},\epsilon) \exp( -C_{9}(\alpha_{1},\epsilon) (m-1)L)
+ 2aLp^{-1} 2^{m-1} \exp(-\mu_{\epsilon_{0}}   (m-1)L/\alpha_{1})     \right]    \\
   +   2\left[       2p^{-1} 2^{m-1} \phi_{0}(a \delta_{0})
 \exp(-\mu_{\epsilon_{0}}   mL/\alpha_{1})     \right].\end{eqnarray*}
Inequalities~(\ref{e:ineg-1-6}) and~(\ref{e:ineg-2-6}) entail assumptions (3)-(4)-(5)-(6) of
Lemma~\ref{l:technique-suite}.
Now observe that, for $L$ large enough, $\sum_{n=1}^{+\infty} a_{n} \exp(tn) < +\infty$ for $t>0$ small enough.
As $L$ goes to infinity,
$\sum_{n=1}^{+\infty} a_{n}$ goes to zero, as can be checked by studying each term in the definition of $(a_{m})_{m \geq 1}$.
Part c) then follows from applying Lemma~\ref{l:lemme15}.
\end{proof}

\begin{lemma}\label{l:deviations-moments}
Let $(Y_{i})_{i \geq 1}$ be a sequence of random variables on a probability space $(O,\H,\T)$, and $(\H_{i})_{i \geq 0}$ an non-decreasing sequence of sub-$\sigma-$algebras of $\H$
 such that $\H_{0}= \{ \emptyset, O \}$. Assume that the following properties hold:

 \begin{itemize}
 \item for all $i \geq 1$, $Y_{i}$ is measurable with respect to $\H_{i}$;
 \item there exists an integer $q \geq 1$ and a constant $0<c_{1}(q)<+\infty$ such that a.s. $\E_{\T}(  Y_{i}^{2q} | \H_{i-1}) \leq c_{1}(q)$.
 \end{itemize}
 Then there exists a constant  $0<c_{2}(q)<+\infty$, depending only on $q$ and $c_{1}(q)$, such that for all $t \geq 0$ and $n \geq 1$,
 $$\T \left( \sup_{k \geq n}   k^{-1} \left|Y_{1}+\cdots Y_{k} -  \sum_{i=1}^{k}\E_{\T}(Y_{i} | \H_{i-1})  \right| \geq t   \right)
 \leq  c_{2}(q) n^{-q}t^{-2q}.$$

\end{lemma}

\begin{proof}

Observe that $\E_{\T}(Y_{i} | \H_{i-1})$ exists and is finite for all $i$ since
 $\E_{\T}(  Y_{i}^{2q} | \H_{i-1}) < + \infty$.
 Now let $Z_{i} := Y_{i} - \E_{\T}(Y_{i} | \H_{i-1})$.
 Observe that, with our assumptions,  $\E_{\T}(  Z_{i} | \H_{i-1}) = 0$ a.s.
Moreover, thanks e.g. to Jensen's inequality, 
 $\E_{\T}(  Z_{i}^{2q} | \H_{i-1}) \leq c_{3}(q)$, where $c_{3}(q)$ depends only on $q$ and $c_{1}(q)$.
 
 We now prove by induction on $\ell$ that, for all $0 \leq \ell \leq q$ there exists a constant $0<c_{4}(\ell)<+\infty$, depending only on $\ell$, $q$ 
 and $c_{1}(q)$, such that, for all $n \geq 1$,
 \begin{equation}\label{e:borne-moments-pairs}\E_{\T}((Z_{1}+\cdots+Z_{n})^{2\ell}) \leq c_{4}(\ell) n^{\ell}.\end{equation}
For $\ell=0$, the result is trivially true for all $n \geq 1$.
Now consider $0 \leq \ell \leq q-1$, assume that the result holds for $\ell$, and let us prove that it holds for $\ell+1$. 
For all $n \geq 1$, $$\E_{\T}((Z_{1}+\cdots+Z_{n+1})^{2\ell+2}) =
\sum_{k=0}^{2\ell+2} \binom{2\ell+2}{k} \E_{\T}((Z_{1}+\cdots+Z_{n})^{2\ell+2-k} Z_{n+1}^{k}).$$
With our assumptions, $\E_{\T}((Z_{1}+\cdots+Z_{n})^{2\ell+1} Z_{n+1})=0$.
Now, by Jensen's inequality,
$\E_{\T}(Z_{n+1}^{2} | \H_{n}) \leq c_{3}(q+1)^{1/(q+1)}$ a.s. 
By our induction hypothesis, we see that 
  $\E_{\T}((Z_{1}+\cdots+Z_{n})^{2\ell}) \leq c_{4}(\ell) n^{\ell}$, with $c_{4}(\ell)$ depending only on $q,\ell$, and $c_{1}(q)$. 
  As a consequence,
 $\E_{\T}((Z_{1}+\cdots+Z_{n})^{2\ell} Z_{n+1}^{2}) \leq c_{4}(\ell)  c_{3}(q)^{1/(q+1)} n^{q}$.
On the other hand, by Jensen's inequality, for $k \geq 3$,   $  \E_{\T} \left| (Z_{1}+\cdots+Z_{n})^{2\ell+2-k}  \right|  \leq
 \E_{\T}((Z_{1}+\cdots+Z_{n})^{2\ell})^{(2\ell+2-k)/2\ell} \leq (c_{4}(\ell)n^{\ell})^{(2\ell+2-k)/2\ell}$.
Similarly, $\E_{\T} \left( \left|  Z_{n+1}^{k} \right|   | \H_{n} \right) \leq c_{3}(q)^{k/2q}$ a.s., so that
$ \left| \E_{\T}((Z_{1}+\cdots+Z_{n})^{2\ell+2-k} Z_{n+1}^{k})  \right| \leq  c_{3}(q)^{k/2q)} (c_{4}(\ell)n^{\ell})^{(2\ell+2-k)/2\ell}$.
Putting these estimates together, we obtain that 
 \begin{eqnarray*}\E_{\T}((Z_{1}+\cdots+Z_{n+1})^{2\ell+2}) -  \E_{\T}((Z_{1}+\cdots+Z_{n})^{2\ell+2})
  \leq \\
  \binom{2\ell+2}{2} c_{4}(\ell)  c_{3}(q)^{1/q} n^{\ell} 
 \\
  +  \sum_{k=3}^{2\ell} \binom{2\ell+2}{2\ell+2-k}  c_{3}(q)^{k/2q} (c_{4}(\ell)n^{\ell})^{(2\ell+2-k)/2\ell}. 
\end{eqnarray*}
Since the are only terms of order $n^{\ell}$ or less in the r.h.s. of the above inequality, summing, we easily deduce that   
$\E_{\T}((Z_{1}+\cdots+Z_{n})^{2\ell+2})  \leq c_{4}(\ell+1) n^{\ell+1}$ for all $n \geq 1$, with a constant $c_{4}(\ell+1)$ depending only on $\ell$, $q$, and 
$c_{1}(q)$, so the induction step from $q$ to $q+1$ is complete.

Now observe that the sequence $(M_{k})_{k \geq 0}$ defined by $M_{0}:=0$ and
$M_{k}:=k^{-1}(Z_{1}+\cdots+Z_{k})$ is a martingale with respect to $(\H_{k})_{k \geq 0}$.
As a consequence, using the maximal inequality for martingales and Inequality~(\ref{e:borne-moments-pairs}), we see that,
for all integers $n \geq 1$ and $\ell \geq 0$,
$$\T \left(   \sup_{2^{\ell} n \leq k \leq 2^{\ell+1}n}  |  M_{k} | \geq t    \right) \leq  c_{4}(q) \left(2^{\ell+1} n \right)^{-q}t^{-2q}.$$
By the union bound,
$$\T \left( \sup_{k \geq n}   k^{-1} \left|Y_{1}+\cdots Y_{k} -  \sum_{i=1}^{k}\E_{\T}(Y_{i} | \H_{i-1})  \right| \geq t   \right) $$
is bounded above by
$$\sum_{\ell=0}^{+\infty} \T \left(   \sup_{2^{\ell} n \leq k \leq 2^{\ell+1}n}  |  M_{k} | \geq t    \right) $$
and so by
$$\sum_{\ell=0}^{+\infty}  c_{4}(q) \left(2^{\ell+1} n \right)^{-q}t^{-2q}.$$
The conclusion follows.
\end{proof}

\begin{lemma}\label{l:deviations-exp}
Let $(Y_{i})_{i \geq 1}$ be a sequence of non-negative integer-valued random variables on a probability space $(O,\H,\T)$, and $(\H_{i})_{i \geq 0}$ an non-decreasing sequence of sub-$\sigma-$algebras of $\H$
 such that $\H_{0}= \{ \emptyset, O \}$. Assume that the following properties hold:
  \begin{itemize}
 \item for all $i \geq 1$, $Y_{i}$ is measurable with respect to $\H_{i}$;
\item there exists $0 < c_{1}, c_{2} < +\infty$ such that for all $i \geq 1$ and $k \geq 0$, $\T(Y_{i} \geq t | \H_{i-1}) \leq c_{1} \exp(-c_{2}k)$.
 \end{itemize}

 Then there exists $c_{3}$ depending only on $c_{1},c_{2}$ such that, for all $t > c_{3}$, there exist  $0 < c_{5}, c_{6} < +\infty$ such that, for all $1 \leq n \leq m$,
 $\T( Y_{1}+\cdots+Y_{n} \geq mt) \leq   c_{5} \exp(-c_{6} m)$.
\end{lemma}

\begin{proof}

For $0 < \lambda < c_{2}$, one has a.s.
\begin{eqnarray*}\E_{\T}( \exp(\lambda Y_{i}) | \H_{i-1}  ) &\leq&  1 + \sum_{k=1}^{+\infty}   (e^{\lambda k} - e^{\lambda (k-1)} )  \T(Y_{i} \geq k | \H_{i-1}) \\
&\leq&  1 + c_{1} (1-e^{-\lambda})  \frac{e^{\lambda-c_{2}}}{1-e^{\lambda-c_{2}}}.\end{eqnarray*}
Letting $j(\lambda):=c_{1} (1-e^{-\lambda})  \frac{e^{\lambda-c_{2}}}{1-e^{\lambda-c_{2}}}$,
we deduce that $$\E_{\T}( \exp(\lambda (Y_{1}+\cdots + Y_{m})) ) \leq      \left( 1 + j(\lambda) \right)^{m}.$$

Then, by Markov's inequality,
 $$\T( Y_{1}+\cdots+Y_{m} \geq mt)  \leq \exp(-m \lambda t) \E_{\T} (\exp( \lambda   (Y_{1}+\cdots + Y_{m}   )),$$
 so that
 \begin{equation}\label{e:encore-cramer}\T( Y_{1}+\cdots+Y_{m} \geq mt)   \leq 
  \exp \left[ -m  \left( \lambda t  + \log (  1 + j(\lambda)) \right) \right].\end{equation}
As $\lambda$ goes to zero, we see that $j(\lambda)=c_{3} \lambda + o(\lambda)$, 
with $c_{3}:= \frac{e^{-c_{2}}}{1-e^{-c_{2}}}$. Choosing $\lambda$ small enough in~(\ref{e:encore-cramer}) yields the result when $n=m$.
For $n \leq m$, observe that by assumption $Y_{1}+\cdots+Y_{n} \leq Y_{1}+\ldots+Y_{m}$.




\end{proof}

\begin{lemma}\label{l:estimation-finale-unif}
For $L\geq L_{0}$, there exists $0 < C_{74}, C_{75} < +\infty$ such that, for all $0 \leq \epsilon \leq \epsilon_{0}$,  
and all $k \geq 1$, 
\begin{itemize}

\item[a)] $\Q^{\epsilon,\theta}_{\I_{0}}(\hat{r}_{S_{k}} >   k C_{75} +  u , \, K > k ) \leq  C_{74} k^{2} u^{-4}$;

\item[b)] $ \Q^{\epsilon,\theta}_{a \delta_{0}}(  \hat{r}_{S_{k}} >   k C_{75} +  u,  \,  U=+\infty, \, K > k)  \leq  C_{74} k^{2} u^{-4}.$

\end{itemize}

\end{lemma}

\begin{proof}

Fix $L \geq L_{0}$. Observe that, for any $k \geq 1$, on $\{K > k\}$, 
\begin{equation}\label{decoupons-r}\hat{r}_{S_{k}} = \hat{r}_{0} + (\hat{r}_{S_{1}}-\hat{r}_{0}) + \sum_{j=1}^{k-1}   \left( \hat{r}_{S_{j+1}}-\hat{r}_{D_{j}}+\hat{r}_{D_{j}}-\hat{r}_{S_{j}}  \right) \un(  K \geq j ).\end{equation}
Observe that, for $w=w_{\hat{r}_{S_{j}}}$ with $1  \leq j \leq K$, denoting $w=(F,r,A)$, the three conditions  
$w \in \L'_{\theta}$, $\phi_{r-L}(w) \leq p$, and $m_{r-\floor{L^{1/4}}, r}(w) \geq a \floor{L^{1/4}}/2$ are satisfied.  
As a consequence, by Lemma~\ref{l:lemme11-unif} and the strong Markov property, for all $1 \leq j \leq k-1$, and all $t>0$, a.s.
 $\Q^{\epsilon,\theta}_{\I_{0}}(\hat{r}_{D_{j}} - \hat{r}_{S_{j}}  > t , \, K \geq j | \F_{S_{j}}  ) \leq  C_{14} \left(t^{-M'} +   L \exp(-C_{15} t) \right)$.
  Now letting, for $j \geq 1$, $Y_{j}:=\left(\hat{r}_{D_{j}}-\hat{r}_{S_{j}}  \right) \un(  K \geq j )$, and $\H_{j}:=\F_{S_{j+1}}$, we see that the assumptions of Lemma~\ref{l:deviations-moments}
are satisfied with $q=2$, since $M' = a+8$.  

Thanks to the above observation on $w=w_{\hat{r}_{S_{j}}}$, and to the fact that, on $\{  K \geq j   \}$, 
$\hat{r}_{S_{j+1}}-\hat{r}_{D_{j}} = L J_{\hat{r}_{D_{j}}}$, 
we see that, by Lemma~\ref{l:lemme17-unif} b) and the strong Markov property, for all $1 \leq j \leq k-1$, and all $t>0$, a.s.
$\Q^{\epsilon,\theta}_{\I_{0}}(\hat{r}_{S_{j+1}} - \hat{r}_{D_{j}}  > t , \, K \geq j| \F_{S_{j}}  ) \leq C_{26} (\floor{L^{-1}t})^{3-M'}$.
Similarly, thanks to Lemma~\ref{l:lemme17-unif} a), one also has that, for all $t>0$, a.s.
 $\Q^{\epsilon,\theta}_{\I_{0}}(\hat{r}_{S_{1}} - \hat{r}_{0}  > t , \, K \geq j| \F_{S_{j}}  ) \leq C_{26} (\floor{L^{-1}t})^{3-M'}$.
 
 Now letting $Y_{1}:=\hat{r}_{S_{1}} - \hat{r}_{0}$, and, for  $j \geq 2$, $Y_{j}:=\left(\hat{r}_{S_{j}}-\hat{r}_{D_{j-1}}  \right) \un(  K \geq j )$, 
 and $\H_{j}:=\F_{S_{j}}$, we see that the assumptions of Lemma~\ref{l:deviations-moments}
are again satisfied with $q=2$.  
  Applying Lemma~\ref{l:deviations-moments}, we deduce the existence of two constants $C_{751}, C_{741}$ not depending on $\epsilon$ such that for all $k \geq 1$ and $u>0$,
 $$\Q^{\epsilon,\theta}_{\I_{0}} \left(     \sum_{j=1}^{k-1}   \left(\hat{r}_{D_{j}}-\hat{r}_{S_{j}}  \right) \un(  K \geq j )    >   k C_{751} +  u , \, K > k \right) \leq  C_{741} k^{2} u^{-4},$$
 and  
 $$\Q^{\epsilon,\theta}_{\I_{0}} \left( \hat{r}_{S_{1}} - \hat{r}_{0} +   \sum_{j=1}^{k-1} \left(  \hat{r}_{S_{j+1}}-\hat{r}_{D_{j}}  \right) \un(  K \geq j )    >   k C_{751} +  u , \, K > k \right) \leq  C_{741} k^{2} u^{-4}.$$
  Part a) of the lemma then follows from the two above inequalities, (\ref{decoupons-r}), and the union bound.
 
 To prove part b), we note that, for all $k \geq 1$,   on $\{K > k, \, U=+\infty\}$, 
 \begin{equation}\label{decoupons-r-2}\hat{r}_{S_{k}} = \hat{r}_{0} + (\hat{r}_{S_{1}}-\hat{r}_{0}) \un(  U = +\infty      ) 
 + \sum_{j=1}^{k-1}   \left( \hat{r}_{S_{j+1}}-\hat{r}_{D_{j}}+\hat{r}_{D_{j}}-\hat{r}_{S_{j}}  \right) \un(  K \geq j ).\end{equation}
We can use the same argument as in the proof of part a) to deal with $ \sum_{j=1}^{k-1}   \left(\hat{r}_{D_{j}}-\hat{r}_{S_{j}}  \right)\un(  K \geq j ) $ and 
$\sum_{j=1}^{k-1} \left(  \hat{r}_{S_{j+1}}-\hat{r}_{D_{j}}  \right) \un(  K \geq j )   $.
 To deal with the remaining term $(\hat{r}_{S_{1}}-\hat{r}_{0}) \un(  U = +\infty )$, observe that $\hat{r}_{S_{1}}-\hat{r}_{0} = L J_{0}$, and apply Lemma~\ref{l:lemme17-unif} c).
\end{proof}

\begin{lemma}\label{l:estimation-finale-exp}

For all $0 \leq \epsilon \leq \epsilon _{0}$, and $L \geq L_{1}(\epsilon)$, there exist $0<C_{97}(\epsilon),C_{98}(\epsilon), C_{99}(\epsilon)<+\infty$ such that, 
for all $k \leq m$,
\begin{itemize}
\item[a)] $\Q^{\epsilon,\theta}_{\I_{0}}(\hat{r}_{S_{k}} >   m C_{97}(\epsilon) , \, K>k ) \leq  C_{98}(\epsilon) \exp(-C_{99}(\epsilon) m);$
and
\item[b)] $\Q^{\epsilon,\theta}_{a \delta_{0}}(  \hat{r}_{S_{k}} >   m C_{97}(\epsilon) ,  \,  U=+\infty, \, K>k)  \leq    C_{98}(\epsilon) \exp(-C_{99}(\epsilon) m).$
\end{itemize}
\end{lemma}

\begin{proof}
Adapt the proof of Lemma~\ref{l:estimation-finale-unif}, using Lemma~\ref{l:deviations-exp} instead of Lemma~\ref{l:deviations-moments}, and Lemma~\ref{l:lemme17-exp} instead
of  Lemma~\ref{l:lemme17-unif}.
\end{proof}

\begin{prop}\label{p:moments-unif}
For all $L \geq L_{0}$, there exists $0<C_{29}<+\infty$ not depending on $\epsilon$ such that, for all $0 \leq \epsilon \leq \epsilon_{0}$,
\begin{itemize}
\item[a)] $\E^{\epsilon,\theta}_{\I_{0}}(  \kappa^{2} ) \leq C_{29}, \, \E^{\epsilon,\theta}_{\I_{0}}(  (\hat{r}_{\kappa})^{2} ) \leq C_{29};$
\item[b)] $\E^{\epsilon,\theta}_{a \delta_{0}}(  \kappa^{2} | U=+\infty) \leq C_{29}, \, \E^{\epsilon,\theta}_{a \delta_{0}}(  (\hat{r}_{\kappa})^{2}| U=+\infty) \leq C_{29}.$
\end{itemize}
\end{prop}

\begin{proof}[Proof of Proposition~\ref{p:moments-unif}]

Observe that, for any integer $\ell \geq 1$,
$$\{ \kappa > t \} \subset   \{  K> \ell   \}  \cup \bigcup_{k=1}^{\ell} \{ K=k, \,  S_{k} > t \},$$
whence
\begin{equation}\label{e:decomp-kappa}\{ \kappa > t \} \subset   \{  K> \ell   \}
\cup \bigcup_{k=1}^{\ell} \{ K=k, \,  \hat{r}_{S_{k}} \geq \floor{\alpha_{1} t} \} \cup  \{  \cup_{s \geq t}  \hat{r}_{s} < \floor{\alpha_{1}s}  \}  .\end{equation}

By the union bound, 
\begin{equation}\label{e:une-reunion-de-plus}\Q^{\epsilon,\theta}_{\I_{0}}(\kappa > t  )  \leq \Q^{\epsilon,\theta}_{\I_{0}}(K>\ell)  +
 \sum_{k=1}^{\ell} \Q^{\epsilon,\theta}_{\I_{0}}(\hat{r}_{S_{k}} \geq \floor{\alpha_{1} t}, K=k)+
   \Q^{\epsilon,\theta}_{\I_{0}}( \cup_{s \geq t}  \hat{r}_{s} < \floor{\alpha_{1}s}).\end{equation}

Now remember $\delta_{3}$ defined in Corollary~\ref{c:corollaire2-bis}  and let
$\ell:=  -4 \log \left( (1-\delta_{3})^{-1}  \ceil{t} \right)$.
By~(\ref{e:L}), $\phi_{r-L}(\I_{0}) \leq p$ so that  $\Q^{\epsilon,\theta}_{\I_{0}}(D<+\infty) \leq 1 - \delta_{3}$. 
Moreover, for all $j \geq 1$, on $K \geq j$, $\phi_{r-L}(w_{\hat{r}_{S_{j}}})$, so that, by the strong Markov property, we have 
a.s. $\Q^{\epsilon,\theta}_{\I_{0}}(D<+\infty | \F_{S_{j}}) \leq 1 - \delta_{3}$. 
We deduce that  
\begin{equation}\label{e:et-un-morceau-un}\Q^{\epsilon,\theta}_{\I_{0}}(K>\ell) \leq (1-\delta_{3})^{\ell} \leq t^{-4}.\end{equation}

Now observe that, for large enough $t$ (not depending on $\epsilon$),
$ \floor{\alpha_{1}t} \geq \ell C_{75} + \alpha_{1}t/2$.
Using  Lemma~\ref{l:estimation-finale-unif} a), we deduce that, for all $1 \leq k\leq \ell$,
\begin{equation}\label{e:et-un-deuxieme-morceau}\Q^{\epsilon,\theta}_{\I_{0}}(\hat{r}_{S_{k}} >   \floor{\alpha_{1} t} , \, K>k ) \leq  C_{74} k^{2} (\alpha_{1}t/2)^{-4}.\end{equation}

Finally, by Lemma~\ref{l:lemme7-addendum-unif},
\begin{equation}\label{e:et-un-troisieme-et-dernier}\Q^{\epsilon,\theta}_{\I_{0}}\left[    \cup_{s \geq t}     \hat{r}_{s} < \floor{\alpha_{1}s} \right] \leq C_{45} t^{-M'}.\end{equation}

Plugging~(\ref{e:et-un-morceau-un}),~(\ref{e:et-un-deuxieme-morceau}) and~(\ref{e:et-un-troisieme-et-dernier}) into~(\ref{e:une-reunion-de-plus}),
 we deduce the conclusion of part a) regarding $\E^{\epsilon,\theta}_{\I_{0}}(  \kappa^{2} )$. 
 The conclusion for  $\E^{\epsilon,\theta}_{\I_{0}}(  (\hat{r}_{\kappa})^{2} )$ follows by an application of Lemma~\ref{l:lemme10}.

As for part b), observe that the estimate in~(\ref{e:et-un-morceau-un}) is still valid when $\I_{0}$ is replaced by $a \delta_{0}$.
On the other hand, the estimate obtained in~(\ref{e:et-un-deuxieme-morceau}) follows from  Lemma~\ref{l:estimation-finale-unif} b). 
Then, by definition, the event $U=+\infty$ rules out the event $\cup_{s \geq t}     \hat{r}_{s} < \floor{\alpha_{1}s}$. 
Part b) is then proved exactly as part a), noting that,  $Q_{a \delta_{0}}(U=+\infty) \geq 1-\delta_{2}$.

\end{proof}

\begin{prop}\label{p:moments-exp}
For all $0<\epsilon \leq \epsilon_{0}$, and $L \geq L_{1}(\epsilon)$, there exists $0<C_{30}(\epsilon), C_{31}(\epsilon)<+\infty$ such that
\begin{itemize}
\item[a)] $\E^{\epsilon,\theta}_{\I_{0}}(  \exp(  -C_{30}(\epsilon) \kappa    ) ) \leq C_{31}(\epsilon), \, \E^{\epsilon,\theta}_{\I_{0}}( \exp(  -C_{30}(\epsilon) \hat{r}_{\kappa} ) \leq C_{31}(\epsilon);$
\item[b)] $\E^{\epsilon,\theta}_{a \delta_{0}}(  \exp(  -C_{30}(\epsilon) \kappa          | U=+\infty) \leq C_{31}(\epsilon), \, \E^{\epsilon,\theta}_{a \delta_{0}}(   \exp(  -C_{30}(\epsilon) \hat{r}_{\kappa}| U=+\infty) \leq C_{31}(\epsilon).$
\end{itemize}
\end{prop}

\begin{proof}[Proof of Proposition~\ref{p:moments-exp}]

The proof is very similar to the proof of Proposition~\ref{p:moments-unif}, but 
this time, we use  $\ell:= \floor{ (1/2) C_{97}(\epsilon)^{-1} \alpha_{1} t  }$, so that the r.h.s. of~(\ref{e:et-un-morceau-un}) now decays exponentially as $t \to +\infty$.

We then use Lemma~\ref{l:estimation-finale-exp} instead of Lemma~\ref{l:estimation-finale-unif}, noting that, 
for large enough $t$, $\floor{\alpha_{1}t} \geq \ell C_{97}(\epsilon)$.
Finally, we use Lemma~\ref{l:lemme7-addendum-exp} instead of Lemma~\ref{l:lemme7-addendum-unif}, and the 
conclusion follows as in the proof of  Proposition~\ref{p:moments-unif}.

\end{proof}
\bibliographystyle{plain}
\bibliography{large-deviations}

\end{document}